\documentclass[reqno]{amsart}
\usepackage{mathtools,amssymb,amsthm,mathrsfs,graphicx}
\usepackage{empheq}

\usepackage{geometry}
\usepackage[numbers,sort&compress]{natbib}

\usepackage[draft=false,hidelinks]{hyperref}
\hypersetup{
  colorlinks   = true, 
  urlcolor     = blue, 
  linkcolor    = red, 
  citecolor   = blue 
}
\usepackage{esint}

\numberwithin{equation}{section}

\theoremstyle{plain}
\newtheorem{thm}{Theorem}[section]
\newtheorem{lem}[thm]{Lemma}
\newtheorem*{lem*}{Lemma}
\newtheorem*{cor}{Corollary}
\newtheorem{proposition}[thm]{Proposition}

\theoremstyle{definition}
\newtheorem{defn}[thm]{Definition}

\theoremstyle{remark}
\newtheorem{remark}[thm]{Remark}

\allowdisplaybreaks

\newcommand{\mynegspace}{\hspace{-0.12em}}
\newcommand{\lvvvert}{\rvert\mynegspace\rvert\mynegspace\rvert}
\newcommand{\rvvvert}{\rvert\mynegspace\rvert\mynegspace\rvert}
\DeclarePairedDelimiter{\vvvert}{\lvvvert}{\rvvvert}

\newcommand{\R}{\mathbb{R}}
\newcommand{\C}{\mathbb{C}}

\DeclareMathOperator*{\esssup}{ess\,sup}

\begin{document}

\title[Elliptic and parabolic equations with degenerate weights]{Elliptic and parabolic equations with rough boundary data in Sobolev spaces with degenerate weights}

\author[B. Bekmaganbetov]{Bekarys Bekmaganbetov}
\address{Division of Applied Mathematics, Brown University, 182 George Street, Providence, RI 02912, USA}
\email{bekarys\_bekmaganbetov@brown.edu}

\author[H. Dong]{Hongjie Dong}
\address{Division of Applied Mathematics, Brown University, 182 George Street, Providence, RI 02912, USA}
\email{hongjie\_dong@brown.edu}
\thanks{H. Dong and B. Bekmaganbetov were partially supported by the NSF under agreements DMS-2055244 and DMS-2350129.}

\keywords{Parabolic equations with inhomogeneous Dirichlet boundary conditions, divergence form, singular-degenerate coefficients, traces, weighted Sobolev spaces}
\subjclass[2020]{35K65; 35K67; 35K20; 35D30; 46E35}

\begin{abstract}
    We investigate the inhomogeneous boundary value problem for elliptic and parabolic equations in divergence form in the half space $\{x_d > 0\}$, where the coefficients are measurable, singular or degenerate, and depend only on $x_d$. The boundary data are considered in Besov spaces of distributions with negative orders of differentiability in the range $(-1,0]$. The solution spaces are weighted Sobolev spaces with power weights that decay rapidly near the boundary, and are outside the Muckenhoupt $A_p$ class. Sobolev spaces with such weights contain functions that are very singular near the boundary and do not possess a trace on the boundary. Consequently, solutions may not exist for arbitrarily prescribed boundary data and right-hand sides of the equations. We establish a natural structural condition on the right-hand sides of the equations under which the boundary value problem is well-posed. 
\end{abstract}
\maketitle

\section{Introduction}

There is significant interest in studying elliptic and parabolic equations with singular or degenerate coefficients due to their broad range of applications in both pure and applied mathematics. These applications span fields such as probability theory, geometric PDEs, porous media, mathematical finance, and mathematical biology; see \cite{DongPhan2021TAMS} for references to related problems. Additionally, such equations are closely linked to extension operators for the fractional Laplace and fractional heat equations; see, for example, \cite{CaffarelliSilvestre2007CPDE, StingaTorrea2017SIAM}.

In this paper, we address the inhomogeneous Dirichlet problem for elliptic and parabolic equations in divergence form in the upper-half space $\R^d_+$, whose coefficients contain the factor $x_d^{\alpha}$ with $\alpha \in (-\infty,1)$, which is singular or degenerate near the boundary $\partial \R^d_+$ when $\alpha < 0$ or $\alpha > 0$, respectively. We establish the solvability of such equations in weighted Sobolev spaces with degenerate power-type weights $x_d^{\beta}$, where $\beta$ is a parameter. The boundary data in this context is taken from spaces with \emph{negative} orders of regularity, which are considered ``rough.'' To highlight our main results for elliptic equations, we present the following theorem, which serves as a special case of the main theorems in Section \ref{sec:mainres}.

Throughout the paper, let $d \geq 2$ be fixed. We denote the upper half space by $\mathbb{R}^d_+ = \{x=(x',x_d) \colon x' \in \mathbb{R}^{d-1}, x_d > 0 \}$, where $x'=(x_1,\dots,x_{d-1})$.  By $L_p(\R^d_+,w dx)$, we denote the weighted Lebesgue space with a given weight $w$. The weak partial derivative in the $x_i$ direction is denoted by $D_i = \partial_{x_i}$ for $i = 1,\dots,d$, and $Du = (D_1 u, \dots, D_d u)$ represents the spatial gradient of a function $u$. Lastly, $B^s_{p,q}(\R^{d-1})$ denotes the inhomogeneous Besov space on $\R^{d-1}$.

\begin{thm} \label{thm:simplethmell}
    Let $\alpha \in (-\infty,1)$, $p \in (1,\infty)$, $\beta \in [p-1,2p-1)$, and $\lambda > 0$. Denote $\Vert \cdot \Vert = \Vert \cdot \Vert_{L_{p}(\R^d_+,x_d^{\beta}dx)}$ and  $N=N(d,\alpha,p,\beta)>0$.
    \begin{itemize}
        \item[(i)] For any $F_1,F_2,\dots,F_{d-1}, f \in L_p(\R^d_+, x_d^{\beta}dx)$ and $U \in B^{1-\frac{\beta+1}{p}}_{p,p}(\R^{d-1})$, there exists a unique weak solution $u=u(x) \in W^1_p(\R^d_+, x_d^{\beta}dx)$ to the equation
        \begin{equation} \label{eq:simplebvpell}
            \begin{cases}
                -\sum_{i=1}^d D_i(x_d^{\alpha}D_i u) + \lambda x_d^{\alpha} u = -\sum_{i=1}^{d-1}D_i (x_d^{\alpha}F_i) + x_d^{\alpha}f & \text{in} \; \R^d_+, \\
                u(x',0) = U(x') & \text{on} \; \R^{d-1}.
            \end{cases}
        \end{equation}
        Moreover, $u$ satisfies
        \begin{equation*}
            \Vert Du \Vert + \sqrt{\lambda} \Vert u \Vert \leq N\Vert F \Vert + \frac{N}{\sqrt{\lambda}}\Vert f \Vert + \Tilde{N}\Vert U \Vert_{B^{1-\frac{\beta+1}{p}}_{p,p}(\R^{d-1})},
        \end{equation*}
        where $\Tilde{N} = \Tilde{N}(d,\alpha,p,\beta,\lambda)>0$.
        \item[(ii)] Suppose that $u \in W^1_p(\R^d_+,x_d^{\beta}dx)$ satisfies {\rm (}weakly{\rm )}
        \begin{equation*}
            D_d(x_d^{\alpha} D_d u) = \sum_{i=1}^{d-1}D_i(x_d^{\alpha}F_i) + x_d^{\alpha}f \qquad \text{in} \; \R^d_+
        \end{equation*}
        for some $F_1,\dots,F_{d-1},f \in L_p(\R^d_+,x_d^{\beta}dx)$. Then $u$ has a boundary trace $U(x')=u(x',0) \in B^{1-\frac{\beta+1}{p}}_{p,p}(\R^{d-1})$ in the following sense:
        \begin{equation*}
            u(\cdot,x_d) \to U \quad \text{strongly in $B^{1-\frac{\beta+1}{p}-\varepsilon}_{p,p}(\R^{d-1})$ as $x_d \to 0^+$,}
        \end{equation*}
        where $\varepsilon=0$ if $\beta > p-1$, and $\varepsilon>0$ is an arbitrary positive real if $\beta=p-1$. Moreover,
        \begin{equation*}
            \Vert U \Vert_{B^{1-\frac{\beta+1}{p}}_{p,p}(\R^{d-1})} \leq N(\Vert u \Vert + \Vert Du \Vert + \Vert F \Vert + \Vert f \Vert).
        \end{equation*}
        \item[(iii)] Let $\beta \in (p-1,2p-1)$. For any $F_1,\dots,F_d,f \in L_p(\R^d_+,x_d^{\beta}dx)$, there exists a weak solution $u$ to the equation
        \begin{equation*} 
            -\sum_{i=1}^d D_i(x_d^{\alpha}D_i u) + \lambda x_d^{\alpha} u = -\sum_{i=1}^{d}D_i (x_d^{\alpha}F_i) + x_d^{\alpha}f \qquad \text{in} \; \R^d_+,
        \end{equation*}
        that satisfies
        \begin{equation*}
            \Vert Du \Vert + \sqrt{\lambda} \Vert u \Vert \leq N\Vert F \Vert + \frac{N}{\sqrt{\lambda}}\Vert f \Vert.
        \end{equation*}
    \end{itemize}
\end{thm}

In the case of parabolic equations, we consider the boundary value problem
\begin{equation} \label{eq:simplebvppar}
    \begin{cases}
        x_d^{\alpha}u_t-\sum_{i=1}^d D_i(x_d^{\alpha}D_i u) + \lambda x_d^{\alpha} u = -\sum_{i=1}^{d-1}D_i (x_d^{\alpha}F_i) + x_d^{\alpha}f & \text{in} \; \R^{d+1}_+, \\
        u(t,x',0) = U(t,x') & \text{on} \; \R^{d},
    \end{cases}
\end{equation}
for which we obtain similar results. Here
\[
\R^{d+1}_+=\{(t,x): t \in \R, x \in \R^d_+\}, \quad \R^d = \partial \R^{d+1}_+ = \{(t,x') \in \R^d\}.
\]
On the level of function spaces, in the parabolic case we work with
\[
u,Du,F,f \in L_p(\R^{d+1}_+,x_d^{\beta}dxdt),
\]and
$U(t,x')=u(t,x',0)$ now belongs to the time-space anisotropic Besov space:
\[
U \in B^{\frac{1}{2}-\frac{\beta+1}{2p},1-\frac{\beta+1}{p}}_{p,p}(\R_t \times \R^{d-1}_{x'}).
\]
Parts (i), (ii), and (iii) of Theorem \ref{thm:simplethmell} represent special cases of respectively Theorems \ref{thm:solvell}, \ref{thm:traceell}, and \ref{thm:existsolell}, in which more general equations with measurable coefficients depending on $x_d$ are considered. The corresponding results for parabolic equations are presented in Theorems \ref{thm:solvpar}, \ref{thm:tracepar}, and \ref{thm:existsolpar}.

Theorem \ref{thm:simplethmell}, being a special case, already captures several key features of our main results. One such feature is that the boundary data $U(x') = u(x',0)$ are taken from Besov spaces $B^{1-\frac{\beta+1}{p}}_{p,p}(\R^{d-1})$. We focus on weight powers $\beta \in [p-1,2p-1)$, where the order of regularity $1-\frac{\beta+1}{p}$ lies in $(-1,0]$. This leads to the treatment of ``rough" boundary data $U(x')$, which are in general distributions rather than functions. Using spaces with weights degenerating near the boundary is a natural approach to treat rough boundary data. Notably, the weight $w_{\beta}(x)=x_d^\beta$ belongs to the Muckenhoupt class $A_p(\R^d_+)$ if and only if $\beta \in (-1,p-1)$, which places our work in a non-$A_p$ framework. Solvability results for equations in non-divergence form in the non-$A_p$ range $\beta \in [p-1,2p-1)$ are known from \cite{KozlovNazarovMathNach2009, LindemulderVeraar2020JDE, DongPhan2023JFA}. Extending these results to equations in divergence form was one of the main motivations of our work. 

Our results represent a partial extension of \cite{DongPhan2021TAMS} to the non-$A_p$ setting, where divergence form equations with the zero boundary condition and $\beta \in (\alpha p-1, p-1)$ (along with more general families of weights and coefficients) were investigated. Another distinctive aspect of our results is that in Parts (i) and (ii) of Theorem \ref{thm:simplethmell}, when the boundary trace $U(x')=u(x',0)$ exists, the right-hand side does not contain the term $D_d(x_d^{\alpha}F_d)$, which is typically allowed to be non-zero in the theory of equations in divergence form. However, according to Part (iii) of Theorem \ref{thm:simplethmell}, the existence of solutions with $F_d \neq 0$ still holds, but boundary values can no longer be prescribed. Further explanation of this can be found in the discussion following (\ref{dirichletparabolic}).

The study of regularity theory for equations with singular-degenerate coefficients attracted much attention in recent years. Equations in divergence and non-divergence form with the zero Dirichlet boundary condition were considered in \cite{DongPhan2021TAMS} and \cite{DongPhan2023JFA}, respectively. Additionally, equations with the conormal boundary condition were considered in \cite{DongPhanIndiana2023}. In these papers, equations with partially BMO type coefficients were studied, and solvability and estimates of solutions in weighted Sobolev spaces were obtained. Recent advances in boundary Schauder-type estimates for parabolic equations with conormal boundary conditions were made in \cite{audrito2024schauder, audrito2024higher}, while similar Schauder estimates for elliptic equations can be found in \cite{SireTerraciniVitaCPDE2021, SireTerraciniVitaMathEng2021, DongJeonVita2023arxiv}. Another general types of equations with singular-degenerate coefficients were studied in \cite{metafune2023JDE, metafune2024regularity}. There, the authors allow for different weights in the tangential and normal derivative components of the second order operators, assuming that the uniformly elliptic part has constant coefficients. For a comprehensive study of anisotropic weighted Sobolev spaces, which arise naturally in the analysis of such equations, see \cite{MetafuneNegroSpina2023TokyoJMath}.

Numerous significant results have been obtained in various directions for uniformly elliptic and parabolic equations (i.e., when $\alpha = 0$ in our notation). In \cite{KozlovNazarovMathNach2009}, the authors derived weighted $L_p$ estimates with power weights for parabolic equations in non-divergence form with the zero boundary condition. Their approach relied on Green's function estimates for equations with discontinuous coefficients depending only on time. Weighted $L_p$ estimates can also be obtained by a different approach based on $H^{\infty}$--calculus, which also allows treating Banach space-valued solutions. We refer to \cite{LindemulderVeraar2020JDE} and the more recent work \cite{LLRVarxiv2024} for an extensive study of various inhomogeneous (initial-)boundary value problems for the Laplace (and heat) equation in weighted spaces using $H^{\infty}$-calculus. In \cite{KimKHKimWoo2024}, the authors explored parabolic equations with solution spaces different from ours, which are based on weighted Sobolev spaces $H^{\gamma}_{p,\theta}$ introduced by Krylov in \cite{KrylovCPDE1999}. Another recent work, \cite{JungKim2024arxiv}, focused on parabolic equations in divergence form in an unweighted setting, where the right-hand side involves a half-derivative in time term, $D_t^{1/2}h$. This term naturally emerges when transforming an equation with an inhomogeneous boundary condition into an equivalent one with zero boundary condition by subtracting an appropriate function. However, we do not employ this approach in this paper. For a comprehensive treatment of boundary value problems for elliptic equations in half spaces using layer potential methods, we also refer to the monograph \cite{BartonMayboroda}.

To give a formal description of our main results for general equations, we introduce some notation. Assume that the coefficients $\{a_0(x_d), a_{ij}(x_d), c(x_d)\}_{i,j = 1}^d:\mathbb{R}_+ \to \mathbb{R}$ are given measurable functions of $x_d$ only, and satisfy the uniform ellipticity and boundedness condition: there exists $\kappa \in (0,1]$ such that
\begin{equation} \label{eq:ellipticitycond}
\kappa |\xi|^2 \leq \sum_{i,j=1}^d a_{ij}(x_d)\xi_i \xi_j, \quad |a_{ij}(x_d)| \leq \kappa^{-1}, \quad \kappa \leq a_0(x_d), c(x_d) \leq \kappa^{-1}
\end{equation}
for all $\xi=(\xi_1,\xi_2,\dots,\xi_d) \in \mathbb{R}^d$ and $x_d \in (0,\infty)$. We study the following elliptic equation in divergence form with singular or degenerate coefficients subject to the inhomogeneous Dirichlet boundary condition:
\begin{empheq}[left={\empheqlbrace}]{align}
	\label{equationnonhomweighted}
    -\sum_{i,j=1}^d D_i(x_d^\alpha a_{ij}(x_d) D_j u)+\lambda x_d^\alpha c(x_d) u &= -\sum_{i=1}^d D_i(x_d^\alpha F_i) + x_d^\alpha f && \text{in} \; \mathbb{R}^d_+, \\
    \label{dirichlet}
    u(x',0) &= U(x') && \text{on} \; \mathbb{R}^{d-1},
\end{empheq}
and the corresponding parabolic equation
\begin{empheq}[left={\empheqlbrace}]{align}
	\label{equationparabolicdivnonhom}
	x_d^{\alpha}a_0(x_d) u_t- \sum_{i,j=1}^d D_i(x_d^\alpha a_{ij}(x_d) D_j u)+\lambda x_d^\alpha c(x_d) u & = -\sum_{i=1}^d D_i(x_d^\alpha F_i) + x_d^\alpha f &&  \text{in} \; \R^{d+1}_+ \\
	\label{dirichletparabolic}
	u(t,x',0) & = U(t,x') && \text{on} \; \R^d.
\end{empheq}

Our main results concerning parabolic equations can be summarized as follows. In Theorem \ref{thm:solvpar}, we establish the existence and uniqueness of solutions to the parabolic problem (\ref{equationparabolicdivnonhom})-(\ref{dirichletparabolic}) under the condition that $F_d = 0$. In Theorem \ref{thm:tracepar}, we present a trace theorem for solution spaces of the parabolic equation (\ref{equationparabolicdivnonhom}) with $F_d=0$.
Furthermore, in Theorem \ref{thm:existsolpar}, we prove the existence of solutions to equation (\ref{equationparabolicdivnonhom}) for general $F_d \neq 0$. The corresponding elliptic results for equations (\ref{equationnonhomweighted})-(\ref{dirichlet}) are covered by Theorems \ref{thm:solvell}, \ref{thm:traceell}, and \ref{thm:existsolell}.

To clarify the condition $F_d=0$, we recall the characterization of trace spaces for weighted Sobolev spaces $W^1_p(\R^d_+,x_d^{\beta}dx)$. Specifically, when $\beta \in (-1,p-1)$, the trace space of $W^1_p(\R^d_+,x_d^{\beta}dx)$ is $B^{1-\frac{\beta+1}{p}}_{p,p}(\R^{d-1})$; see \cite[2.9.2]{Triebel}. On the other hand, when $\beta \leq -1$, any function $u \in W^1_p(\R^d_+,x_d^{\beta}dx)$ automatically satisfies $u(x',0)=0$; see \cite[Theorem 1.3]{Grisvard}. If $\beta \geq p-1$, functions in $W^1_p(\R^d_+,x_d^{\beta}dx)$ may blow-up near the boundary and hence may not admit a well-defined boundary trace. This can be illustrated by the following simple example. Take
\begin{equation} \label{eq:exnotrace}
    u(x',x_d) =
    \begin{cases}
        \varphi(x') \log\log\frac{1}{x_d}, & \text{if} \; \beta = p-1, \\
        \varphi(x') x_d^{-\delta}, & \text{if} \; \beta > p-1.
    \end{cases}
\end{equation}
Here $\varphi$ is an arbitrary fixed non-zero function in $C_0^\infty(\mathbb{R}^{d-1})$ and $\delta$ is any number in $(0,\frac{\beta+1}{p}-1)$ (when $\beta>p-1$). Then $u \in W^1_{p}(\mathbb{R}^{d-1}\times (0,\frac{1}{2}), x_d^{\beta}dx)$ and $u(\cdot,x_d)$ does not have a limit in $\mathcal{D}'(\mathbb{R}^{d-1})$, i.e., in the sense of distributions, as $x_d \to 0^+$. This observation leads to the following counterexample. 

\begin{proposition} \label{prop:counterex}
Let $p \in (1,\infty)$ and $\beta \in [p-1, 2p-1)$. There exists $\bar{F}_d \in L_p(\R^d_+,x_d^{\beta}dx) \cap C^{\infty}(\R^d_+)$ such that the following assertion holds. Let $u \in W^1_p(\R^d_+,x_d^{\beta}dx)$ be a weak solution to
\[
-\Delta u + u = -\sum_{i=1}^{d-1}D_i F_i - D_d \bar{F}_d + f \qquad \text{in \, $\R^d_+$}
\]
for some $F_1,\dots,F_{d-1},f \in L_p(\R^d_+,x_d^{\beta}dx)$. Then $u(\cdot,x_d)$ does not have a limit in $\mathcal{D}'(\R^{d-1})$ as $x_d \to 0^+$.
\end{proposition}

Indeed, take $\bar{F}_d = D_d \bar{u}$, where $\bar{u} \in W^1_p(\R^d_+,x_d^{\beta}dx)$ and $\bar{u}(\cdot,x_d)$ does not have a limit in $\mathcal{D}'(\R^{d-1})$ as $x_d \to 0^+$, as given in (\ref{eq:exnotrace}). By the construction, one can verify that the difference $u-\bar{u}$ satisfies
\[
D^2_d(u-\bar{u})=\sum_{i=1}^{d-1} D_i \Tilde{F}_i + \Tilde{f}
\]
for some $\Tilde{F}_i, \Tilde{f} \in L_p(\R^d_+,x_d^{\beta}dx)$. Then by Part (ii) of Theorem \ref{thm:simplethmell} with $\alpha = 0$, $(u-\bar{u})(\cdot,0) \in B^{1-\frac{\beta+1}{p}}_{p,p}(\R^{d-1})$ is well defined, and therefore $u(\cdot,x_d)$ does not have a limit in $\mathcal{D}'(\R^{d-1})$ as $x_d \to 0^+$. This counterexample, along with the positive results from Theorems \ref{thm:solvell}, \ref{thm:traceell}, \ref{thm:solvpar}, and \ref{thm:tracepar}, indicates that $F_d=0$ is a natural structural condition for the well-posedness of the boundary value problems (\ref{equationnonhomweighted})-(\ref{dirichlet}) and (\ref{equationparabolicdivnonhom})-(\ref{dirichletparabolic}) in the degenerate range $\beta \in [p-1,2p-1)$.


The main ideas of our proofs are as follows. In equations (\ref{equationnonhomweighted}) and (\ref{equationparabolicdivnonhom}), we first simplify the situation by reducing the coefficients to $a_{dj}=\delta_{dj}$, that is $a_{dd}=1$ and $a_{dj}=0$ for $j=1,\dots, d-1$, by a change of variables introduced in Section \ref{subsec:varchange}. Then, in Section \ref{subsec:extop}, we construct extension operators for Besov spaces that are suitable for equations with singular-degenerate coefficients. Using $a_{dj}=\delta_{dj}$ and the independence of the coefficients in $t$ and $x'$, we reformulate equations (\ref{equationnonhomweighted}) and (\ref{equationparabolicdivnonhom}) with $F_d=0$ as differentiated versions of related non-divergence form equations, which are studied in Section \ref{subsec:nondiv}. For $F_d\neq 0$, a different change of variables is employed to obtain an equivalent equation, enabling us to apply the solvability results from \cite{DongPhanIndiana2023}, which deal with equations under the conormal boundary condition. In the proofs of trace theorems, Theorems \ref{thm:traceell} and \ref{thm:tracepar}, we use a duality argument based on the isomorphism $B^s_{p,p}=(B^{-s}_{p',p'})^*$. To estimate the action $\langle u(\cdot, 0), \varphi \rangle$ on a Schwartz function $\varphi$, we introduce a 1D function $z(x_d) = \langle u(\cdot,x_d), v(\cdot,x_d)\rangle$, where $v(\cdot,x_d)$ is an appropriately chosen extension of $\varphi$. Using the equation satisfied by~$u$, we derive an ODE satisfied by $z$. Analyzing this ODE allows us to establish the existence of the boundary trace $u(\cdot,0)$, and to estimate the rate at which $u(\cdot,x_d)$ converges to $u(\cdot,0)$ as $x_d \to 0^+$.

At first glance, our proofs in both elliptic and parabolic settings may seem very similar. However, treating parabolic equations necessitates substantially more preparation and additional estimates. Specifically, for the existence of solutions in the parabolic case, we require certain estimates of the action of the non-local operator $(\partial_t-\Delta_{x'}+1)^{1/2}$. Additionally, for the convergence estimates of $u(\cdot, x_d)$ to $u(\cdot,0)$, in the duality argument we need to construct a family of test functions due to the presence of the $a_0(x_d)u_t$ term, whereas in the elliptic case a single test function suffices.

The remaining parts of the paper are organized as follows. In Section \ref{sec:mainres}, we introduce notation, define function spaces, and present the main results. In Section \ref{sec:auxtools}, we introduce Besov and Bessel-potential spaces, Hardy's inequalities, a change of variables, and auxiliary lemmas. In Section \ref{sec:nondiv}, we construct extension operators for Besov spaces and consider equations in non-divergence form. In Section \ref{sec:proofsell}, we provide the proofs of Theorems \ref{thm:solvell} and \ref{thm:traceell}. In Section \ref{sec:proofspar}, we prove Theorems \ref{thm:solvpar}, \ref{thm:tracepar},  \ref{thm:existsolpar}, and \ref{thm:existsolell}. Finally, Appendix \ref{appendix} contains the proof of Lemma \ref{lem:solutionoperatorwelldefined}.

\section{Function spaces and main results} \label{sec:mainres}
\subsection{Notation, function spaces and differential operators} For Banach spaces $X$ and $Y$, we write $X=Y$ if $X$ and $Y$ coincide as vector spaces and their norms are equivalent. We write $X \hookrightarrow Y$ or $Y \hookleftarrow X$ if $X \subset Y$ and there exists a constant $C>0$ such that $\Vert x \Vert_Y \leq C \Vert x \Vert_X$ for all $x \in X$. By $N_{a,b,\dots}$, we denote a positive constant that depends only on $a,b,\ldots$. We write $A \sim_{a,b,\dots} B$ if $N_{a,b,\dots}^{-1} A \leq B \leq N_{a,b,\dots} A$ for some $N_{a,b,\dots}>0$. For a real number $a$, we denote $a^+ = \max(a,0), a^- = \max(-a,0)$.

Let $p \in [1,\infty)$ and $\alpha, \beta \in \R$. We denote
\begin{align*}
    L_{p,\beta}(\R^d_+) & = L_p(\R^d_+, x_d^\beta dx), \\
    \Vert u \Vert_{ L_{p,\beta}(\R^d_+)}^p & = \int_{\R^d_+} |u(x)|^p x_d^{\beta} dx; \\
    L_{p,\beta}(\R^{d+1}_+) & = L_p(\R^{d+1}_+, x_d^\beta dxdt), \\
    \Vert u \Vert_{ L_{p,\beta}(\R^{d+1}_+)}^p & = \int_{\R^{d+1}_+} |u(t,x)|^p x_d^{\beta} dxdt.
\end{align*}
For an open set $\mathcal{O}$, we say $u \in L_{p,\text{loc}}(\mathcal{O})$ if $u \in L_p(K)$ for every compact $K \subset \mathcal{O}$.

We introduce the following weighted Sobolev spaces with first-order spatial derivatives:
\begin{align*}
W^1_{p,\beta}(\mathbb{R}^{d}_+) & = \{u\colon u, Du \in L_{p,\beta}(\mathbb{R}^{d}_+)\}, \\
\Vert u \Vert_{W^1_{p,\beta}(\mathbb{R}^{d}_+)} & = \Vert u \Vert_{L_{p,\beta}(\mathbb{R}^{d}_+)}+\Vert Du \Vert_{L_{p,\beta}(\mathbb{R}^{d}_+)}; \\
W^{0,1}_{p,\beta}(\R^{d+1}_+) & = \{u \colon u, Du \in L_{p,\beta}(\R^{d+1}_+)\} = L_p(\mathbb{R};W^1_{p,\beta}(\mathbb{R}^d_+)), \\
\Vert u \Vert_{W^{0,1}_{p,\beta}(\R^{d+1}_+)} & = \Vert u \Vert_{L_{p,\beta}(\R^{d+1}_+)} + \Vert Du \Vert_{L_{p,\beta}(\R^{d+1}_+)}.
\end{align*}

We also define the following weighted Sobolev spaces with second-order spatial derivatives and first-order time derivative (for parabolic spaces):
\begin{align*}
W^2_{p,\alpha,\beta}(\mathbb{R}^{d}_+) & = \{u \in W^1_{p,\beta}(\mathbb{R}^{d}_+): D_{x'}D u, D^2_d u + \alpha x_d^{-1} D_d u \in L_{p,\beta}(\mathbb{R}^{d}_+)\}, \\
    \Vert u \Vert_{W^2_{p,\alpha,\beta}(\mathbb{R}^{d}_+)} & =  \Vert u \Vert_{W^1_{p,\beta}(\mathbb{R}^{d}_+)} + \Vert D_{x'} D u\Vert_{L_{p,\beta}(\mathbb{R}^{d}_+)} + \Vert D^2_d u + \alpha x_d^{-1} D_d u\Vert_{L_{p,\beta}(\mathbb{R}^{d}_+)}; \\
W^2_{p,\beta,\mathcal{N}}(\mathbb{R}^{d}_+)  & = \{u \in W^2_{p,\beta}(\mathbb{R}^{d}_+): x_d^{-1} D_d u \in L_{p,\beta}(\mathbb{R}^{d}_+)\}, \\
    \Vert u \Vert_{W^2_{p,\beta,\mathcal{N}}(\mathbb{R}^{d}_+)} & = \Vert u \Vert_{W^2_{p,\beta}(\mathbb{R}^{d}_+)} + \Vert x_d^{-1}D_d u\Vert_{L_{p,\beta}(\mathbb{R}^{d}_+)}; \\
W^{1,2}_{p,\alpha,\beta}(\R^{d+1}_+) & = \{u: u_t, u, Du, D_{x'}Du, D^2_d u + \alpha x_d^{-1} D_d u \in L_{p,\beta}(\R^{d+1}_+)\}, \\
    \Vert u \Vert_{W^{1,2}_{p,\alpha,\beta}(\R^{d+1}_+)} & =  \Vert u_t \Vert_{L_{p,\beta}(\R^{d+1}_+)} + \Vert u \Vert_{L_{p,\beta}(\R^{d+1}_+)} + \Vert D u\Vert_{L_{p,\beta}(\R^{d+1}_+)} \\
    & + \Vert D_{x'}D u \Vert_{L_{p,\beta}(\R^{d+1}_+)} + \Vert D^2_d u + \alpha x_d^{-1} D_d u \Vert_{L_{p,\beta}(\R^{d+1}_+)}; \\
    W^{1,2}_{p,\beta,\mathcal{N}}(\R^{d+1}_+) & = \{u \in W^{1,2}_{p,\beta}(\R^{d+1}_+): x_d^{-1} D_d u \in L_{p,\beta}(\R^{d+1}_+)\}, \\
    \Vert u \Vert_{W^{1,2}_{p,\beta,\mathcal{N}}(\R^{d+1}_+)} & = \Vert u \Vert_{W^{1,2}_{p,\beta}(\R^{d+1}_+)} + \Vert x_d^{-1}D_d u\Vert_{L_{p,\beta}(\R^{d+1}_+)}.
\end{align*}
When $\alpha=0$, we write $W^{2}_{p,\beta}(\R^d_+)$ and $W^{1,2}_{p,\beta}(\R^{d+1}_+)$ instead of $W^{2}_{p,0,\beta}(\R^d_+)$ and $W^{1,2}_{p,0,\beta}(\R^{d+1}_+)$.

As shown in \cite{DongPhan2023JFA}, see also Section \ref{sec:nondiv} below, the spaces $W^2_{p,\alpha,\beta}(\mathbb{R}^{d}_+)$ and $W^{1,2}_{p,\alpha,\beta}(\R^{d+1}_+)$ are natural solution spaces for equations in non-divergence form. The letter $\mathcal{N}$ in $W^2_{p,\beta,\mathcal{N}}$ and $W^{1,2}_{p,\beta,\mathcal{N}}$ above stands for the Neumann boundary condition, since the condition $x_d^{-1}D_d u \in L_{p,\beta}(\mathbb{R}^{d}_+)$ is an analogue of the Neumann boundary condition $D_d u(\cdot,0)=0$ for functions in $W^2_{p,\beta}(\mathbb{R}^{d}_+)$ and $W^{1,2}_{p,\beta}(\R^{d+1}_+)$. The spaces $W^2_{p,\alpha,\beta}(\mathbb{R}^{d}_+)$, $W^2_{p,\beta,\mathcal{N}}(\R^d_+)$, and more general anisotropic weighted Sobolev spaces (with slightly different notation), were studied in detail in \cite{MetafuneNegroSpina2023TokyoJMath}.

We also define
\begin{align*}
	\widehat{W}^{-1}_{p,\beta}(\mathbb{R}^{d}_+) &= \{g \in \mathcal{D}'(\mathbb{R}^{d}_+)\colon \exists F_i, f \in L_{p,\beta}(\mathbb{R}^{d}_+) \; \text{s.t.} \; g = \sum_{i=1}^{d-1} D_i F_i+ f\}, \\
	\Vert g \Vert_{\widehat{W}^{-1}_{p,\beta}(\mathbb{R}^{d}_+)} &= \inf \left\{\sum_{i=1}^{d-1}\Vert F_i \Vert_{L_{p,\beta}(\mathbb{R}^{d}_+)} + \Vert f \Vert_{L_{p,\beta}(\mathbb{R}^{d}_+)} \colon g = \sum_{i=1}^{d-1} D_i F_i+ f \right\}; \\
\widehat{\mathbb{H}}^{-1}_{p,\beta}(\R^{d+1}_+) &= \{g \in \mathcal{D}'(\R^{d+1}_+)\colon \exists F_i, f \in L_{p,\beta}(\R^{d+1}_+) \; \text{s.t.} \; g = \sum_{i=1}^{d-1} D_i F_i+ f\}, \\
\Vert g \Vert_{\widehat{\mathbb{H}}^{-1}_{p,\beta}(\R^{d+1}_+)} &= \inf \left\{\sum_{i=1}^{d-1}\Vert F_i \Vert_{L_{p,\beta}(\R^{d+1}_+)} + \Vert f \Vert_{L_{p,\beta}(\R^{d+1}_+)} \colon g = \sum_{i=1}^{d-1} D_i F_i+ f \right\}.
\end{align*}
Here $\mathcal{D}'(\mathcal{O})$ denotes the space of distributions on an open set $\mathcal{O}$. In the definitions above, we use the hat notation $\widehat{W}$ and $\widehat{\mathbb{H}}$ to emphasize that differentiation in $x_d$ is not included.

By $\mathcal{E}$ and $\mathcal{P}$, we denote differential operators of the forms
\begin{gather} \label{eq:operatorE}
\mathcal{E} u = x_d^{-\alpha}\sum_{j=1}^d D_{d}(x_d^\alpha a_{dj}(x_d) D_{j}u), \\
\label{eq:operatorP}
\mathcal{P}u = a_0(x_d)u_t - x_d^{-\alpha}\sum_{j=1}^d D_{d}(x_d^\alpha a_{dj}(x_d) D_{j}u).
\end{gather}
The coefficients $\{a_{dj}(x_d)\}_{j=1}^d$ and $a_0(x_d)$ are assumed to be measurable functions on $\R_+$ satisfying the non-degeneracy and boundedness conditions
\begin{equation} \label{eq:ellipticityP}
    \kappa \leq a_0(x_d), a_{dd}(x_d) \leq \kappa^{-1}, \qquad |a_{dj}(x_d)| \leq \kappa^{-1}
\end{equation}
for all $x_d \in (0,\infty)$, $1 \leq j \leq d-1$, and some $\kappa \in (0,1]$.

For operators $\mathcal{E}$ and $\mathcal{P}$ as in (\ref{eq:operatorE}) and (\ref{eq:operatorP}), we define
\begin{align*}
W^1_{p,\beta}(\mathcal{E}, \R^d_+) & = \{u \in W^1_{p,\beta}(\mathbb{R}^d_+): \mathcal{E}u \in \widehat{W}^{-1}_{p,\beta}(\mathbb{R}^d_+)\}, \\
\Vert u \Vert_{W^1_{p,\beta}(\mathcal{E}, \R^d_+)} & = \Vert u \Vert_{W^1_{p,\beta}(\mathbb{R}^d_+)} + \Vert \mathcal{E}u \Vert_{\widehat{W}^{-1}_{p,\beta}(\mathbb{R}^d_+)}; \\
\mathcal{H}^1_{p,\beta}(\mathcal{P},\R^{d+1}_+) & = \{u \in W^{0,1}_{p,\beta}(\R^{d+1}_+): \mathcal{P}u \in \widehat{\mathbb{H}}^{-1}_{p,\beta}(\R^{d+1}_+)\}, \\
\Vert u \Vert_{\mathcal{H}^1_{p,\beta}(\mathcal{P},\R^{d+1}_+)} & = \Vert u \Vert_{W^{0,1}_{p,\beta}(\R^{d+1}_+)} + \Vert \mathcal{P}u \Vert_{\widehat{\mathbb{H}}^{-1}_{p,\beta}(\R^{d+1}_+)}.
\end{align*}
When $\mathcal{E} = x_d^{-\alpha}D_{d}(x_d^{\alpha}D_d\,) = D^2_d + \frac{\alpha}{x_d}D_d$, that is $a_{dd}=1$ and $a_{dj}=0$ for $1 \leq j \leq d-1$, we write $W^1_{p,\alpha,\beta}(\R^d_+)$ instead of $W^1_{p,\beta}(\mathcal{E}, \R^d_+)$.

As will be shown in Theorems \ref{thm:solvell} and \ref{thm:solvpar}, $W^1_{p,\beta}(\mathcal{E}, \R^d_+)$ and $\mathcal{H}^1_{p,\beta}(\mathcal{P},\R^{d+1}_+)$ are solution spaces for equations in divergence form with $F_d=0$ (\ref{equationnonhomweighted}) and (\ref{equationparabolicdivnonhom}). Theorems \ref{thm:traceell} and \ref{thm:tracepar} are trace theorems for these spaces in the degenerate range $\beta \in [p-1,2p-1)$.

We say that a function $a:\R_+\to \C$ is bounded uniformly continuous on $\R_+$, denoted as $a \in BUC(\R_+)$, if $a \in L_{\infty}(\R_+)$, and $\lim_{ \varepsilon \to 0^+} \omega_a(\varepsilon) = 0$, where
\begin{equation*}
\omega_a(\varepsilon) := \sup_{|s_1-s_2| \leq \varepsilon}|a(s_1)-a(s_2)|
\end{equation*}
is the modulus of continuity of $a$.

Throughout the paper, we understand the boundary conditions (\ref{dirichlet})-(\ref{dirichletparabolic}) in the sense of distributions.
\begin{defn} \label{def:trace}
Let $\mathcal{O}$ be an open subset of $\R^m$ for some $m \geq 1$. Suppose that $A$ is a Banach space continuously embedded in $\mathcal{D}'(\mathcal{O})$.

Let $u(y) \in C((0,\infty);A)$ and $U \in \mathcal{D}'(\mathcal{O})$. We say that $u(0)=U$ if 
\[
\lim_{y \to 0^+} u(y) = U \qquad \text{in} \; \mathcal{D}'(\mathcal{O}).
\]
\end{defn}
To apply Definition \ref{def:trace} to Sobolev spaces, we use the following
\begin{lem*}
Let $\mathcal{O}$ be an open subset in $\R^m$ for some $m \geq 1$, $I = (a, b)$, where $-\infty \leq a < b \leq +\infty$. Suppose $u(z,y) \in L_{1,\operatorname{loc}}(\mathcal{O} \times I)$ with $D_y u \in  L_{1,\operatorname{loc}}(\mathcal{O} \times I)$. Then $u$ has a version $v$, such that for a.e. $z \in \mathcal{O}$, the function
\begin{equation*}
y \mapsto v(z,y)
\end{equation*}
is absolutely continuous on all compact subintervals of $I$. The classical derivative $D_y v$, which exists a.e., agrees a.e. with $D_y u$.
\end{lem*}
See \cite[9.10.9]{BogachevSmolyanov} for the proof. As a corollary of the above lemma and the fundamental theorem of calculus, elements $u \in W^1_{p,\beta}(\R^d_+)$ have versions $u \in C((0,\infty); L_p(\R^{d-1}))$. We always use such identification if needed to define $u(\cdot, x_d)$ for \textit{all} $x_d > 0$.

\subsection{Main results: elliptic equations} We now state our main results for elliptic equations. Let $F, f \in L_{1,\text{loc}}(\R^d_+)$. We say that $u \in W^1_{1,\operatorname{loc}}(\mathbb{R}^d_+)$ is a weak solution to \eqref{equationnonhomweighted} if
\begin{equation} \label{eq:weaksolutionellipticdefinition}
    \int_{\mathbb{R}^d_+} x_d^\alpha(a_{ij}D_j u D_i \varphi + \lambda c u \varphi) dx =  \int_{\mathbb{R}^d_+} x_d^{\alpha}(F_i D_i\varphi + f \varphi) dx
\end{equation}
for all $\varphi \in C_0^\infty(\mathbb{R}^d_+)$. By $u \in W^1_{1,\operatorname{loc}}(\mathbb{R}^d_+)$ we mean $u,Du \in L_{1,\operatorname{loc}}(\R^d_+)$.

Our first main result is about the solvability of elliptic equations in divergence form. Everywhere below in Theorems \ref{thm:solvell}--\ref{thm:existsolpar} we have $N=N(d,\alpha,p,\beta, \kappa) > 0$.
\begin{thm} \label{thm:solvell}
Let $\alpha \in (-\infty,1), p \in (1,\infty), \beta \in [p-1, 2p-1)$, and $\kappa \in (0,1]$. Suppose that \eqref{eq:ellipticitycond} is satisfied. Then, for any $\lambda > 0$, $F \in L_{p,\beta}(\mathbb{R}^d_+; \mathbb{R}^d)$ with $F_d=0$, $f \in L_{p,\beta}(\R^d_+)$, and $U \in B^{1-\frac{\beta+1}{p}}_{p,p}(\mathbb{R}^{d-1})$, there exists a unique weak solution $u \in W^1_{p,\beta}(\mathbb{R}^d_+)$ to \eqref{equationnonhomweighted}--\eqref{dirichlet}. Moreover, $u$ satisfies
\begin{equation} \label{eq:estsolell}
\begin{gathered}
\Vert Du \Vert_{L_{p,\beta}(\mathbb{R}^{d}_+)} + \sqrt{\lambda} \Vert u \Vert_{L_{p,\beta}(\mathbb{R}^{d}_+)} \leq N \Vert F \Vert_{L_{p,\beta}(\mathbb{R}^{d}_+)} + \frac{N}{\sqrt{\lambda}} \Vert f \Vert_{L_{p,\beta}(\mathbb{R}^{d}_+)} \\
+ N\bigl(\lambda^{1-\frac{\beta+1}{2p}} \Vert (\lambda-\Delta_{x'})^{-1/2} U \Vert_{L_p(\R^{d-1})}+[(\lambda-\Delta_{x'})^{-1/2} U]_{B^{2-\frac{\beta+1}{p}}_{p,p}(\mathbb{R}^{d-1})}\bigr). \\
\end{gathered}
\end{equation}
\end{thm}
See Section \ref{subsec:besselbesovspaces} for a precise definition of Besov spaces $B^s_{p,p}(\R^{d-1})$.

Our next result is the trace theorem for the spaces $W^1_{p,\beta}(\mathcal{E}, \R^d_+)$.
\begin{thm} \label{thm:traceell}
Let $\alpha \in (-\infty,1), p \in (1,\infty), \beta \in [p - 1, 2p-1)$, and $\kappa \in (0, 1]$. Let $\mathcal{E}$ be as in \eqref{eq:operatorE} and suppose that \eqref{eq:ellipticityP} holds. Then, for every $u \in W^1_{p,\beta}(\mathcal{E}, \R^d_+)$ there exists a boundary trace $U = u(\cdot, 0) \in B^{1-\frac{\beta+1}{p}}_{p,p}(\mathbb{R}^{d-1})$. It satisfies
\[\Vert U \Vert_{B^{1-\frac{\beta+1}{p}}_{p,p}(\mathbb{R}^{d-1})} \leq N \Vert u \Vert_{W^1_{p,\beta}(\mathcal{E}, \R^d_+)}.
\]
Moreover,
\begin{itemize}
\item[(i)] if $\beta \in (p-1,2p-1)$, or $\beta = p-1, \alpha \in (-\infty, \frac{p-1}{2p-1}]$, then
\begin{equation*}
u(\cdot,x_d) \to U \quad \text{strongly in $B^{1-\frac{\beta+1}{p}}_{p,p}(\mathbb{R}^{d-1})$ as $x_d \to 0^+$};
\end{equation*}
\item[(ii)] if $\beta = p-1, \alpha \in (\frac{p-1}{2p-1}, 1)$, then for any $\varepsilon > 0$
\begin{equation*}
u(\cdot,x_d) \to U \quad \text{strongly in $B^{-\varepsilon}_{p,p}(\mathbb{R}^{d-1})$ as $x_d \to 0^+$}.
\end{equation*}
\end{itemize}
\end{thm}

The next result is about the existence of solutions to elliptic equations when $F_d \neq 0$. In view of the counterexample from Proposition \ref{prop:counterex}, boundary condition is not imposed.
\begin{thm} \label{thm:existsolell}
    Let $\alpha \in (-\infty, 1), p \in (1,\infty), \beta \in (p - 1, 2p-1)$, and $\kappa \in (0,1]$. Suppose that \eqref{eq:ellipticitycond} is satisfied. Then, for any $\lambda > 0, F \in L_{p,\beta}(\mathbb{R}^d_+;\R^d)$, and $f \in L_{p,\beta}(\R^d_+)$, there exists a weak solution $u \in W^1_{p,\beta}(\mathbb{R}^d_+)$ to \eqref{equationnonhomweighted}, that satisfies
\begin{equation*}
    \Vert Du \Vert_{L_{p,\beta}(\mathbb{R}^d_+)} + \sqrt{\lambda} \Vert u \Vert_{L_{p,\beta}(\mathbb{R}^d_+)} \leq N\Vert F \Vert_{L_{p,\beta}(\mathbb{R}^d_+)} + N \frac{1}{\sqrt{\lambda}} \Vert f \Vert_{L_{p,\beta}(\mathbb{R}^d_+)}.
\end{equation*}
\end{thm}

\subsection{Main results: parabolic equations} We now state the corresponding results for parabolic equations. Let $F,f \in L_{1,\text{loc}}(\R^{d+1}_+)$. We say that $u \in W^{0,1}_{1,\text{loc}}(\R^{d+1}_+)$ is a weak solution to \eqref{equationparabolicdivnonhom} if
\begin{equation} \label{eq:weaksolell} 
    \int_{\R^{d+1}_+} x_d^\alpha (-a_0 u \varphi_t + \lambda c u \varphi) dxdt + \int_{\R^{d+1}_+}x_d^\alpha (a_{ij}D_j u- F_i)D_i\varphi dxdt = \int_{\R^{d+1}_+} x_d^\alpha f\varphi dxdt
\end{equation}
for all $\varphi \in C_0^\infty(\R^{d+1}_+)$. By $u \in W^{0,1}_{1,\text{loc}}(\R^{d+1}_+)$ we mean $u,Du \in L_{1,\text{loc}}(\R^{d+1}_+)$.

Our next theorem is the solvability result for parabolic equations in divergence form.
\begin{thm} \label{thm:solvpar}
Let $\alpha \in (-\infty, 1), p \in (1,\infty), \beta \in [p-1,2p-1)$, and $\kappa \in (0,1]$. Suppose that \eqref{eq:ellipticitycond} is satisfied. Then, for any $\lambda > 0$, $F\in  L_{p,\beta}(\mathbb{R}^{d+1}_+; \mathbb{R}^d)$ with $F_d=0$, $f \in  L_{p,\beta}(\mathbb{R}^{d+1}_+)$, and $U \in B^{\frac{1}{2}-\frac{\beta+1}{2p}, 1-\frac{\beta+1}{p}}_{p,p}(\mathbb{R}^{d})$, there exists a unique weak solution $u \in W^{0,1}_{p,\beta}(\mathbb{R}^{d+1}_+)$ to \eqref{equationparabolicdivnonhom}--\eqref{dirichletparabolic}. Moreover, $u$ satisfies
\begin{equation} \label{eq:estsolpar}
\begin{gathered}
\Vert Du \Vert_{L_{p,\beta}(\mathbb{R}^{d+1}_+)} + \sqrt{\lambda} \Vert u \Vert_{L_{p,\beta}(\mathbb{R}^{d+1}_+)} \leq N \Vert F \Vert_{L_{p,\beta}(\mathbb{R}^{d+1}_+)} + \frac{N}{\sqrt{\lambda}} \Vert f \Vert_{L_{p,\beta}(\mathbb{R}^{d+1}_+)} \\
+ N\bigl(\lambda^{1-\frac{\beta+1}{2p}} \Vert (\partial_t-\Delta_{x'}+\lambda)^{-\frac{1}{2}} U \Vert_{L_p(\R^d)}+[(\partial_t-\Delta_{x'}+\lambda)^{-\frac{1}{2}} U]_{B^{1-\frac{\beta+1}{2p}, 2-\frac{\beta+1}{p}}_{p,p}(\mathbb{R}^{d})}\bigr).
\end{gathered}
\end{equation}
\end{thm}
See Section \ref{subsec:besselbesovspaces} for precise definitions of Besov spaces $B^{s/2,s}_{p,p}(\R^{d})$ and the operator $(\partial_t-\Delta_{x'}+\lambda)^{-\frac{1}{2}}$.

Our next result is the trace theorem for parabolic spaces $\mathcal{H}^1_{p,\beta}(\mathcal{P},\mathbb{R}^{d+1}_+)$.
\begin{thm} \label{thm:tracepar}
Let $\alpha \in (-\infty, 1), p \in (1,\infty), \beta \in [p-1,2p-1)$, and $\kappa \in (0, 1]$. Let $\mathcal{P}$ be as in \eqref{eq:operatorP} and suppose that \eqref{eq:ellipticityP} holds. Then, for every $u \in \mathcal{H}^1_{p,\beta}(\mathcal{P},\mathbb{R}^{d+1}_+)$ there exists a boundary trace $U = u(\cdot,0) \in B^{\frac{1}{2}-\frac{\beta+1}{2p},1-\frac{\beta+1}{p}}_{p,p}(\mathbb{R}^{d})$. It satisfies 
\[
\Vert U \Vert_{B^{\frac{1}{2}-\frac{\beta+1}{2p},1-\frac{\beta+1}{p}}_{p,p}(\mathbb{R}^{d})} \leq N \Vert u \Vert_{\mathcal{H}^1_{p,\beta}(\mathcal{P},\mathbb{R}^{d+1}_+)}.
\]
Moreover,
\begin{itemize}
\item[(i)] if $\beta \in (p-1,2p-1)$, then
\begin{equation*}
u(\cdot,x_d) \to U \quad \text{weakly in $B^{\frac{1}{2}-\frac{\beta+1}{2p},1-\frac{\beta+1}{p}}_{p,p}(\mathbb{R}^{d})$ as $x_d \to 0^+$};
\end{equation*}
\item[(ii)] if $\beta = p-1$, then for any $\varepsilon > 0$
\begin{equation*}
u(\cdot,x_d) \to U \quad \text{weakly in $B^{-\varepsilon/2,-\varepsilon}_{p,p}(\mathbb{R}^{d})$ as $x_d \to 0^+$}.
\end{equation*}
\end{itemize}
If $a_0', a_{dd}' \in BUC(\R_+)$, the convergence in {\rm (i)} and {\rm (ii)} is strong.
\end{thm}

Our last main result is about the solvability of parabolic equations when $F_d \neq 0$.
\begin{thm} \label{thm:existsolpar}
    Let $\alpha \in (-\infty, 1), p \in (1,\infty), \beta \in (p - 1, 2p-1)$, and $\kappa \in (0,1]$. Suppose that \eqref{eq:ellipticitycond} is satisfied. Then, for any $\lambda > 0$, $F \in L_{p,\beta}(\R^{d+1}_+;\R^d)$ and $f \in L_{p,\beta}(\R^{d+1}_+)$, there exists a weak solution $u \in W^{0,1}_{p,\beta}(\R^{d+1}_+)$ to \eqref{equationparabolicdivnonhom}, that satisfies
    \[
    \Vert Du \Vert_{L_{p,\beta}(\R^{d+1}_+)} + \sqrt{\lambda} \Vert u \Vert_{L_{p,\beta}(\R^{d+1}_+)} \leq N\Vert F \Vert_{L_{p,\beta}(\R^{d+1}_+)} + N \frac{1}{\sqrt{\lambda}} \Vert f \Vert_{L_{p,\beta}(\R^{d+1}_+)}. 
    \]
\end{thm}

\section{Auxiliary tools} \label{sec:auxtools}
In this section, we introduce some definitions and auxiliary results.
\subsection{Distributions and Fourier transform} Let $m \in \mathbb{N}$.
By $\mathcal{S}(\mathbb{R}^m)$ we denote the Schwartz space of smooth rapidly decreasing functions on $\R^m$. Its continuous dual space, i.e., the space of tempered distributions on $\mathbb{R}^m$, is denoted by $\mathcal{S}'(\mathbb{R}^m)$. We denote by
\begin{gather*}
\langle \cdot, \cdot \rangle: \mathcal{S}'(\mathbb{R}^{m}) \times  \mathcal{S}(\mathbb{R}^{m}) \rightarrow \mathbb{C}, \\
\langle f, \varphi \rangle = f(\varphi),
\end{gather*}
the dual pairing of $\mathcal{S}'(\mathbb{R}^{m})$ and $\mathcal{S}(\mathbb{R}^{m})$. For a pair of functions $f,g$ on $\R^m$, we denote
\[
\langle f , g \rangle = \int_{\R^m} fg,
\]
whenever $fg \in L_1(\R^m)$.
We define the Fourier transform on $\mathcal{S}(\mathbb{R}^{d-1}_{x'})$ or $\mathcal{S}(\mathbb{R}^{d}_{t,x'})$ by
\begin{align*}
    \mathcal{F}_{x'\to\xi'}f(\xi') & = \int_{\mathbb{R}^{d-1}}f(x')e^{-ix'\cdot \xi'}dx', \\
    \mathcal{F}_{(t,x')\to (\tau,\xi')}f(\tau,\xi') & = \int_{\mathbb{R}^{d}}f(t,x')e^{-i(t\tau+x'\cdot \xi')}dtdx'.
\end{align*}
The Fourier transform is extended to $\mathcal{S}'$ by $\langle \mathcal{F}f, \varphi \rangle = \langle f, \mathcal{F} \varphi \rangle.$

For $\mu > 0$ and $s \in \R$, let
\begin{align*}
    (\mu - \Delta_{x'})^s f & = \mathcal{F}^{-1}[(\mu+|\xi'|^2)^{s} \mathcal{F}f], \\
    (\partial_t - \Delta_{x'} + \mu)^{s} f & = \mathcal{F}^{-1}[(i\tau+|\xi'|^2+\mu)^{s} \mathcal{F}f].
\end{align*}

\subsection{Bessel-potential and Besov spaces} \label{subsec:besselbesovspaces}
Let $p \in (1, \infty), s \in \R$. The isotropic Bessel potential spaces are defined by
\begin{align*}
    H^s_{p}(\R^{d-1}) & = (1-\Delta_{x'})^{-s/2}L_p(\R^{d-1}), \\
    \Vert f \Vert_{H^s_{p}(\R^{d-1})} & = \Vert (1-\Delta_{x'})^{s/2} f \Vert_{L_{p}(\R^{d-1})}.
\end{align*}
Let $E$ be a Banach space. The parabolic $E$-valued Bessel potential spaces are defined by
\begin{align*}
    H^{s/2,s}_p(\R^d;E) & = (\partial_t - \Delta_{x'} + 1)^{-s/2}L_p(\R^d;E), \\
    \Vert u \Vert_{H^{s/2,s}_p(\R^d;E)} & = \Vert (\partial_t - \Delta_{x'} + 1)^{s/2} u \Vert_{L_p(\R^d;E)}.
\end{align*}
Here $L_p(\R^d;E)$ is the $L_p$ space of strongly measurable $E$-valued functions on $\R^d$. In the scalar case $E=\mathbb{C}$, we write $H^{s/2,s}_{p}(\R^{d})$ instead of $H^{s/2,s}_p(\R^d;E)$. It is well-known that when $s=m \in \mathbb{N}_0$, Bessel spaces coincide with Sobolev spaces: $H^m_p(\R^{d-1})=W^m_p(\R^{d-1})$. When $s = 2m, m \in \mathbb{N}_0$, it holds that $H^{m,2m}(\R^d)=W^{m,2m}_p(\R^d)$.

Let $p,q \in [1,\infty]$. For $s \in (0,2)$, we define isotropic Besov spaces by
\begin{gather*}
    B^{s}_{p,q}(\mathbb{R}^{d-1}) = \{f \in L_p(\mathbb{R}^{d-1}): [f]_{B^{s}_{p,q}(\mathbb{R}^{d-1})} < \infty\}, \\
    \Vert f \Vert_{B^{s}_{p,q}(\mathbb{R}^{d-1})} = \Vert f \Vert_{L_p(\mathbb{R}^{d-1})} + [f]_{B^{s}_{p,q}(\mathbb{R}^{d-1})},
\end{gather*}
where
\begin{equation*}
    [f]_{B^{s}_{p,q}(\mathbb{R}^{d-1})} = \left\{\int_{\R^{d-1}}\left[\frac{\Vert f(x'+h)-2f(x')+f(x'-h) \Vert_{L_p^{x'}(\R^{d-1})}}{|h|^s}\right]^q\frac{dh}{|h|^{d-1}}\right\}^{1/q}
\end{equation*}
for $1 \leq q < \infty$, and for $q = \infty$
\begin{equation*}
    [f]_{B^{s}_{p,\infty}(\mathbb{R}^{d-1})} = \esssup_{h \in \R^{d-1}} \left\{|h|^{-s}\Vert f(x'+h)-2f(x')+f(x'-h) \Vert_{L_p(\R^{d-1})}\right\}.
\end{equation*}
For $s \in (0,2)$, we define parabolic Besov spaces by 
\begin{gather*}
    B^{s/2,s}_{p,q}(\R^d) = \{f \in L_p(\R^d): [f]_{B^{s/2,s}_{p,q}(\R^d)} < \infty\}, \\
    \Vert f \Vert_{B^{s/2,s}_{p,q}(\R^d)} = \Vert f \Vert_{L_p(\R^d)} + [f]_{B^{s/2,s}_{p,q}(\R^d)},
\end{gather*}
where
\begin{gather*}
[f]_{B^{s/2,s}_{p,q}(\R^d)} = \left\{\int_{\R^{d-1}}\left[\frac{\Vert f(t,x'+h)-2f(t, x')+f(t, x'-h) \Vert_{L_p^{t,x'}(\R^{d})}}{|h|^s}\right]^q\frac{dh}{|h|^{d-1}}\right\}^{1/q} \\
+ \left\{\int_{\R}\left[\frac{\Vert f(t+\tau,x')-f(t, x') \Vert_{L_p(\R^{d})}}{|\tau|^{s/2}}\right]^q\frac{d\tau}{|\tau|}\right\}^{1/q}
\end{gather*}
for $1 \leq q < \infty$, and a natural modification for $q = \infty$.

For $s \leq 0$ or $s \geq 2$, $p \in (1,\infty)$, and $q \in [1,\infty]$, we set
\begin{align*}
	B^{s}_{p,q}(\R^{d-1}) & = (1-\Delta_{x'})^{m/2}B^{s+m}_{p,q}(\R^{d-1}), \\
    \Vert f \Vert_{B^{s}_{p,q}(\R^{d-1})} & = \Vert (1-\Delta_{x'})^{-m/2}f \Vert_{B^{s+m}_{p,q}(\R^{d-1})}, \\
    B^{s/2,s}_{p,q}(\mathbb{R}^{d}) & = (\partial_t - \Delta_{x'} + 1)^{m/2}B^{(s+m)/2,s+m}_{p,q}(\mathbb{R}^{d}), \\
    \Vert f \Vert_{B^{s/2,s}_{p,q}(\mathbb{R}^{d})} & = \Vert (\partial_t - \Delta_{x'} + 1)^{-m/2}f \Vert_{B^{(s+m)/2,s+m}_{p,q}(\mathbb{R}^{d})},
\end{align*}
where $m$ is the integer such that $s+m \in (0,1]$. It should be noted that Besov spaces with $s \geq 2$ can also be characterized by norms based on higher-order differences, similar to the $s \in (0,2)$ case. Moreover, the Littlewood-Paley theory allows for unified definitions of norms on $B^s_{p,q}$ and $B^{s/2,s}_{p,q}$, with single formulae valid for all $s \in \R$. See \cite[Chapter 2]{Triebel} for isotropic spaces and \cite[Chapter 5]{Triebel3} for general anisotropic spaces.

In the remaining part of this subsection, we state some properties of Besov and Bessel spaces.
\begin{lem} \label{lem:isomorphisms}
Let $s, r \in \R, p \in (1,\infty), q \in [1,\infty], \mu > 0$, and $X \in \{B_{p,q}, H_p\}$. Operators
\begin{gather*}
(\mu-\Delta_{x'})^{r/2}:X^s(\R^{d-1}) \to X^{s-r}(\R^{d-1}), \\
(\partial_t-\Delta_{x'}+\mu)^{r/2}:X^{s/2,s}(\R^{d}) \to X^{(s-r)/2,s-r}(\R^{d})
\end{gather*}
are isomorphisms of Banach spaces.
\end{lem}

For $p \in [1,\infty]$, its adjoint H\"older exponent $p'$ is defined by $\frac{1}{p}+\frac{1}{p'}=1$.  
\begin{lem} \label{lem:dualitybesov}
Let $s \in \mathbb{R}$, $p \in (1,\infty)$. Then
\begin{gather*}
	(B^{s}_{p,p}(\mathbb{R}^{d-1}))^* = B^{-s}_{p',p'}(\mathbb{R}^{d-1}), \\
    (B^{s/2,s}_{p,p}(\R^d))^* = B^{-s/2,-s}_{p',p'}(\R^d).
\end{gather*}
\end{lem}
Lemma \ref{lem:dualitybesov} needs to be understood in the following sense. A distribution $f \in \mathcal{S}'(\R^d)$ belongs to $B^{-s/2,-s}_{p',p'}(\R^d)$ if and only if
\begin{equation*}
    |\langle f, \varphi \rangle| \leq c \Vert \varphi \Vert_{B^{s/2,s}_{p,p}(\R^d)}
\end{equation*}
holds for all $\varphi \in \mathcal{S}(\R^d)$ with a constant $c$ independent of $\varphi$. The smallest constant $c$ for which the above inequality holds is comparable to the norm $\Vert f \Vert_{B^{-s/2,-s}_{p',p'}(\R^d)}$. Similarly for $B^s_{p,p}(\R^{d-1})$.

Let $[\;, \;]_{\theta}$ and $(\;, \;)_{\theta, q}$ be respectively the complex and real interpolation functors. Here $\theta \in (0,1)$, $q \in [1,\infty]$. We refer to \cite{Triebel} for precise definitions and results of interpolation theory.

\begin{lem} \label{lem:interpbesov}
Let $p \in (1,\infty), s_0, s_1 \in \R$, and $\theta \in (0,1)$. Then
\begin{gather*}
[B^{s_0}_{p,p}(\R^{d-1}), B^{s_1}_{p,p}(\R^{d-1})]_{\theta} = B^{s_{\theta}}_{p,p}(\R^{d-1}), \\
[B^{s_0/2,s_0}_{p,p}(\R^{d}), B^{s_1/2,s_1}_{p,p}(\R^{d})]_{\theta} = B^{s_{\theta}/2, s_{\theta}}_{p,p}(\R^{d}),
\end{gather*}
where $s_{\theta}=(1-\theta)s_0+\theta s_1$.
\end{lem}

For the proofs of Lemmas \ref{lem:isomorphisms}, \ref{lem:dualitybesov}, \ref{lem:interpbesov}, and further properties of Besov and Bessel spaces, see \cite{Triebel} for isotropic spaces, and \cite{Grubb2018JFA}, \cite[Chapter 5]{Triebel3} for parabolic spaces.

We will need the following lemma in Section \ref{sec:nondiv}.
\begin{lem} \label{lem:interpbessel}
    Let $\{E_0,E_1\}$ be an interpolation couple of Banach spaces, $s \in \R$, $p \in (1,\infty), \theta \in (0,1)$. Then
    \begin{equation*}
        [H^{s/2,s}_p(\R^d;E_0),H^{s/2,s}_p(\R^d;E_1)]_{\theta} = H^{s/2,s}_p(\R^d;[E_0,E_1]_{\theta}).
    \end{equation*}
\end{lem}
\begin{proof} Let $\Theta = (\partial_t-\Delta_{x'}+1)^{1/2}$. Then
    \begin{gather*}
        [H^{s/2,s}_p(\R^d;E_0),H^{s/2,s}_p(\R^d;E_1)]_{\theta} = [\Theta^{-s}L_p(\R^d;E_0),\Theta^{-s}L_p(\R^d;E_1)]_{\theta} \\
        = \Theta^{-s}[L_p(\R^d;E_0), L_p(\R^d;E_1)]_{\theta}= [\text{by \cite[1.18.4]{Triebel}}] \\
        = \Theta^{-s}L_p(\R^d;[E_0,E_1]_{\theta}) = H^{s/2,s}_p(\R^d;[E_0,E_1]_{\theta}).
    \end{gather*}
\end{proof}

\subsection{Hardy's inequalities and weighted spaces}
We will frequently use Hardy's inequalities in the following form. Let $G = \mathbb{R}^{d-1}$ or $G = \mathbb{R}^{d}$.
\begin{lem} \label{lem:Hardy}
    Let $p \in [1,\infty)$, $\alpha \in \mathbb{R}$.
    \begin{enumerate}
        \item[(i)] Suppose that $\beta > (1+\alpha)p-1$. Let $u \in L_{p,\beta}(G \times \mathbb{R}_+)$ with $D_d u+\frac{\alpha}{x_d}u \in L_{p,\beta}(G \times \mathbb{R}_+)$. Then $u/x_d \in L_{p,\beta}(G \times \mathbb{R}_+)$, and
    \begin{equation*}
        \Vert u/x_d \Vert_{L_{p,\beta}(G \times \mathbb{R}_+)} \leq \frac{1}{\frac{\beta-\alpha p+1}{p}-1}\Vert D_d u+\frac{\alpha}{x_d}u \Vert_{L_{p,\beta}(G \times \mathbb{R}_+)}.
    \end{equation*}
    \item[(ii)] Suppose that $\beta < (1+\alpha)p-1$. Let $u \in L_{p,\beta}(G \times \mathbb{R}_+)$ with $D_d u +\frac{\alpha}{x_d}u\in L_{p,\beta}(G \times \mathbb{R}_+)$ and $\lim_{x_d \to 0^+} x_d^{\alpha}u(\cdot,x_d)=0$ a.e. on $G$. Then $u/x_d\in L_{p,\beta}(G \times \mathbb{R}_+)$, and
    \begin{equation*}
        \Vert u/x_d \Vert_{L_{p,\beta}(G \times \mathbb{R}_+)} \leq \frac{1}{1-\frac{\beta-\alpha p+1}{p}}\Vert D_d u+\frac{\alpha}{x_d}u \Vert_{L_{p,\beta}(G \times \mathbb{R}_+)}.
    \end{equation*}
    Suppose $\beta \leq \alpha p - 1$. If $u \in L_{p,\beta}(G \times \mathbb{R}_+)$ and $D_d u +\frac{\alpha}{x_d}u\in L_{p,\beta}(G \times \mathbb{R}_+)$, then the condition $\lim_{x_d \to 0^+} x_d^{\alpha}u(\cdot,x_d)=0$ a.e. on $G$ holds automatically.
    \end{enumerate}
\end{lem}
When $\alpha=0$, the inequalities above are classical Hardy's inequalities; see \cite{Grisvard}. The case of general $\alpha$ reduces to $\alpha=0$, since the conditions on $u$ can be rewritten as $x_d^{\alpha}u, D_d(x_d^{\alpha} u) \in L_{p,\beta-\alpha p}(G\times \R_+)$.

\begin{cor}
In the setting of Lemma {\rm \ref{lem:Hardy}}, both $D_d u$ and $u/x_d$ are controlled by $D_d u+\frac{\alpha}{x_d}u$:
\begin{equation} \label{eq:controlderbymixed}
\Vert D_d u \Vert_{L_{p,\beta}(G \times \mathbb{R}_+)} + \Vert u/x_d \Vert_{L_{p,\beta}(G \times \mathbb{R}_+)} \leq N_{\alpha,p,\beta}\Vert D_d u+\frac{\alpha}{x_d}u \Vert_{L_{p,\beta}(G \times \mathbb{R}_+)}.
\end{equation}
\end{cor}
\begin{lem} \label{lem:structureW1}
	Let $\alpha, \beta \in \R$.
    \begin{enumerate}
    \item[(i)] Let $p \in [1,\infty)$, $\beta \leq \alpha p - 1$ or $\beta > (1+\alpha)p - 1$. Then
    \begin{equation*}
    \begin{gathered}
        W^2_{p,\alpha,\beta}(\mathbb{R}^{d}_+) = W^2_{p,\beta,\mathcal{N}}(\mathbb{R}^{d}_+), \\
        W^{1,2}_{p,\alpha,\beta}(\R^{d+1}_+) = W^{1,2}_{p,\beta,\mathcal{N}}(\R^{d+1}_+),
    \end{gathered}
    \end{equation*}
    with the equivalence of norms.
    \item[(ii)] Let $p \in (1, \infty)$. Then
\begin{equation*}
W^1_{p,\alpha,\beta}(\R^d_+) = (1-\Delta_{x'})^{\frac{1}{2}}W^2_{p,\alpha,\beta}(\R^d_+).
\end{equation*}
The norms $\Vert \cdot \Vert_{W^1_{p,\alpha,\beta}(\R^d_+)}$ and $\Vert (1-\Delta_{x'})^{-\frac{1}{2}}\cdot \Vert_{W^2_{p,\alpha,\beta}(\R^d_+)}$ are equivalent.
    \end{enumerate}
\end{lem}
Recall that by definition $u \in W^1_{p,\alpha,\beta}(\R^d_+)$ means $u \in W^1_{p,\beta}(\R^d_+)$ with $D^2_d u + \frac{\alpha}{x_d}D_d u = \sum_{i=1}^{d-1}D_i F_i + f$ for some $F_i, f \in L_{p,\beta}(\R^d_+)$.
\begin{proof}
Let $v \in W^2_{p,\alpha,\beta}(\mathbb{R}^{d}_+)$ or $v \in W^{1,2}_{p,\alpha,\beta}(\R^{d+1}_+)$. Part (i) follows immediately from Lemma \ref{lem:Hardy} and (\ref{eq:controlderbymixed}) applied with $u=D_d v$. Part (ii) follows from the definitions of the spaces involved, and the boundedness of $D_{x'}(1-\Delta_{x'})^{-\frac{1}{2}}$ on $L_p(\R^{d-1})$.
\end{proof}

\subsection{A useful change of variables} \label{subsec:varchange}
We now introduce a change of variables $\Phi:\R^d_+ \to \R^d_+$ such that in the new variables equations (\ref{equationnonhomweighted}) and (\ref{equationparabolicdivnonhom}) have $a_{dj} = \delta_{dj}$. That is,
\[
a_{dj} = 0 \quad \text{and} \quad a_{dd} = 1
\]
for $j=1,\dots,d-1$. Here we assume that $\alpha < 1$ and (\ref{eq:ellipticitycond}) hold.

The mapping $y=\Phi(x)$ is given by
\begin{equation} \label{eq:coordchange}
\begin{gathered}
    y_j = x_j - \int_{0}^{x_d}\frac{a_{dj}(s)}{a_{dd}(s)} ds, \quad 1 \leq j \leq d-1, \\
    y_d = \left((1-\alpha)\int_{0}^{x_d}\frac{s^{-\alpha}}{a_{dd}(s)} ds\right)^{\frac{1}{1-\alpha}}.
\end{gathered}
\end{equation}
The following identities hold:
\begin{equation*}
\begin{gathered}
    D_{x_j} = D_{y_j}, \quad 1 \leq j \leq d-1, \\
    D_{x_d} = -\sum_{j=1}^{d-1}\frac{a_{dj}(x_d)}{a_{dd}(x_d)}D_{y_j}
    +\frac{y_d^{\alpha}}{x_d^{\alpha}}\frac{1}{a_{dd}(x_d)}D_{y_d}.
\end{gathered}
\end{equation*}

It is then straightforward to verify that the mapping $\Phi=\Phi(x)$ is a one-to-one mapping from $\R^d_+$ onto $\R^d_+$, bi-Lipschitz, and is an identity on the boundary: 
\begin{gather*}
\Vert D\Phi \Vert_{L_{\infty}(\R^d_+)}+\Vert D\Phi^{-1} \Vert_{L_{\infty}(\R^d_+)} \leq N_{d,\alpha, \kappa}, \\
\Phi(x',0) = (x',0).
\end{gather*}
In the $y$-variables, equation (\ref{equationparabolicdivnonhom}) reads
\begin{equation*} 
    y_d^\alpha \Tilde{a}_0(y_d) u_t - \sum_{i,j=1}^d D_{y_i}(y_d^\alpha \Tilde{a}_{ij}(y_d) D_{y_j}u)+\lambda y_d^\alpha \Tilde{c}(y_d) u = -\sum_{j=1}^d D_{y_i}(y_d^\alpha \Tilde{F}_i) + y_d^\alpha \Tilde{f},
\end{equation*}
and similarly for elliptic equation without the $u_t$ term, where
\begin{equation*}
\begin{array}{rlr}
    \Tilde{a}_{ij}(y_d) = & (x_d/y_d)^{2\alpha}(a_{ij}a_{dd}-a_{id}a_{dj}), & \quad 1 \leq i,j \leq d - 1, \\
    \Tilde{a}_{id}(y_d) = & (x_d/y_d)^{\alpha}(a_{id}-a_{di}), &  1 \leq i \leq d - 1, \\
    \Tilde{a}_{dj}(y_d) = & 0, & 1 \leq j \leq d - 1, \\
    \Tilde{a}_{dd}(y_d) = & 1, & \\
    \Tilde{a}_{0}(y_d) = & (x_d/y_d)^{2\alpha} a_{dd}(x_d) a_0(x_d), & \\
    \Tilde{c}(y_d) = & (x_d/y_d)^{2\alpha} a_{dd}(x_d) c(x_d), & \\
    \Tilde{F}_i = & (x_d/y_d)^{2\alpha} (a_{dd} F_i - a_{di}F_d), &  1 \leq i \leq d-1, \\
    \Tilde{F}_d = & (x_d/y_d)^{\alpha} F_d, & \\
    \Tilde{f} = & (x_d/y_d)^{2\alpha} a_{dd}f. & 
\end{array}
\end{equation*}
It is straightforward to verify that if $a_0'(x_d), a_{dd}'(x_d) \in BUC(\R_+)$, then $\Tilde{a}_0'(y_d) \in BUC(\R_+)$. This justifies the use of change of variables $\Phi$ in the proof of Theorem \ref{thm:tracepar}.

\section{Extension operators, equations in non-divergence form} \label{sec:nondiv}
\subsection{Extension operators} \label{subsec:extop}
In this subsection, we construct some extension operators for Besov spaces $B^{r}_{p,q}(\R^{d-1})$ and $B^{r/2,r}_{p,q}(\R^d)$.

For $\alpha < 1, \tau > 0, x_d > 0$, and $z' = (z_1,\dots,z_{d-1}) \in \R^{d-1}$, denote
\begin{equation*}
Q^{\alpha}_{x_d}(\tau,z') = \frac{1}{4^{(d-\alpha)/2}\pi^{(d-1)/2}\Gamma((1-\alpha)/2)}\cdot \frac{x_d^{1-\alpha}}{\tau^{(d-\alpha)/2+1}} e^{-(x_d^2+|z'|^2)/(4\tau)}.
\end{equation*}
For $\alpha=0$, we write $Q_{x_d}$ instead of $Q_{x_d}^0$.
For $(t,x) \in \R^{d+1}_+$ and $g=g(t,x')$, denote
\begin{equation*}
[Q^{\alpha}_{x_d} \ast g](t,x') = \int_0^{\infty} \int_{\R^{d-1}} Q^{\alpha}_{x_d}(\tau,z') g(t-\tau,x'-z') dz' d\tau.
\end{equation*}
A direct computation shows that the kernel $Q^\alpha_{x_d}$ satisfies
\begin{gather*}
\partial_{\tau}Q^\alpha_{x_d} = \Delta_{z'} Q^\alpha_{x_d} + D_{x_d}^2 Q^\alpha_{x_d} + \frac{\alpha}{x_d} D_{x_d} Q^\alpha_{x_d} \quad \text{for $\tau,x_d > 0, z' \in \R^{d-1}$}, \\
\int_0^{\infty} \int_{\R^{d-1}} Q^{\alpha}_{x_d}(\tau,z') dz'd\tau = 1 \quad \text{for all $x_d > 0$},\\
\lim_{x_d \to 0^+} \int_{|z'|+\tau \geq \delta} Q^{\alpha}_{x_d}(\tau,z') dz'd\tau = 0 \quad \text{for any $\delta > 0$}.
\end{gather*}

The following proposition is an immediate corollary of the above properties.
\begin{proposition} \label{prop:extensionparabolic}
Let $\alpha \in (-\infty,1), p \in [1,\infty], g \in L_p(\R^d)$. Set $u(t,x) = [Q^{\alpha}_{x_d} \ast g](t,x')$. Then $u \in C^{\infty}(\R^{d+1}_+)$ satisfies 
\begin{equation*}
\begin{cases}
x_d^{\alpha}u_t = \sum\limits_{i=1}^d D_i(x_d^{\alpha} D_i u) & \text{in $\R^{d+1}_+$}, \\
\lim\limits_{x_d \to 0^+} u(t,x',x_d) = g(t,x') & \text{for a.e. $(t,x') \in \R^d$}.
\end{cases}
\end{equation*}
For $p \in [1,\infty)$, the convergence $\lim_{x_d \to 0^+} u(t,x',x_d) = g(t,x')$ holds in $L_p(\R^d)$.
\end{proposition}
See \cite{StingaTorrea2017SIAM} for a thorough study of relations and applications of the extension operator $Q^{\alpha}_{x_d} \ast g$ to the fractional heat operator $(\partial_t - \Delta_{x'})^{\frac{1-\alpha}{2}} g$ for $-1 < \alpha < 1$.

The following result was proved in \cite[Theorem 3]{JonesJMM1968}. It states, in particular, that the parabolic Besov norm of a function can be equivalently expressed in terms of integrability properties of derivatives of its caloric extension ($\alpha = 0$). We denote $L_q^* L_p (\R^{d+1}_+) = L_q(\R_+,\frac{dx_d}{x_d}; L_p(\R^{d}_{t, x'}))$.
\begin{thm} \label{thm:equivnormsparbesov}
Let $p,q \in [1,\infty], g \in L_p(\R^{d})$. Denote $u(t,x) = [Q_{x_d} \ast g](t,x')$. Let
\begin{gather*}
[g]_{B^{r/2,r}_{p,q}(\R^{d})}^{(1)} = \Vert x_d^{1-r} D u \Vert_{L_q^* L_p (\R^{d+1}_+)} + \Vert x_d^{2-r} \partial_t u \Vert_{L_q^* L_p(\R^{d+1}_+)}, \\
[g]_{B^{r/2,r}_{p,q}(\R^{d})}^{(2)} = \Vert x_d^{2-r} D^2 u \Vert_{L_q^* L_p (\R^{d+1}_+)} + \Vert x_d^{2-r} \partial_{t} u \Vert_{L_q^* L_p (\R^{d+1}_+)}, \\
[g]_{B^{r/2,r}_{p,q}(\R^{d})}' = \left\Vert \frac{\Vert g(t-\tau,x'-z')-g(t,x')\Vert_{L_p^{t,x'}(\R^d)}}{(|z'|+|\tau|^{1/2})^r} \right\Vert_{L_q(\R^d,\frac{dz' d\tau}{(|z'|+|\tau|^{1/2})^{d+1}})}.
\end{gather*}
For $r \in (0,1)$, the seminorms $[g]_{B^{r/2,r}_{p,q}}$, $[g]_{B^{r/2,r}_{p,q}}^{(1)}$, $[g]_{B^{r/2,r}_{p,q}}^{(2)}$, and $[g]_{B^{r/2,r}_{p,q}}'$ are equivalent. \\
For $r \in (0,2)$, the seminorms $[g]_{B^{r/2,r}_{p,q}}$ and $[g]_{B^{r/2,r}_{p,q}}^{(2)}$ are equivalent. The implied constants depend only on $d$ and $r$.
\end{thm}
The case $p=q=\infty$ with similar results was also considered in \cite{StingaTorrea2017SIAM}. The case $p=q=\infty$ is of a particular interest since for $r \in (0,\infty)$ the spaces $B^{r/2,r}_{\infty,\infty}(\R^d)$ coincide with the parabolic H\"older-Zygmund spaces $C^{r/2,r}(\R^d)$.
See Theorems 7.1 and 7.2 in \cite{StingaTorrea2017SIAM}.

We will need the following partial extension of Theorem \ref{thm:equivnormsparbesov} to the whole range $\alpha < 1$, which is of independent interest on its own.
\begin{thm} \label{thm:extbesovweightedpar}
Let $\alpha \in (-\infty, 1), p, q \in [1,\infty], g \in L_{p}(\R^d)$. Denote $u(t,x) = [Q^{\alpha}_{x_d} \ast g](t,x')$. For $0 < r < \min(1-\alpha,1)$, we have
\begin{equation*}
C_{d,\alpha, r}^{-1} [g]_{B^{r/2,r}_{p,q}(\R^d)} \leq \Vert x_d^{1-r} D u \Vert_{L_q^* L_p(\R^{d+1}_+)} + \Vert x_d^{2-r} \partial_t u \Vert_{L_q^* L_p(\R^{d+1}_+)} \leq C_{d,\alpha, r} [g]_{B^{r/2,r}_{p,q}(\R^d)}.
\end{equation*}

\end{thm}
\begin{proof}
The first inequality is a standard trace estimate and uses only $u(t,x',0)=g(t,x')$. See the argument on p. 394 in \cite{JonesJMM1968}.

To prove the second inequality, we show an equivalent bound
\begin{equation} \label{eq:besovextparupperbound}
\Vert x_d^{1-r} D u \Vert_{L_q^* L_p(\R^{d+1}_+)} + \Vert x_d^{2-r} \partial_t u \Vert_{L_q^* L_p(\R^{d+1}_+)} \leq C_{d,\alpha, r} [g]'_{B^{r/2,r}_{p,q}(\R^d)},
\end{equation}
where the seminorm $[g]'$ is defined in Theorem \ref{thm:equivnormsparbesov}.
Using
\begin{equation*}
\int_0^\infty \int_{\R^{d-1}}\partial_{\ell} Q^{\alpha}_{x_d}(\tau,z')dz'd\tau = 0,
\end{equation*}
where $\ell \in \{z_1,\dots, z_{d-1}, x_d, \tau\}$, we have
\begin{equation*}
\partial_{\ell}u(t,x) = \int_0^\infty \int_{\R^{d-1}} \partial_{\ell}Q^{\alpha}_{x_d}(\tau,z')[g(t-\tau,x'-z')-g(t,x')] dz'd\tau.
\end{equation*}
Let $\nabla = (\partial_{z_1}, \dots, \partial_{z_{d-1}}, \partial_{x_d})$. Applying Minkowski's inequality, we get
\begin{equation} \label{eq:extensionintermediatebound}
\begin{gathered}
x_d^{1-r} \Vert D u(\cdot,\cdot,x_d) \Vert_{L_p(\R^d)} + x_d^{2-r} \Vert \partial_t u(\cdot,\cdot,x_d) \Vert_{L_p(\R^d)} \\
\leq \int_0^\infty \int_{\R^{d-1}} L_{x_d}(\tau,z') \frac{\Vert g(t-\tau,x'-z')-g(t,x')\Vert_{L_p(\R^d)}}{(|z'|+\tau^{1/2})^r} dz' d\tau,
\end{gathered}
\end{equation}
where
\begin{equation*}
L_{x_d}(\tau,z') = \left\{x_d^{1-r}|\nabla Q^{\alpha}_{x_d}(\tau,z')|+x_d^{2-r}|\partial_{\tau} Q^{\alpha}_{x_d}(\tau,z')|\right\}\cdot(|z'|+\tau^{1/2})^r.
\end{equation*}
It is straightforward to verify the following derivative estimates of the kernel $Q^{\alpha}_{x_d}$:
\begin{equation} \label{eq:Qgradestimate}
|\nabla Q^{\alpha}_{x_d}| \leq C_{d,\alpha}\left(\frac{1}{x_d}+\frac{x_d+|z'|}{\tau}\right)Q^{\alpha}_{x_d}, \quad
|\partial_{\tau} Q^{\alpha}_{x_d}| \leq C_{d,\alpha} \left(\frac{1}{\tau}+\frac{x_d^2+|z'|^2}{\tau^2}\right)Q^{\alpha}_{x_d}.
\end{equation}
Using (\ref{eq:Qgradestimate}), the definition of $Q^{\alpha}_{x_d}$, and the condition $0 < r < \min(1-\alpha,1)$, one can verify that
\begin{equation*}
\int_0^\infty \int_{\R^{d-1}} L_{x_d}(\tau,z') dz' d\tau \leq C_{d,\alpha,r} < \infty
\end{equation*}
for all $x_d>0$. Applied to (\ref{eq:extensionintermediatebound}), this implies (\ref{eq:besovextparupperbound}) for $q=\infty$. One can also check that
\begin{equation*}
\int_0^{\infty} L_{x_d}(\tau,z') \frac{(|z'|+\tau^{1/2})^{d+1}}{x_d} dx_d \leq C_{d,\alpha,r} < \infty
\end{equation*}
for a.e. $(\tau,z') \in \R^d$. This implies (\ref{eq:besovextparupperbound}) for $q = 1$.
For $q \in (1,\infty)$, estimate (\ref{eq:besovextparupperbound}) now follows by interpolation.
\end{proof}

In the following theorem, we construct extension operators, which are suitable for equations in non-divergence form.
For $x_d > 0, z' \in \R^{d-1}$, let
\begin{equation} \label{eq:Poissonkernel}
P_{x_d} (z') = \int_0^{\infty} Q_{x_d} (\tau, z') d\tau = \frac{\Gamma(d/2)}{\pi^{d/2}}\cdot \frac{x_d}{(x_d^2+|z'|^2)^{d/2}}.
\end{equation}
Then for $g=g(x')$, the function
\begin{equation*}
[P_{x_d} \ast g](x') = \int_{\R^{d-1}} P_{x_d}(z')g(x'-z')dz'
\end{equation*}
is the harmonic extension of $g$. Let $\chi(y)$ be a fixed cut-off function with
\begin{equation} \label{eq:chicutoff}
\begin{gathered}
\chi(y) \in C_0^{\infty}(\mathbb{R}), \quad \text{$0 \leq \chi(y) \leq 1$ for all $y$,} \\
\text{$\chi(y) = 1$ for $|y| \leq 1$, \quad $\chi(y) = 0$ for $|y| \geq 2$.}
\end{gathered}
\end{equation}
\begin{proposition} \label{prop:extensionoperator}
    Let $p \in (1,\infty)$ and $\beta \in (-1,2p-1)$.
    \begin{enumerate}
    \item[(i)]
        Let $g=g(x')$. The operator $(E_{\operatorname{ell}}g)(x) = \chi(x_d)[(2P_{x_d}-P_{2x_d}) \ast g](x')$ satisfies
        \begin{gather*}
            E_{\operatorname{ell}}: B^{2-\frac{\beta+1}{p}}_{p,p}(\mathbb{R}^{d-1}) \to W^2_{p,\beta,\mathcal{N}}(\mathbb{R}^{d}_+), \\
        	E_{\operatorname{ell}}g(x',0) = g(x');
        \end{gather*}
    \item[(ii)]
        Let $g=g(t,x')$. The operator $(E_{\operatorname{par}}g)(t,x) = \chi(x_d)[(2Q_{x_d}-Q_{2x_d}) \ast g](t,x')$ satisfies
        \begin{gather}
        E_{\operatorname{par}}: B^{1-\frac{\beta+1}{2p}, 2-\frac{\beta+1}{p}}_{p,p}(\R^d) \to W^{1,2}_{p,\beta,\mathcal{N}}(\R^{d+1}_+), \label{eq:extopparnondivmapping} \\
    	E_{\operatorname{par}}g(t,x',0) = g(t,x'). \label{eq:extopparnondivbc}
        \end{gather}
    \end{enumerate}
\end{proposition}
\begin{proof}
We first give the proof for $E_{\operatorname{par}}$. Let $v(t,x) = [Q_{x_d} \ast g](t,x')$ denote the caloric extension of $g$. Denote $u(t,x) = (E_{\operatorname{par}} g)(t,x) = \chi(x_d) (2v(t,x',x_d)-v(t,x',2x_d))$. The boundary condition (\ref{eq:extopparnondivbc}) holds by Proposition \ref{prop:extensionparabolic}. We now show that
\begin{equation} \label{eq:extopparnondivmappingweak}
E_{\operatorname{par}}: B^{1-\frac{\beta+1}{2p}, 2-\frac{\beta+1}{p}}_{p,p}(\mathbb{R}^{d}) \to W^{1,2}_{p,\beta}(\mathbb{R}^{d+1}_+).
\end{equation}
The proof of (\ref{eq:extopparnondivmappingweak}) is standard, but we give it here for completeness. Throughout the proof, $C=C_{d,p,\beta}$ unless otherwise stated. From Theorem \ref{thm:equivnormsparbesov} with $q=p$ and $r=2-(\beta+1)/p$, it follows that
\begin{equation*}
\Vert D^2 v \Vert_{L_{p,\beta}(\R^{d+1}_+)} + \Vert \partial_t v \Vert_{L_{p,\beta}(\R^{d+1}_+)} \leq C_{d,p, \beta} [g]_{B^{1-\frac{\beta+1}{2p}, 2-\frac{\beta+1}{p}}_{p,p}(\R^d)}.
\end{equation*}
A direct computation shows that
\begin{equation*}
\int_0^{\infty} \int_{\R^{d-1}} Q_{x_d}(\tau,z') dz'd\tau = 1, \quad \int_0^{\infty} \int_{\R^{d-1}} |\nabla Q_{x_d}(\tau,z')| dz'd\tau \leq C_d x_d^{-1}
\end{equation*}
for all $x_d > 0$. Therefore, by Young's inequality for all $x_d > 0$,
\begin{equation*}
\Vert v(\cdot,\cdot,x_d) \Vert_{L_p(\R^d)} \leq \Vert g \Vert_{L_p(\R^d)}, \quad \Vert D v(\cdot,\cdot,x_d) \Vert_{L_p(\R^d)} \leq C_d x_d^{-1} \Vert g \Vert_{L_p(\R^d)}.
\end{equation*}
Combining the above estimates with $u(t,x) = \chi(x_d) (2v(t,x',x_d)-v(t,x',2x_d))$ gives
\begin{gather*}
\Vert u \Vert_{L_{p,\beta}(\R^{d+1}_+)} \leq C_{p,\beta} \Vert g \Vert_{L_{p}(\R^d)}, \\
\Vert D^2 u \Vert_{L_{p,\beta}(\R^{d+1}_+)} + \Vert \partial_t u \Vert_{L_{p,\beta}(\R^{d+1}_+)} \leq C \Vert g \Vert_{B^{1-\frac{\beta+1}{2p}, 2-\frac{\beta+1}{p}}_{p,p}(\R^d)}.
\end{gather*}
The term $Du$ is estimated by Hardy's inequality at infinity; see Lemma \ref{lem:Hardy} (i):
\begin{equation*}
\Vert Du \Vert_{L_{p,\beta}(\R^{d+1}_+)} \leq C \Vert D_d Du \Vert_{L_{p,\beta+p}(\R^{d+1}_+)} \leq C \Vert D_d Du \Vert_{L_{p,\beta}(\R^{d+1}_+)} \leq C \Vert g \Vert_{B^{1-\frac{\beta+1}{2p}, 2-\frac{\beta+1}{p}}_{p,p}(\R^d)},
\end{equation*}
where the middle inequality holds since $\chi(y) = 0$ for $y \geq 2$. From the above estimates, (\ref{eq:extopparnondivmappingweak}) follows.

We now show (\ref{eq:extopparnondivmapping}). For $\beta \in (-1,p-1)$, $D_d v(t,x',0) \in B^{\frac{1}{2}-\frac{\beta+1}{2p}, 1-\frac{\beta+1}{p}}_{p,p}(\R^d)$ is well defined. Therefore, $D_d u(t,x',0)=0$ by the construction of $u$. Now, (\ref{eq:extopparnondivmapping}) follows by Hardy's inequality at zero, see Lemma \ref{lem:Hardy} (ii). For $\beta\in (p-1,2p-1)$, (\ref{eq:extopparnondivmapping}) holds since $W^{1,2}_{p,\beta,\mathcal{N}}(\mathbb{R}^{d+1}_+) = W^{1,2}_{p,\beta}(\mathbb{R}^{d+1}_+)$ by Lemma~\ref{lem:Hardy} (i) (i.e., $\Vert D_d u/x_d \Vert_{L_{p,\beta}(\R^{d+1}_+)} \leq \Vert D^2_d u \Vert_{L_{p,\beta}(\R^{d+1}_+)}$). For $\beta = p-1$, we use interpolation. Fix any $-1 < \beta_0 < p-1 < \beta_1 < 2p-1$ such that $p-1=(\beta_0+\beta_1)/2$. Then
    \begin{gather*}
        E_{\operatorname{par}}: B^{1/2,1}_{p,p}(\mathbb{R}^{d}) = [B^{1-\frac{\beta_0+1}{2p}, 2-\frac{\beta_0+1}{p}}_{p,p}(\mathbb{R}^{d}), B^{1-\frac{\beta_1+1}{2p}, 2-\frac{\beta_1+1}{p}}_{p,p}(\mathbb{R}^{d})]_{\frac{1}{2}} \to \\
        \to [W^{1,2}_{p,\beta_0,\mathcal{N}}(\mathbb{R}^{d+1}_+),W^{1,2}_{p,\beta_1,\mathcal{N}}(\mathbb{R}^{d+1}_+)]_{\frac{1}{2}} \hookrightarrow W^{1,2}_{p,p-1,\mathcal{N}}(\mathbb{R}^{d+1}_+),
    \end{gather*}
    where the last embedding is a simple consequence of the Stein-Weiss theorem on complex interpolation of weighted Lebesgue spaces. See a similar argument in \cite[Theorem 3.3]{CwikelEinavJFA2019}.

The argument for $E_{\operatorname{ell}}$ is the same, except that an elliptic version of Theorem \ref{thm:equivnormsparbesov} needs to be used. See, for example, \cite[Theorem 4]{Taibleson1964Indiana}.
\end{proof}
\subsection{Non-divergence form equations} \label{subsec:nondiv}
In this subsection, we prove some auxiliary results concerning the boundary value problem for equations in non-divergence form. Following the discussion in Section \ref{subsec:varchange}, throughout this subsection we assume that
\begin{equation} \label{eq:structuralcond}
a_{dj} = 0 \quad \text{and} \quad a_{dd}=1, \quad j = 1, \dots, d-1.
\end{equation}
For compactness of notation, we denote
\begin{gather}
\label{eq:operatornondivell}
\mathcal{E}_{\text{non-div}} = -D^2_d-\frac{\alpha}{x_d}D_d - \sum_{i=1}^{d-1}\sum_{j=1}^d a_{ij}(x_d)D_i D_j + \lambda c(x_d), \\
\label{eq:operatornondivpar}
\mathcal{P}_{\text{non-div}} = a_0(x_d)\partial_t-D^2_d-\frac{\alpha}{x_d}D_d - \sum_{i=1}^{d-1}\sum_{j=1}^d a_{ij}(x_d)D_i D_j + \lambda c(x_d).
\end{gather}
We consider the parabolic boundary value problem
\begin{equation} \label{eq:nondivpar}
\begin{cases}
\mathcal{P}_{\text{non-div}}u = f & \text{in $\R^{d+1}_+$}, \\
   u(t,x',0) = U(t,x') & \text{on $\R^d$}, 
\end{cases}
\end{equation}
and, in the elliptic case,
\begin{equation} \label{eq:nondivell}
\begin{cases}
\mathcal{E}_{\text{non-div}}u = f & \text{in $\R^{d}_+$}, \\
   u(x',0) = U(x') & \text{on $\R^{d-1}$},
\end{cases}
\end{equation}

	Let $f \in L_{1, \text{loc}}(\R^{d+1}_+)$ and $U \in \mathcal{D}'(\R^d)$. We say that $u \in W^{1,2}_{1,\text{loc}}(\R^{d+1}_+)$ is a strong solution to (\ref{eq:nondivpar}) if the first equation in (\ref{eq:nondivpar}) holds a.e. in $\R^{d+1}_+$ and the second equation holds in the sense of Definition \ref{def:trace}.

Let $f \in L_{1, \text{loc}}(\R^{d}_+)$ and $U \in \mathcal{D}'(\R^{d-1})$. We say that $u \in W^2_{1,\text{loc}}(\mathbb{R}^d_+)$ is a strong solution to (\ref{eq:nondivell}) if the first equation in (\ref{eq:nondivell}) holds a.e. in $\mathbb{R}^d_+$ and the second equation holds in the sense of Definition~\ref{def:trace}.

As a corollary of Proposition \ref{prop:extensionoperator}, we get the following result on the solvability of inhomogeneous non-divergence form equations (\ref{eq:nondivpar}) in $W^{1,2}_{p,\alpha,\beta}(\R^{d+1}_+)$.
\begin{thm} \label{thm:solvnondivpar}
    Let $\alpha \in (-\infty,1), p \in (1,\infty), \beta \in (\alpha^+ p-1, 2p-1)$, and $\kappa \in (0,1]$. Suppose that \eqref{eq:ellipticitycond} and \eqref{eq:structuralcond} are satisfied. Then, for any $\lambda > 0, f \in L_{p,\beta}(\R^{d+1}_+)$, and $U \in B^{1-\frac{\beta+1}{2p},2-\frac{\beta+1}{p}}_{p,p}(\R^d)$, there exists a unique strong solution $u \in W^{1,2}_{p,\alpha,\beta}(\R^{d+1}_+)$ to \eqref{eq:nondivpar}. Moreover, $u$ satisfies
\begin{multline} \label{eq:estnondivpar}
    \Vert u_t \Vert + \Vert D_{x'} Du\Vert + \Vert D^2_d u + \alpha x_d^{-1} D_d u\Vert + \lambda^{1/2} \Vert Du \Vert + \lambda \Vert u \Vert \leq \\
    \leq N \Vert f \Vert + N(\lambda^{1-\frac{\beta+1}{2p}}\Vert U \Vert_{L_{p}(\R^d)}+[U]),
\end{multline}
where $\Vert \cdot \Vert = \Vert \cdot \Vert_{L_{p,\beta}(\R^{d+1}_+)}$, $[U] = [U]_{B^{1-\frac{\beta+1}{2p},2-\frac{\beta+1}{p}}_{p,p}(\R^d)}$, and $N = N(d,\alpha,p,\beta,\kappa) > 0$.
\end{thm}
\begin{proof}
For $U=0$ and $\beta \in (\alpha p-1,2p-1)$, Theorem \ref{thm:solvnondivpar} follows from \cite[Theorem 2.3]{DongPhan2023JFA}.

For $U \neq 0$, by considering the extension operator from Proposition~\ref{prop:extensionoperator}, we obtain the existence and uniqueness of solutions for all $\lambda>0$ and estimate (\ref{eq:estnondivpar}) for $\lambda=1$. For general $\lambda>0$, estimate (\ref{eq:estnondivpar}) follows by scaling. For completeness, we present the argument here. Let $\mu = \lambda^{\frac{1}{2}}, (a_0^{\mu},a_{ij}^{\mu},c^{\mu})(x_d) = (a_0,a_{ij},c)(\mu^{-1}x_d), U^{\mu}(t,x')=U(\mu^{-2}t,\mu^{-1}x')$, and $f^{\mu}(t,x)=f(\mu^{-2}t,\mu^{-1}x)$.

Let $u^{\mu} \in W^{1,2}_{p,\alpha,\beta}(\R^{d+1}_+)$ be the strong solution to
\begin{gather*}
    a_0^{\mu}(x_d)u^{\mu}_t-D^2_d u^{\mu}-\frac{\alpha}{x_d}D_d u^{\mu} - \sum_{i=1}^{d-1}\sum_{j=1}^d a_{ij}^{\mu}(x_d)D_i D_j u^{\mu} + c^{\mu}(x_d) u^{\mu} 
    = \mu^{-2}f, \\
    u^{\mu}(t,x',0) = U^{\mu}(t,x').
\end{gather*}

The function $u^{\mu}$ satisfies
\begin{gather*}
    \Vert u^{\mu}_t \Vert + \Vert D_{x'} D u^{\mu} \Vert + \Vert D^2_d u^{\mu} +\alpha x_d^{-1}D_d u^{\mu} \Vert + \Vert D u^{\mu} \Vert + \Vert u^{\mu} \Vert \leq N\mu^{-2}\Vert f^{\mu}\Vert \\
    + N\bigl(\Vert U^{\mu} \Vert_{L_p(\R^{d})} + [U^{\mu}]\bigr).
\end{gather*}
Rewriting this estimate in terms of $u(t,x) = u^\mu(\mu^2 t, \mu x), f$, and $U$, gives  (\ref{eq:estnondivpar}).
\end{proof}
The following is an elliptic counterpart of Theorem \ref{thm:solvnondivpar}. 
\begin{thm} \label{thm:solvnondivell}
    Let $\alpha \in (-\infty,1), p \in (1,\infty), \beta \in (\alpha^+ p-1, 2p-1)$, and $\kappa \in (0,1]$. Suppose that \eqref{eq:ellipticitycond} and \eqref{eq:structuralcond} are satisfied. Then, for any $\lambda > 0, f \in L_{p,\beta}(\mathbb{R}^d_+)$, and $U \in B^{2-\frac{\beta+1}{p}}_{p,p}(\mathbb{R}^{d-1})$, there exists a unique strong solution $u \in W^2_{p,\alpha,\beta}(\mathbb{R}^{d}_+)$ to \eqref{eq:nondivell}. Moreover, $u$ satisfies
\begin{equation*}
    \Vert D_{x'} Du\Vert + \Vert D^2_d u + \alpha x_d^{-1} D_d u\Vert + \lambda^{1/2} \Vert Du \Vert + \lambda \Vert u \Vert
    \leq N \Vert f \Vert + N(\lambda^{1-\frac{\beta+1}{2p}}\Vert U \Vert_{L_{p}(\mathbb{R}^{d-1})}+[U]),
\end{equation*}
    where $\Vert \cdot \Vert = \Vert \cdot \Vert_{L_{p,\beta}(\mathbb{R}^d_+)}$, $[U] = [U]_{B^{2-\frac{\beta+1}{p}}_{p,p}(\mathbb{R}^{d-1})}$, and $N = N(d,\alpha,p,\beta,\kappa) > 0$.
\end{thm}

The proof is similar to the proof of Theorem \ref{thm:solvnondivpar}, except that \cite[Theorem 2.4]{DongPhan2023JFA} should be used instead of \cite[Theorem 2.3]{DongPhan2023JFA}.

Through the end of this section, we fix $\alpha \in (-\infty,1)$, $p \in (1, \infty)$, $\kappa \in (0,1]$, $\lambda>0$, and the coefficients $a_0(x_d),a_{ij}(x_d),c(x_d)$ satisfying (\ref{eq:ellipticitycond}) and (\ref{eq:structuralcond}). For $\beta \in (\alpha^+ p-1, 2p-1)$, let
\begin{equation} \label{eq:solop}
\begin{gathered}
    S:B^{1-\frac{\beta+1}{2p},2-\frac{\beta+1}{p}}_{p,p}(\mathbb{R}^d) \to W^{1,2}_{p,\alpha,\beta}(\mathbb{R}^{d+1}_+), \\
    S:U \mapsto u
\end{gathered}
\end{equation}
be the solution operator from Theorem \ref{thm:solvnondivpar} to the boundary value problem
\begin{equation} \label{eq:nondivzerorhs}
\begin{cases}
    \mathcal{P}_{\text{non-div}} u = 0 & \text{in $\mathbb{R}^{d+1}_+$}, \\
    u(t,x',0) = U(t,x') & \text{on $\mathbb{R}^d$}.
\end{cases}
\end{equation}

Formally, the operator $S$ depends on both $p$ and $\beta$. The next proposition shows that the operator $S$ is well defined on the union of trace spaces $\bigcup_{\beta \in (\alpha^+ p-1, 2p-1)} B^{1-\frac{\beta+1}{2p},2-\frac{\beta+1}{p}}_{p,p}(\mathbb{R}^d)$.

\begin{lem} \label{lem:solutionoperatorwelldefined}
Let $p_i \in (1,\infty), \beta_i \in (\alpha^+ p_i-1, 2p_i-1)$ and suppose that $U \in B^{1-\frac{\beta_i+1}{2p_i},2-\frac{\beta_i+1}{p_i}}_{p_i,p_i}(\mathbb{R}^d)$, where $i=1,2$. Let $u_i \in W^{1,2}_{p_i,\alpha,\beta_i}(\R^{d+1}_+)$ be the strong solutions to \eqref{eq:nondivzerorhs}. Then $u_1=u_2$.
\end{lem}
The proof of Lemma \ref{lem:solutionoperatorwelldefined} is given in Appendix \ref{appendix}.

For $\alpha \in (-\infty, 1)$ and $p \in (1,\infty)$ we define an open interval
\begin{equation} \label{eq:rangeIp}
I_p(\alpha)=
    \begin{cases}
        (-1,2p-1) & \text{for $\alpha \in (-\infty,-1]\cup\{0\}$}, \\
        ((1+\alpha)p-1,2p-1) & \text{for $\alpha \in (-1,1)\setminus \{0\}$}.
    \end{cases}
\end{equation}
The following lemma is one of the main tools in the proof of the existence of solutions in Theorem~\ref{thm:tracepar}.
\begin{lem} \label{lem:embedsoltnop}
    For any $\beta \in I_p(\alpha)$ holds
    \begin{equation} \label{eq:mapbesovtomixed}
        S:B^{1-\frac{\beta+1}{2p},2-\frac{\beta+1}{p}}_{p,p}(\mathbb{R}^d) \to H^{1/2,1}_p(\mathbb{R}^d; W^1_{p,\beta}(\mathbb{R}_+)).
    \end{equation}
    The operator norm of $S$ in \eqref{eq:mapbesovtomixed} is bounded by $N=N(d,\alpha,p,\beta,\kappa,\lambda)>0$.
\end{lem}

\begin{proof}
    By Lemma \ref{lem:Hardy}, for any $\beta \in I_p(\alpha)$, we have $W^{1,2}_{p,\alpha,\beta}(\mathbb{R}^{d+1}_+) \hookrightarrow W^{1,2}_{p,\beta}(\mathbb{R}^{d+1}_+).$
    Note that
    \begin{gather*}
        W^{1,2}_{p,\beta}(\mathbb{R}^{d+1}_+) \hookrightarrow W^{1,2}_p(\mathbb{R}^d; L_{p,\beta}(\mathbb{R}_+)) \cap L_p(\mathbb{R}^d; W^{2}_{p,\beta}(\mathbb{R}_+)) \\
        \hookrightarrow [W^{1,2}_p(\mathbb{R}^d; L_{p,\beta}(\mathbb{R}_+)), L_p(\mathbb{R}^d; W^{2}_{p,\beta}(\mathbb{R}_+))]_{\frac{1}{2}} = [\text{by \cite[Theorem 4.5.5]{amann2019linear}}] \\
        = H^{1/2,1}_p(\mathbb{R}^d; [L_{p,\beta}(\mathbb{R}_+), W^{2}_{p,\beta}(\mathbb{R}_+)]_{\frac{1}{2}}).
    \end{gather*}
    In \cite[Proposition 3.17]{LindemulderVeraar2020JDE}, it was shown that for $\beta \in (-p-1,2p-1)\setminus \{-1,p-1\}$, we have
    \begin{equation*}
        [L_{p,\beta}(\mathbb{R}_+), W^{2}_{p,\beta}(\mathbb{R}_+)]_{\frac{1}{2}} = W^1_{p,\beta}(\mathbb{R}_+).
    \end{equation*}
    From (\ref{eq:solop}) and the argument above, (\ref{eq:mapbesovtomixed}) holds for $\beta \in I_p(\alpha)\setminus \{p-1\}$. For $\alpha \in (0,1)$, we are done.
    
    For $\alpha \in (-\infty,0]$ and $\beta = p-1$, we use interpolation. Fix any $\beta_0, \beta_1 \in I(\alpha)$ such that $\beta_0 < p-1 <\beta_1$ and $p-1 = (\beta_0+\beta_1)/2$. Taking the $[\;,\,]_{\frac{1}{2}}$ interpolation of (\ref{eq:mapbesovtomixed}) with $\beta_0$ and $\beta_1$, we get
    \begin{gather*}
        S:B^{1/2,1}_{p,p}(\mathbb{R}^d) = [B^{1-\frac{\beta_0+1}{2p},2-\frac{\beta_0+1}{p}}_{p,p}(\mathbb{R}^d),B^{1-\frac{\beta_1+1}{2p},2-\frac{\beta_1+1}{p}}_{p,p}(\mathbb{R}^d)]_{\frac{1}{2}} \\
        \to [H^{1/2,1}_p(\mathbb{R}^d; W^1_{p,\beta_0}(\mathbb{R}_+)),H^{1/2,1}_p(\mathbb{R}^d; W^1_{p,\beta_1}(\mathbb{R}_+))]_{\frac{1}{2}} = [\text{by Lemma \ref{lem:interpbessel}}]\\
        = H^{1/2,1}_p(\mathbb{R}^d;[W^1_{p,\beta_0}(\mathbb{R}_+),W^1_{p,\beta_1}(\mathbb{R}_+)]_{\frac{1}{2}}) \hookrightarrow H^{1/2,1}_p(\mathbb{R}^d; W^1_{p,p-1}(\mathbb{R}_+)).
    \end{gather*}
    The bound
    \begin{equation*}
        \Vert S \Vert_{B^{1-\frac{\beta+1}{2p},2-\frac{\beta+1}{p}}_{p,p}(\mathbb{R}^d) \to H^{1/2,1}_p(\mathbb{R}^d; W^1_{p,\beta}(\mathbb{R}_+))} \leq N(d,\alpha,p,\beta,\kappa,\lambda)
    \end{equation*}
    holds because the norms of all inclusions above depend only on the parameters listed in $N$.
\end{proof}

\section{Proofs of main results: elliptic equations} \label{sec:proofsell}
\subsection{Solvability of equations with \texorpdfstring{$F_d=0$}{mathformtext}}
\begin{proof}[Proof of Theorem {\rm \ref{thm:solvell}}] Without loss of generality, assume that $a_{dj}=0$ for $1 \leq j \leq d-1$ and $a_{dd}=1$. The general case follows from the change of variables \eqref{eq:coordchange}.

For compactness of notation, let $\mathcal{E}_{\text{non-div}}$ be as in (\ref{eq:operatornondivell}), and 
\begin{equation*}
\mathcal{E}_{\text{div}} = -D^2_d-\frac{\alpha}{x_d}D_d - \sum_{i=1}^{d-1}\sum_{j=1}^d D_i(a_{ij}(x_d)D_j\;) + \lambda c(x_d).
\end{equation*}

We first show the uniqueness of solutions. From the linearity of the problem, we need to show that if $u,Du \in L_{p,\beta}(\R^d_+)$ and $\mathcal{E}_{\text{div}}u = 0$ with $u(x',0)=0$, then $u=0$. For $ \varepsilon>0$, let
\begin{equation*}
u^{(\varepsilon)}(x',x_d) = u \ast \eta^{(\varepsilon)} = \int_{\R^{d-1}} u(x'-z',x_d) \eta^{(\varepsilon)}(z') dz'
\end{equation*}
be the mollification of $u$ in the $x'$ variables. Here $\eta^{(\varepsilon)}(x')=\varepsilon^{-d+1}\eta(x'/\varepsilon)$, where $\eta(x')$ is a fixed function such that
\[
\eta \in C_0^{\infty}(\R^{d-1}), \quad \int \eta = 1, \quad \eta \geq 0.
\]
Since the coefficients $a_{ij}$ and $c$ depend only on $x_d$, $u^{(\varepsilon)}$ satisfies $\mathcal{E}_{\text{div}}u^{(\varepsilon)} = 0$ and $u^{(\varepsilon)}(x',0)=0$. From the properties of mollification, we have
\begin{equation} \label{eq:integrabilitymollifell}
u^{(\varepsilon)}, Du^{(\varepsilon)}, D_{x'} Du^{(\varepsilon)} \in L_{p,\beta}(\R^d_+).
\end{equation}
From (\ref{eq:integrabilitymollifell}) and the equation $\mathcal{E}_{\text{div}}u^{(\varepsilon)} = 0$, it follows that $u^{(\varepsilon)}$ is a strong solution to $\mathcal{E}_{\text{non-div}}u^{(\varepsilon)} = 0$ with $u^{(\varepsilon)}(x',0)=0$. From Theorem \ref{thm:solvnondivell}, we have $u^{(\varepsilon)} = 0$. Since $\varepsilon > 0$ is arbitrary, $u=0$.

We now show the existence of solutions for $\lambda=1$. As $F_d=0$ and $a_{dj}=\delta_{dj}$, equations (\ref{equationnonhomweighted})-(\ref{dirichlet}) after multiplication by $x_d^{-\alpha}$ are given by 
\begin{equation} \label{eq:ellipticnonhom}
\begin{cases}
\mathcal{E}_{\text{div}}u = - \sum_{i=1}^{d-1}D_i F_i+f & \text{in $\R^{d}_+$,} \\
u(x',0)=U(x') & \text{on $\R^{d-1}$}.
\end{cases}
\end{equation}
To construct a solution to (\ref{eq:ellipticnonhom}), we differentiate solutions to non-divergence form equation~(\ref{eq:nondivell}).

Let $W=(1-\Delta_{x'})^{-\frac{1}{2}}U \in B^{2-\frac{\beta+1}{p}}_{p,p}(\R^{d-1})$. Using Theorem \ref{thm:solvnondivell}, we find unique strong solutions $u_0,u_1,\dots,u_{d-1},w \in W^2_{p,\alpha,\beta}(\R^d_+)$ to
\begin{equation*}
\begin{cases}
\mathcal{E}_{\text{non-div}} w = 0, \\
w(x',0) = W(x'),
\end{cases} \quad
\begin{cases}
\mathcal{E}_{\text{non-div}} u_0 = f, \\
u_i(x',0)=0,
\end{cases} \quad
\begin{cases}
\mathcal{E}_{\text{non-div}} u_i = F_i, \\
u_i(x',0)=0,
\end{cases}
\end{equation*}
where $i = 1,\dots,d-1$. Then
\begin{equation*}
u = -\sum_{i=1}^{d-1}D_i u_i + u_0 + (1-\Delta_{x'})^{\frac{1}{2}} w \in W^1_{p,\beta}(\R^d_+)
\end{equation*}
is a weak solution to (\ref{eq:ellipticnonhom}) satisfying (\ref{eq:estsolell}).

For general $\lambda > 0$, the existence of solutions and estimate (\ref{eq:estsolell}) follow by scaling.
\end{proof}

We also state a corresponding result for the non-degenerate case $\beta \in (\alpha^+ p-1, p-1)$.

\begin{proposition} \label{prop:solvabilityellnondeg}
Let $\alpha \in (-\infty,1), p \in (1,\infty), \beta \in (\alpha^+ p-1, p-1)$, and $\kappa \in (0,1]$. Suppose that \eqref{eq:ellipticitycond} is satisfied. Then, for any $\lambda > 0, F\in  L_{p,\beta}(\mathbb{R}^d_+; \mathbb{R}^d), f \in  L_{p,\beta}(\mathbb{R}^d_+)$, and $U \in B^{1-\frac{\beta+1}{p}}_{p,p}(\mathbb{R}^{d-1})$, there exists a unique weak solution $u \in W^1_{p,\beta}(\mathbb{R}^d_+)$ to \eqref{equationnonhomweighted}--\eqref{dirichlet}. Moreover, $u$ satisfies
\begin{equation} \label{eq:estimatesolutionellipticdivergenceA_p}
    \Vert Du \Vert + \sqrt{\lambda} \Vert u \Vert 
    \leq N \Vert F \Vert + \frac{N}{\sqrt{\lambda}} \Vert f \Vert + N(\lambda^{\frac{1}{2}-\frac{\beta+1}{2p}}\Vert U \Vert_{L_{p}(\mathbb{R}^{d-1})}+[U]),
\end{equation}
where $\Vert \cdot \Vert = \Vert \cdot \Vert_{L_{p,\beta}(\mathbb{R}^{d}_+)}, [U] = [U]_{B^{1-\frac{\beta+1}{p}}_{p,p}(\mathbb{R}^{d-1})}$, and $N=N(d,\alpha,p,\beta, \kappa) > 0$.
\end{proposition}
For $U=0$ and $\beta \in (\alpha p-1, p-1)$, Proposition \ref{prop:solvabilityellnondeg} follows from \cite[Theorem 2.8]{DongPhan2021TAMS} (also see Remark \ref{rmk:refremark} below for clarification). The case $U \neq 0$ reduces to a problem with $U=0$ by considering an extension operator $B^{1-\frac{\beta+1}{p}}_{p,p}(\mathbb{R}^{d-1}) \to W^1_{p,\beta}(\R^d_+)$. For example, $U(x') \mapsto \chi(x_d) [P_{x_d} \ast U](x')$, where $P_{x_d}$ is the Poisson kernel (\ref{eq:Poissonkernel}) and $\chi$ is from (\ref{eq:chicutoff}). When estimate (\ref{eq:estimatesolutionellipticdivergenceA_p}) is available for $\lambda = 1$, the case of general $\lambda > 0$ follows by scaling.

\subsection{Existence of trace}

Let $p \in (1,\infty)$ be fixed. For a weight parameter $\beta \in \R$, we define its dual exponent $\gamma$ (with respect to $p$) by
\begin{equation} \label{eq:defgamma}
\frac{\beta}{p}+\frac{\gamma}{p'}=1.
\end{equation}
Equivalently, $\gamma = \frac{p-\beta}{p-1}$. Note that
\begin{equation*} 
\beta \in [p-1, 2p-1) \iff \gamma \in (-1, p'-1].
\end{equation*}
The following lemma is a corollary of H\"older's inequality.
\begin{lem} \label{lem:holderweighted}
	Let $m \in \mathbb{N}$, $p \in (1,\infty)$, $\beta, \gamma \in \mathbb{R}$ satisfying $\frac{\beta}{p}+\frac{\gamma}{p'}=1$. Suppose that $u \in L_{p,\beta}(\mathbb{R}^m_+)$ and $v \in L_{p',\gamma}(\mathbb{R}^m_+)$. For $y>0$, define
    \[
    	z(y) =  \langle u(\cdot,y), v(\cdot,y) \rangle = \int_{\mathbb{R}^{m-1}}u(\cdot,y)v(\cdot,y).
    \]
    Then $yz(y) \in L_1(0,\infty)$ with
    \begin{gather*}
    	\Vert yz(y) \Vert_{L_1(0,\infty)} \leq \Vert u \Vert_{L_{p,\beta}(\mathbb{R}^m_+)} \Vert v \Vert_{L_{p',\gamma}(\mathbb{R}^m_+)}.
    \end{gather*}
\end{lem}

Furthermore, we use the notation $(\tau_h u)(\cdot, x_d) = u(\cdot, x_d+h)$. The following lemma states that shift operators $\tau_h$ are continuous in $L_{p,\beta}$ for $\beta \geq 0$.
\begin{lem} \label{lem:contshiftslp}
Let $p \in [1,\infty), \beta \geq 0$. Let $L_{p,\beta} = L_{p,\beta}(\R^d_+)$ or $L_{p,\beta} = L_{p,\beta}(\R^{d+1}_+)$. Then for any $f \in L_{p,\beta}$,
\begin{equation*}
\Vert \tau_h f \Vert_{L_{p,\beta}} \leq \Vert f \Vert_{L_{p,\beta}} \quad \text{for all $h > 0$}, \quad \lim_{h \to 0^+} \Vert \tau_h f - f \Vert_{L_{p,\beta}} \to 0^+.
\end{equation*}
\end{lem}
The proof follows by a standard argument using the density of $C_0$ in $L_{p,\beta}$.

The following lemma, which is of independent interest, is a key tool for proving the existence of trace in Theorems \ref{thm:traceell} and \ref{thm:tracepar}.
\begin{lem} \label{lem:1dtrace}
    Let $\alpha \in (-\infty, 1)$ and let $h \in [0,\infty)$ be fixed. Suppose that a function $w \in W^2_{1,\operatorname{loc}}(0,\infty)$ satisfies $\liminf_{y \to +\infty} |w(y)| = 0$ and
    \begin{equation} \label{eq:lemmaL1condition}
    y(w''(y)+\frac{\alpha}{h+y}w'(y)) \in L_1(0,\infty).
    \end{equation}
    Then, for all $y > 0$,
    \begin{equation} \label{eq:solution1dlemma}
    w(y) = \frac{1}{1-\alpha}\int_y^{\infty} \bigl((h+s)^{1-\alpha}-(h+y)^{1-\alpha}\bigr)(h+s)^{\alpha}(w''(s)+\frac{\alpha}{h+s}w'(s)) ds.
    \end{equation}
    Moreover,
    \begin{equation} \label{eq:formula1dlemma}
        w(0) := \lim_{y \to 0^+}w(y) = \begin{cases}
            \frac{1}{1-\alpha}\int_0^{\infty}K_{\alpha}(s/h)s(w''(s)+\frac{\alpha}{h+s}w'(s))ds, & h > 0, \\
            \frac{1}{1-\alpha}\int_0^\infty s(w''(s)+\frac{\alpha}{s}w'(s))ds, & h=0,
        \end{cases}
    \end{equation}
    where
    \[
        K_{\alpha}(s) := 1 - \frac{(1+s)^\alpha-1}{s}.
    \]
    Finally, the following estimate holds:
    \begin{equation} \label{eq:estimate1dlemma}
        |w(0)| \leq C_{\alpha} \int_0^{\infty} s\Bigl|w''(s)+\frac{\alpha}{h+s}w'(s)\Bigr|ds,
    \end{equation}
    where
    \begin{equation*}
        C_{\alpha} =
        \begin{cases}
            \max(1,(1-\alpha)^{-1}), & h > 0, \\
            (1-\alpha)^{-1}, & h = 0.
        \end{cases}
    \end{equation*}
\end{lem}

\begin{proof}
    Let $w$ be as in the statement of the lemma. Define a function $z(y)$ by
    \[
        z(y) = \frac{1}{1-\alpha}\int_y^{\infty} \bigl((h+s)^{1-\alpha}-(h+y)^{1-\alpha}\bigr)(h+s)^{\alpha}(w''(s)+\frac{\alpha}{h+s}w'(s)) ds.
    \]
    From (\ref{eq:lemmaL1condition}), it follows that $|z(y)| \to 0$ as $y \to +\infty$. By differentiating the integral, we get that for all $y>0$,
    \begin{equation*}
        z'(y) = -\int_y^{\infty}(h+y)^{-\alpha}(h+s)^{\alpha}(w''(s)+\frac{\alpha}{h+s}w'(s)) ds.
    \end{equation*}
    Differentiating once again gives that for a.e. $y>0$,
    \begin{equation*}
        z''(y) = w''(y)+\frac{\alpha}{h+y}w'(y) +\alpha \int_y^{\infty}(h+y)^{-\alpha-1}(h+s)^{\alpha}(w''(s)+\frac{\alpha}{h+s}w'(s)) ds.
    \end{equation*}
    That is, for a.e. $y>0$,
    \begin{equation*}
        z''(y)+\frac{\alpha}{h+y}z'(y) = w''(y)+\frac{\alpha}{h+y}w'(y).
    \end{equation*}
    Then the function $z-w$ satisfies
    \begin{equation*}
        ((h+y)^{\alpha}(z-w)')' = 0,
    \end{equation*}
    from which it follows that
    \begin{equation*}
        z(y)-w(y) = A(h+y)^{1-\alpha}+B
    \end{equation*}
    for some constants $A,B$. From
    \[
    \liminf_{y \to +\infty} |w(y)| = 0 \qquad \text{and} \qquad \lim_{y \to +\infty} |z(y)| = 0
    \]
    it follows that $A=B=0$ and hence $w(y)=z(y)$. This proves (\ref{eq:solution1dlemma}).

	By taking the limit as $y \to 0^+$, we see that (\ref{eq:solution1dlemma}) holds for $y=0$ as well. This proves (\ref{eq:formula1dlemma}).
	
The bound (\ref{eq:estimate1dlemma}) and the formula for $C_{\alpha}$ follow from
\begin{equation*}
\min(1,1-\alpha) \leq K_{\alpha}(s) \leq \max(1,1-\alpha) \qquad \text{for all $s>0$}.
\end{equation*}
\end{proof}

\begin{proof}[Proof of Theorem {\rm \ref{thm:traceell}}] Without loss of generality, assume that $a_{dj} = 0$ for $1 \leq j \leq d-1$ and $a_{dd}=1$. That is, $\mathcal{E}u = u_{x_d x_d}+\alpha x_d^{-1}u_{x_d}$ and $W^1_{p,\beta}(\mathcal{E}, \mathbb{R}^d_+)=W^1_{p,\alpha,\beta}(\mathbb{R}^d_+)$. The general case follows from the change of variables \eqref{eq:coordchange}. Everywhere below in the proof $N=N(d,\alpha,p,\beta)>0$ unless otherwise stated.

{(i)} We first consider the case $\beta \in (p-1,2p-1)$. Recall that $\gamma \in (-1,p'-1)$ is defined by (\ref{eq:defgamma}).

Suppose that $u \in W^1_{p,\alpha,\beta}(\mathbb{R}^d_+)$. Then, $u$ satisfies
\begin{equation} \label{eq:equationutraceelliptic}
\mathcal{E}u = u_{x_d x_d}+\frac{\alpha}{x_d}u_{x_d} = \sum_{i=1}^{d-1}(F_i)_{x_i}+f
\end{equation}
for some $F_i, f \in L_{p,\beta}(\mathbb{R}^d_+)$ with
\[
\sum_{i=1}^{d-1}\Vert F_i \Vert_{L_{p,\beta}(\mathbb{R}^d_+)} + \Vert f \Vert_{L_{p,\beta}(\mathbb{R}^d_+)} \leq \Vert u \Vert_{W^1_{p,\alpha,\beta}(\mathbb{R}^d_+)}.
\]
Fix any $V \in \mathcal{S}(\mathbb{R}^{d-1})$. Let $v \in W^1_{p',\gamma}(\mathbb{R}^d_+)$ be the solution to
\begin{equation} \label{eq:equationvtraceelliptic}
\left\{
\begin{array}{rll}
   -\Delta v + v = & 0 & \text{in $\mathbb{R}^d_+$}, \\
    v(x',0) = & V(x') & \text{on $\mathbb{R}^{d-1}$}.
\end{array}\right.
\end{equation}
By Proposition \ref{prop:solvabilityellnondeg} with $\alpha=0$,
\begin{equation*}
\Vert v \Vert_{W^1_{p',\gamma}(\mathbb{R}^d_+)} \leq N_{d,p,\beta} \Vert V \Vert_{B^{1-\frac{\gamma+1}{p'}}_{p',p'}(\mathbb{R}^{d-1})}.
\end{equation*}
For $h \geq 0$ and $x_d > 0$, define
\[
z^h(x_d) = \langle \tau_h u, v\rangle = \langle u(\cdot, x_d+h), v(\cdot,x_d)\rangle = \int_{\R^{d-1}} u(x',x_d+h)v(x',x_d) dx'.
\]
Using the Leibnitz rule, we compute the first two derivatives of $z^h$:
\begin{gather*}
(z^h)'(x_d) = \langle \tau_h u_{x_d}, v \rangle + \langle \tau_h u, v_{x_d} \rangle, \\
(z^h)''(x_d) = \langle  \tau_h u_{x_d x_d}, v \rangle + 2 \langle  \tau_h u_{x_d}, v_{x_d} \rangle + \langle  \tau_h u, v_{x_d x_d} \rangle.
\end{gather*}
From equations (\ref{eq:equationutraceelliptic}) and (\ref{eq:equationvtraceelliptic}), it follows that $z^h$ satisfies the following ODE:
\begin{equation} \label{eq:ODEz^helliptic}
    (z^h)''(x_d) +\frac{\alpha}{h+x_d}(z^h)'(x_d)= g^h(x_d),
\end{equation}
where
\begin{equation} \label{eq:g^hformulaelliptic}
\begin{gathered}
g^h(x_d)= \langle \tau_h(u_{x_dx_d}+\alpha x_d^{-1}u_{x_d}),v \rangle + \langle \tau_h u, v_{x_d x_d}\rangle + 2 \langle  \tau_h u_{x_d}, v_{x_d} \rangle + \alpha \langle \tau_h (u/x_d), v_{x_d}\rangle \\
= \langle \tau_h (\sum_{i=1}^{d-1}(F_i)_{x_i}+f),v \rangle + \langle \tau_h u, -\Delta_{x'}v+v\rangle + 2 \langle  \tau_h u_{x_d}, v_{x_d} \rangle + \alpha \langle \tau_h (u/x_d), v_{x_d}\rangle \\
= [\text{integration by parts in $x'$}] = \langle  \tau_h f, v \rangle - \sum_{i=1}^{d-1}\langle  \tau_h F_{i}, v_{x_i} \rangle + \langle \tau_h u, v \rangle + \sum_{i=1}^{d-1}\langle  \tau_h u_{x_i}, v_{x_i} \rangle  \\
+ 2 \langle  \tau_h u_{x_d}, v_{x_d} \rangle + \alpha \langle \tau_h (u/x_d), v_{x_d}\rangle.
\end{gathered}
\end{equation}
Using Lemma \ref{lem:holderweighted}, we estimate
\begin{equation} \label{eq:z^h''L1estell}
\begin{gathered}
\Vert x_d g^h(x_d)\Vert_{L_1(0,\infty)} \leq N\Bigl\{\sum_{i=1}^{d-1}\Vert \tau_h F_i \Vert_{L_{p,\beta}(\R^d_+)}+\Vert \tau_h f\Vert_{L_{p,\beta}(\R^d_+)}+\Vert \tau_h u\Vert_{W^1_{p,\beta}(\R^d_+)}\\
+ \Vert \tau_h (u/x_d)\Vert_{L_{p,\beta}(\R^d_+)}\Bigr\}\Vert v \Vert_{W^1_{p',\gamma}(\R^d_+)} \leq N \Vert u \Vert_{W^1_{p,\alpha,\beta}(\mathbb{R}^d_+)} \Vert V \Vert_{B^{1-\frac{\gamma+1}{p'}}_{p',p'}(\mathbb{R}^{d-1})}.
\end{gathered}
\end{equation}
Here we used $\Vert \tau_h \cdot \Vert_{L_{p,\beta}(\mathbb{R}^{d}_+)} \leq \Vert \cdot \Vert_{L_{p,\beta}(\mathbb{R}^{d}_+)}$ ($\beta \geq 0$) and Hardy's inequality ($\beta > p-1$)
$$\Vert u/x_d \Vert_{L_{p,\beta}(\mathbb{R}^{d}_+)} \leq N_{p,\beta} \Vert u_{x_d} \Vert_{L_{p,\beta}(\mathbb{R}^{d}_+)}.$$

For $h=0$, we write $z(s)$ and $g(s)$ instead of $z^0(s)$ and $g^0(s)$. From Lemma \ref{lem:1dtrace} with $h=0$, it follows that for all $V \in \mathcal{S}(\mathbb{R}^{d-1})$,
\begin{equation*}
\langle U, V \rangle := \lim_{x_d\to 0^+}\langle u(\cdot, x_d), V(\cdot, x_d) \rangle = \frac{1}{1-\alpha}\int_0^{\infty} s (z''(s)+\frac{\alpha}{s}z'(s)) ds.
\end{equation*}
Moreover, by (\ref{eq:z^h''L1estell})
\[
| \langle U, V \rangle | \leq N \Vert u \Vert_{W^1_{p,\alpha,\beta}(\mathbb{R}^d_+)} \Vert V \Vert_{B^{1-\frac{\gamma+1}{p'}}_{p',p'}(\mathbb{R}^{d-1})}.
\]
From Lemma \ref{lem:dualitybesov}, it now follows that
\begin{gather*}
U \in B^{1-\frac{\beta+1}{p}}_{p,p}(\mathbb{R}^{d-1}) = (B^{1-\frac{\gamma+1}{p'}}_{p',p'}(\mathbb{R}^{d-1}))^*, \\
\Vert U \Vert_{B^{1-\frac{\beta+1}{p}}_{p,p}(\mathbb{R}^{d-1})} \leq N \Vert u \Vert_{W^1_{p,\alpha,\beta}(\mathbb{R}^d_+)}.
\end{gather*}

We now show the strong convergence $\lim\limits_{h \to 0^+} \Vert u(\cdot, h)-U \Vert_{B^{1-\frac{\beta+1}{p}}_{p,p}(\mathbb{R}^{d-1})} = 0$. Similarly as above for $h=0$, from the ODE (\ref{eq:ODEz^helliptic}) and Lemma \ref{lem:1dtrace}, it follows that for all $h>0$ and $V \in \mathcal{S}(\R^{d-1})$, we have
\begin{equation*}
    \langle u(\cdot, h), V \rangle = z^h(0) = \frac{1}{1-\alpha}\int_0^{\infty}K(s/h)sg^h(s)ds
\end{equation*}
and the uniform estimate
\begin{equation*} 
    \Vert u(\cdot, h) \Vert_{B^{1-\frac{\beta+1}{p}}_{p,p}(\mathbb{R}^{d-1})} \leq N \Vert u \Vert_{W^1_{p,\alpha,\beta}(\mathbb{R}^d_+)},
\end{equation*}
where $K(s) = 1 - \frac{(1+s)^{\alpha}-1}{s}$. We now estimate
\begin{gather*}
    | \langle u(\cdot,h) - U, V \rangle| = \frac{1}{1-\alpha}\left|\int_0^{\infty}s(K(s/h)g^h(s)-g(s))ds\right| \\
    \leq N_{\alpha} \int_0^\infty sK(s/h)|g^h(s)-g(s)| ds + N_{\alpha}\int_0^\infty |K(s/h)-1|s|g(s)| ds.
\end{gather*}
We denote $\mathcal{W}:=(u, u/x_d, Du, F_i, f)_{i=1}^{d-1}$. To estimate the second term above, by Lemma \ref{lem:holderweighted} and (\ref{eq:g^hformulaelliptic}), we have
\begin{gather*}
\int_0^\infty |K(s/h)-1|s|g(s)| ds \leq N \Vert (K(x_d/h)-1)\mathcal{W}\Vert_{L_{p,\beta}(\mathbb{R}^d_+)}\Vert v \Vert_{W^1_{p',\gamma}(\R^d_+)} \\
\leq N \Vert (K(x_d/h)-1)\mathcal{W}\Vert_{L_{p,\beta}(\mathbb{R}^d_+)} \Vert V \Vert_{B^{1-\frac{\gamma+1}{p'}}_{p',p'}(\mathbb{R}^{d-1})}.
\end{gather*}

To estimate the remaining term, note that from (\ref{eq:g^hformulaelliptic}), we have
\begin{equation*}
\begin{gathered}
g^h(s)-g(s) = \langle \tau_h f - f, v \rangle - \sum_{i=1}^{d-1}\langle \tau_h F_{i} - F_i, v_{x_i} \rangle + \langle \tau_h u - u, v \rangle\\
+ \sum_{i=1}^{d-1}\langle \tau_h u_{x_i} - u_{x_i}, v_{x_i} \rangle + 2\langle  \tau_h u_{x_d}-u_{x_d}, v_{x_d} \rangle + \alpha \langle \tau_h (u/x_d)-(u/x_d), v_{x_d}\rangle.
\end{gathered}
\end{equation*}
Using Lemma \ref{lem:holderweighted}, we estimate
\begin{gather*}
    \int_0^{\infty}sK(s/h)|g^h(s)-g(s)|ds \leq N_{\alpha} \int_0^\infty s|g^h(s)-g(s)| ds \\
    \leq N \Vert (\tau_h-\operatorname{id})\mathcal{W}\Vert_{L_{p,\beta}(\mathbb{R}^d_+)} \Vert V \Vert_{B^{1-\frac{\gamma+1}{p'}}_{p',p'}(\mathbb{R}^{d-1})}.
\end{gather*}
Since $V \in \mathcal{S}(\R^{d-1})$ is arbitrary, from the estimates above we have
\begin{equation*}
\Vert u(\cdot, h)-U \Vert_{B^{1-\frac{\beta+1}{p}}_{p,p}(\mathbb{R}^{d-1})} \leq N \Vert (K(x_d/h)-1)\mathcal{W}\Vert_{L_{p,\beta}(\mathbb{R}^d_+)} +  N \Vert (\tau_h-\operatorname{id})\mathcal{W}\Vert_{L_{p,\beta}(\mathbb{R}^d_+)}.
\end{equation*}
By the dominated convergence theorem, $\lim_{h \to 0^+}\Vert (K(x_d/h)-1)\mathcal{W}\Vert_{L_{p,\beta}(\mathbb{R}^d_+)} = 0$. By the continuity of shifts $\tau_h$ in $L_{p,\beta}(\R^d_+)$ with $\beta \geq 0$ (see Lemma \ref{lem:contshiftslp}), $\lim_{h \to 0^+} \Vert (\tau_h-\operatorname{id})\mathcal{W}\Vert_{L_{p,\beta}(\mathbb{R}^d_+)} = 0$. The strong convergence $\lim\limits_{h \to 0^+}\Vert u(\cdot, h)-U \Vert_{B^{1-\frac{\beta+1}{p}}_{p,p}(\mathbb{R}^{d-1})} = 0$ is proved.

{(ii)} Now suppose $\beta = p-1$. Fix a cut-off function $\chi \in C_0^{\infty}(\R)$ such that $\chi(x_d) = 1$ for $|x_d| \leq 1$. Then $u\chi \in W^1_{p,\alpha,p-1+\varepsilon p}(\R^d_+)$ for all $\varepsilon > 0$. Since $u=u\chi$ for $x_d \leq 1$, from Part (i) it follows that
\begin{equation*}
u(\cdot,x_d) \to u(\cdot,0) \quad \text{strongly in $B^{-\varepsilon}_{p,p}(\mathbb{R}^{d-1})$ as $x_d \to 0^+$.}
\end{equation*}

To prove strong convergence in $B^0_{p,p}(\R^{d-1})$, consider a change of variables
\[
u(x',x_d) = v(x', x_d^{1-\alpha}),
\]
so that
\begin{align*}
    D_d u(\cdot,x_d) & = (1-\alpha)x_d^{-\alpha}D_d v(\cdot,x_d^{1-\alpha}), \\
    (D^2_d+\frac{\alpha}{x_d}D_d)u(\cdot,x_d) & = (1-\alpha)^2 x_d^{-2 \alpha}D^2_d v(\cdot,x_d^{1-\alpha}).
\end{align*}
Then, from $u \in W^1_{p,\alpha,p-1}(\R^d_+)$ and (\ref{eq:equationutraceelliptic}), we get
\begin{equation} \label{eq:Vspacecondition}
\Vert v \Vert_{V_2} := \Vert x_d^{\frac{1}{1-\alpha}} v \Vert_{L_p(\R_+, x_d^{-1} dx_d; W^1_p(\R^{d-1}))} + \Vert x_d^{\frac{1-2\alpha}{1-\alpha}} D^2_d v \Vert_{L_p(\R_+, x_d^{-1} dx_d; H^{-1}_p(\R^{d-1}))} \leq N_{\alpha} \Vert u \Vert_{W^1_{p,\alpha,p-1}(\R^d_+)}.
\end{equation}
From the trace method of interpolation (see \cite[Section 1.8]{Triebel}), (\ref{eq:Vspacecondition}) implies
\begin{gather*}
u(x',0) = v(x',0) \in (W^1_p(\mathbb{R}^{d-1}),H^{-1}_{p}(\mathbb{R}^{d-1}))_{\frac{1}{2},p} = B^{0}_{p,p}(\mathbb{R}^{d-1}), \\
\Vert u(x',0) \Vert_{B^{0}_{p,p}(\mathbb{R}^{d-1})} \leq N \Vert v \Vert_{V_2} \leq N \Vert u \Vert_{W^1_{p,\alpha,p-1}(\R^d_+)}.
\end{gather*}

Now suppose $\alpha \in [-(p-1),\frac{p-1}{2p-1}]$, that is $\frac{p}{1-\alpha}-1 \geq 0$ and $\frac{p(1-2\alpha)}{1-\alpha}-1 \geq 0$. In this case,
\begin{equation*}
\lim_{h \to 0^+} \Vert \tau_h v - v \Vert_{V_2} = 0,
\end{equation*}
and therefore,
\begin{equation*}
\lim_{h \to 0^+} \Vert v(\cdot,h) - v(\cdot,0) \Vert_{B^{0}_{p,p}(\mathbb{R}^{d-1})} = \lim_{h \to 0^+} \Vert u(\cdot,h) - u(\cdot,0) \Vert_{B^{0}_{p,p}(\mathbb{R}^{d-1})} = 0.
\end{equation*}

For $\alpha < -(p-1)$, we can reduce to $\alpha=0$. If $\alpha<0$, from Lemma \ref{lem:structureW1} and Hardy's inequality for $\beta > (1+\alpha)p-1$, it follows that
\[
W^1_{p,\alpha,p-1}(\mathbb{R}^{d}_+) = (1-\Delta_{x'})^{\frac{1}{2}}W^2_{p,\alpha,p-1}(\mathbb{R}^{d}_+) \hookrightarrow (1-\Delta_{x'})^{\frac{1}{2}}W^2_{p,p-1}(\mathbb{R}^d_+) = W^1_{p,0,p-1}(\mathbb{R}^{d}_+),
\]
and therefore,
\begin{equation*}
u(\cdot,x_d) \to u(\cdot,0), \quad \text{strongly in $B^{0}_{p,p}(\mathbb{R}^{d-1})$ as $x_d \to 0^+$},
\end{equation*}
holds for all $\alpha \in (-\infty, \frac{p-1}{2p-1}]$. The proof is complete.
\end{proof}

We omit the proof of Theorem \ref{thm:existsolell}, as it is nearly identical to the proof of the parabolic result Theorem \ref{thm:existsolpar}, which is given in the next Section.

\section{Proofs of main results: parabolic equations} \label{sec:proofspar}

\subsection{Existence of solutions with \texorpdfstring{$F_d \neq 0$}{mathformtext}}
In the proof of Theorem \ref{thm:existsolpar}, we will use a result about the solvability of equations (\ref{equationnonhomweighted})-(\ref{equationparabolicdivnonhom}) with the conormal boundary condition
\begin{equation} \label{eq:conormalbc}
    \lim_{x_d\to 0^+}x_d^{\alpha}\Bigl(\sum_{j=1}^d a_{dj}(x_d)D_j u - F_d\Bigr) = 0.
\end{equation}
\begin{defn}
    Let $p \in (1,\infty), \alpha \in \mathbb{R}, \beta < (1+\alpha)p-1$, and $x_d^{\alpha}F, x_d^{\alpha}f \in L_{1,\text{\rm loc}}(\R^{d+1}_+)$ {\rm (}resp. $L_{1,\text{\rm loc}}(\R^{d}_+)${\rm )}. We say that $u \in W^{0,1}_{p,\beta}(\R^{d+1}_+)$ {\rm (}resp. $u \in W^{1}_{p,\beta}(\R^d_+)${\rm )} is a weak solution to \eqref{equationparabolicdivnonhom}-\eqref{eq:conormalbc} {\rm (}resp. \eqref{equationnonhomweighted}-\eqref{eq:conormalbc}{\rm )} if \eqref{eq:weaksolell} {\rm (}resp. \eqref{eq:weaksolutionellipticdefinition}{\rm )} holds for all $\varphi \in C_0^\infty(\overline{\mathbb{R}^{d+1}_+})$ {\rm (}resp. $\varphi \in C_0^{\infty}(\overline{\R^d_+}){\rm )}$.
\end{defn}
\begin{remark}
The condition $\beta < (1+\alpha)p-1$ ensures that if $u, Du \in L_{p,\beta}(\R^d_+)$, then by H\"older's inequality $x_d^{\alpha} u, x_d^{\alpha} Du \in L_{1,\text{loc}}(\overline{\R^d_+})$. Therefore, the integrals in (\ref{eq:weaksolell}) and (\ref{eq:weaksolutionellipticdefinition}) are well defined.
\end{remark}

The following proposition is a special case of \cite[Theorem 2.4]{DongPhanIndiana2023}.

\begin{proposition} \label{prop:solvabilityconormal}
    Let $\alpha \in (-1,\infty)$, $p \in (1,\infty)$, and $\beta \in (-1, (1+\alpha)p-1)$. Then, for any $\lambda > 0$ and $F,f \in L_{p,\beta}(\R^{d+1}_+)$ {\rm (}resp. $F,f \in L_{p,\beta}(\R^d_+)${\rm )}, there exists a unique weak solution $u \in W^{0,1}_{p,\beta}(\R^{d+1}_+)$ {\rm (}resp. $u \in W^1_{p,\beta}(\R^d_+)${\rm )} to \eqref{equationparabolicdivnonhom}-\eqref{eq:conormalbc} {\rm (}resp. \eqref{equationnonhomweighted}-\eqref{eq:conormalbc}{\rm )}. Moreover, $u$ satisfies
    \[
    \Vert Du \Vert_{L_{p,\beta}} + \sqrt{\lambda} \Vert u \Vert_{L_{p,\beta}} \leq N\Vert F \Vert_{L_{p,\beta}} + N \frac{1}{\sqrt{\lambda}} \Vert f \Vert_{L_{p,\beta}}, 
    \]
    where $N=N(d,\alpha,p,\beta,\kappa) > 0$ and $L_{p,\beta}=L_{p,\beta}(\R^{d+1}_+)$ {\rm (}resp. $L_{p,\beta} = L_{p,\beta}(\R^d_+)${\rm )}.
\end{proposition}
\begin{remark} \label{rmk:refremark}
Even though \cite[Theorem 2.4]{DongPhanIndiana2023} is stated for $a_0=c=1$, the same proof with minimal changes works for $a_0=a_0(x_d)$ and $c=c(x_d)$ as well. This also applies to \cite[Theorem 2.5]{DongPhan2021TAMS} and \cite[Theorem 2.8]{DongPhan2021TAMS}, which are used in the proofs of Propositions \ref{prop:solvabilityparnondeg} and \ref{prop:solvabilityellnondeg}, respectively.
\end{remark}

\begin{proof}[Proof of Theorem {\rm \ref{thm:existsolpar}}]
Without loss of generality, we assume $a_{dd} = 1$. The general case can be reduced to $a_{dd} = 1$ by the change of variables $y_d=y_d(x_d)$ from (\ref{eq:coordchange}).
    
    Consider a change of variables $u = x_d^{1-\alpha} w$. Equivalently, $w = x_d^{\alpha-1} u$. Plugging $u = x_d^{1-\alpha} w$ in (\ref{equationparabolicdivnonhom}), and multiplying the equation by $x_d^{1-\alpha}$ give an equivalent equation for $w$ (using $a_{dd}=1$):
\begin{equation} \label{eq:equationw1}
\begin{gathered}
        x_d^{2-\alpha}a_0 w_t - \sum_{i,j=1}^d D_i (x_d^{2-\alpha} a_{ij} D_j w) - (1-\alpha)\sum_{j=1}^{d-1}x_d^{1-\alpha}(a_{jd}-a_{dj})D_j w + \lambda x_d^{2-\alpha} c(x_d) w \\
        = -\sum_{i=1}^d D_i(x_d F_i) + (1-\alpha)F_d + x_d f.
\end{gathered}
\end{equation}
    The first-order derivatives of $w$ can be absorbed into the main symbol:
    \begin{equation*}
        \sum_{i,j=1}^d D_i(x_d^{2-\alpha} a_{ij}D_j w) + (1-\alpha)\sum_{j=1}^{d-1}x_d^{1-\alpha}(a_{jd}-a_{dj})D_j w = \sum_{i,j=1}^d D_i(x_d^{2-\alpha} b_{ij}D_j w),
    \end{equation*}
    where
    \begin{equation*} 
    \begin{array}{rlr}
    b_{ij} = & a_{ij}, & 1 \leq i, j \leq d-1, \\
    b_{dd} = & 1, & \\
    b_{dj}-a_{dj} = & (1-\alpha)x_d^{\alpha-2}\int_0^{x_d}(a_{jd}(s)-a_{dj}(s))s^{1-\alpha} ds, & 1 \leq j \leq d-1, \\
    b_{jd}+b_{dj} = & a_{jd} + a_{dj}, & 1 \leq j \leq d-1.
    \end{array}
    \end{equation*}
    The term $(1-\alpha)F_d$ can be absorbed into the divergence term by introducing
    \[
    G_d(t,x',x_d) = (1-\alpha)x_d^{\alpha-2} \int_{x_d}^{\infty} F_d(t,x',s) ds,
    \]
    so that
    \[
    D_d(x_d^{2-\alpha} G_d) = -(1-\alpha)F_d.
    \]
    From Hardy's inequality for $\beta > p-1$ (Lemma \ref{lem:Hardy} (i) with $\alpha=0$), it follows that
    \[
    \Vert G_d \Vert_{L_{p,\beta+(1-\alpha)p}(\R^{d+1}_+)} = \Vert x_d^{2-\alpha}G_d \Vert_{L_{p,\beta-p}(\R^{d+1}_+)} \leq (1-\alpha)N_{p,\beta} \Vert F_d \Vert_{L_{p,\beta}(\R^{d+1}_+)}.
    \]
    Therefore, we can rewrite the equation for $w$ as
    \begin{equation} \label{eq:weightedparnewcoord}
    \begin{gathered}
        x_d^{2-\alpha}a_0(x_d)w_t - \sum_{i,j=1}^d D_i (x_d^{2-\alpha} b_{ij}(x_d)D_j w) + \lambda x_d^{2-\alpha} c(x_d) w
        = -\sum_{j=1}^d D_j(x_d^{2-\alpha}\Tilde{F}_j) + x_d^{2-\alpha} \Tilde{f},
    \end{gathered}
    \end{equation}
    where
    \begin{equation*}
    \begin{gathered}
        \Tilde{F}_j = x_d^{\alpha-1}F_j, \qquad j = 1, \dots, d-1, \\
        \Tilde{F}_d = x_d^{\alpha-1}F_d+G_d, \qquad \Tilde{f} = x_d^{\alpha-1}f.
    \end{gathered}
    \end{equation*}
    Clearly, it also holds that
        \begin{equation*}
    \begin{aligned}
        \Vert \Tilde{F}_j \Vert_{L_{p,\beta+(1-\alpha)p}(\R^{d+1}_+)} & = \Vert F_j \Vert_{L_{p,\beta}(\R^{d+1}_+)}, \qquad j = 1, \dots, d-1, \\
        \Vert \Tilde{f} \Vert_{L_{p,\beta+(1-\alpha)p}(\R^{d+1}_+)} & = \Vert f \Vert_{L_{p,\beta}(\R^{d+1}_+)}, \\
        \Vert \Tilde{F}_d \Vert_{L_{p,\beta+(1-\alpha)p}(\R^{d+1}_+)} & \leq N_{\alpha,p,\beta} \Vert F_d \Vert_{L_{p,\beta}(\R^{d+1}_+)}.
    \end{aligned}
    \end{equation*}
    From Proposition \ref{prop:solvabilityconormal}, there exists a unique weak solution $w \in W^{0,1}_{p,\beta+(1-\alpha)p}(\R^{d+1}_+)$ to \eqref{eq:weightedparnewcoord} satisfying the conormal boundary condition. Moreover, $w$ satisfies
    \begin{equation} \label{eq:estimatewbyFf}
    \begin{gathered}
    \Vert Dw \Vert_{L_{p,\beta+(1-\alpha)p}(\R^{d+1}_+)} + \sqrt{\lambda} \Vert w \Vert_{L_{p,\beta+(1-\alpha)p}(\R^{d+1}_+)} \leq N\Vert \Tilde{F} \Vert_{L_{p,\beta+(1-\alpha)p}(\R^{d+1}_+)} \\
    + N \frac{1}{\sqrt{\lambda}} \Vert \Tilde{f} \Vert_{L_{p,\beta+(1-\alpha)p}(\R^{d+1}_+)} \leq N \Vert F \Vert_{L_{p,\beta}(\R^{d+1}_+)} + N \frac{1}{\sqrt{\lambda}} \Vert f \Vert_{L_{p,\beta}(\R^{d+1}_+)}.
    \end{gathered}
    \end{equation}
    It is left to estimate $u = x_d^{1-\alpha}w$, which is a weak solution to (\ref{equationparabolicdivnonhom}) by the definition of $w$. Clearly,
    \begin{gather*}
    \Vert u \Vert_{L_{p,\beta}(\R^{d+1}_+)} = \Vert w \Vert_{L_{p,\beta+(1-\alpha)p}(\R^{d+1}_+)}, \quad
    \Vert D_{x'}u \Vert_{L_{p,\beta}(\R^{d+1}_+)} = \Vert D_{x'}w \Vert_{L_{p,\beta+(1-\alpha)p}(\R^{d+1}_+)}.
    \end{gather*}
    Note that $D_d u = x_d^{1-\alpha}(D_d w + (1-\alpha)x_d^{-1}w).$ From Hardy's inequality ($\beta > \alpha p - 1$), we have
    \[
    \Vert w/x_d \Vert_{L_{p,\beta+(1-\alpha)p}(\R^{d+1}_+)} \leq N_{\alpha,p,\beta} \Vert D_d w \Vert_{L_{p,\beta+(1-\alpha)p}(\R^{d+1}_+)},
    \]
    and therefore,
    \[
    \Vert D_{d}u \Vert_{L_{p,\beta}(\R^{d+1}_+)} \leq N_{\alpha,p,\beta} \Vert D_d w \Vert_{L_{p,\beta+(1-\alpha)p}(\R^{d+1}_+)}.
    \]
    From the above estimates and (\ref{eq:estimatewbyFf}), the desired estimate for $u$ follows. The theorem is proved.
\end{proof}

\begin{remark}
According to Proposition \ref{prop:solvabilityconormal}, when
\begin{equation*}
    \alpha > -1, \quad \beta \in (-1, (1+\alpha)p-1),
\end{equation*}
equations (\ref{equationparabolicdivnonhom}) and (\ref{equationnonhomweighted}) can be solved with the conormal boundary condition directly. Therefore, when $\alpha \in (0,1)$, the assertions of Theorems \ref{thm:existsolpar} and \ref{thm:existsolell} hold true for $\beta \in (-1,2p-1)$.
\end{remark}

\subsection{Existence and uniqueness of solutions} We first treat the non-degenerate range $\beta \in (\alpha^+ p-1, p-1)$.
\begin{proposition} \label{prop:solvabilityparnondeg}
Let $\alpha \in (-1,1), p \in (1,\infty), \beta \in (\alpha^+ p-1, p-1)$, and $\kappa \in (0,1]$. Suppose that \eqref{eq:ellipticitycond} is satisfied. Then, for any $\lambda > 0, F\in  L_{p,\beta}(\mathbb{R}^{d+1}_+; \mathbb{R}^d), f \in  L_{p,\beta}(\mathbb{R}^{d+1}_+)$, and $U \in B^{\frac{1}{2}-\frac{\beta+1}{2p}, 1-\frac{\beta+1}{p}}_{p,p}(\mathbb{R}^{d})$, there exists a unique weak solution $u \in W^{0,1}_{p,\beta}(\mathbb{R}^{d+1}_+)$ to \eqref{equationparabolicdivnonhom}--\eqref{dirichletparabolic}. Moreover, $u$ satisfies
\begin{equation} \label{eq:estdivA_p}
    \Vert Du \Vert + \sqrt{\lambda} \Vert u \Vert 
    \leq N \Vert F \Vert + \frac{N}{\sqrt{\lambda}} \Vert f \Vert + N(\lambda^{\frac{1}{2}-\frac{\beta+1}{2p}}\Vert U \Vert_{L_{p}(\mathbb{R}^d)}+[U]),
    \end{equation}
where $\Vert \cdot \Vert = \Vert \cdot \Vert_{L_{p,\beta}(\mathbb{R}^{d+1}_+)}, [U]=[U]_{B^{\frac{1}{2}-\frac{\beta+1}{2p}, 1-\frac{\beta+1}{p}}_{p,p}(\mathbb{R}^{d})}$, and $N=N(d,\alpha,p,\beta, \kappa) > 0$.
\end{proposition}
\begin{proof}
For $U=0$ and $\beta \in (\alpha p-1, p-1)$, Proposition \ref{prop:solvabilityparnondeg} follows from \cite[Theorem 2.5]{DongPhan2021TAMS} (recall Remark \ref{rmk:refremark}). Now consider $U \neq 0$. We further assume $\lambda = 1$ as the general case follows by scaling. We may also assume $a_0(x_d) = 1$. Indeed, consider a change of variables given by (using $\alpha > -1$)
\begin{equation} \label{eq:changevara0}
z_d = \left((1+\alpha)\int_0^{x_d} s^{\alpha} a_0(s) ds\right)^{\frac{1}{1+\alpha}},
\end{equation}
with $x_i$ unchanged for $i=1,\dots,d-1$. It holds that
\begin{equation*}
(z_d/x_d)^{\alpha} a_0(x_d)^{-1}D_{x_d} = D_{z_d}.
\end{equation*}
Multiplying equation (\ref{equationparabolicdivnonhom}) by $(z_d/x_d)^{\alpha} a_0(x_d)^{-1}$ and using the above formula give an equation with $\Tilde{a}_0(z_d) = 1$. We therefore assume further that $a_0(x_d) = 1$.

Consider the following extensions of $U$ to $\R^{d+1}_+$: $v(t,x) = [Q^{\alpha}_{x_d} \ast U](t,x')$ and $u_0(t,x) = \chi(x_d)v(t,x)$, where $\chi$ is from (\ref{eq:chicutoff}). From Theorem \ref{thm:extbesovweightedpar} with $q=p$ and $r = 1-\frac{\beta+1}{p}$, we have
\[
\Vert Dv \Vert_{L_{p,\beta}(\R^{d+1}_+)} \leq N [U]_{B^{\frac{1}{2}-\frac{\beta+1}{2p}, 1-\frac{\beta+1}{p}}_{p,p}(\mathbb{R}^{d})}.
\]
Since $\int_0^\infty \int_{\R^{d-1}} Q^{\alpha}_{x_d}(\tau,z')dz'd\tau=1$ for all $x_d>0$, Young's inequality gives 
\[
\Vert v(t,x',x_d) \Vert_{L_p(\R^d)} \leq \Vert U \Vert_{L_p(\R^d)}
\]
for all $x_d > 0$. The above estimates imply
\begin{equation*}
\Vert u_0 \Vert_{W^{0,1}_{p,\beta}(\R^{d+1}_+)} \leq N \Vert U \Vert_{B^{\frac{1}{2}-\frac{\beta+1}{2p}, 1-\frac{\beta+1}{p}}_{p,p}(\mathbb{R}^{d})}.
\end{equation*}
Moreover, $u_0$ satisfies
\begin{equation*}
x_d^{\alpha} \partial_t u_0 = \sum_{i=1}^d D_i(x_d^{\alpha} D_i u_0) - D_d(x_d^{\alpha} \chi'(x_d) v) - x_d^{\alpha}\chi'(x_d) D_d v,
\end{equation*}
where all $Du_0, \chi' v, \chi' D_d v \in L_{p,\beta}(\R^{d+1}_+)$. Thus it suffices to solve for $u-u_0$ with the zero boundary condition. The proof is complete.
\end{proof}

\begin{remark}
As can be seen from the proof, the assertion of Proposition \ref{prop:solvabilityparnondeg} holds for $\alpha \leq -1$ as well if we assume $a_0 = 1$. The condition $\alpha > -1$ was used to apply the change of variables (\ref{eq:changevara0}) when $a_0 \neq \operatorname{const}$.
\end{remark}

The proof of Theorem \ref{thm:solvpar} is similar to Theorem \ref{thm:solvell}, but modifications need to be made in the existence of solutions.
\begin{proof}[Proof of Theorem {\rm \ref{thm:solvpar}}] Without loss of generality, assume that $a_{dj}=0$ for $1 \leq j \leq d-1$ and $a_{dd}=1$. The general case follows from the change of variables \eqref{eq:coordchange}. Everywhere below in the proof $N=N(d,\alpha,p,\beta,\kappa)>0$ unless otherwise stated.

For compactness of notation, let $\mathcal{P}_{\text{non-div}}$ be as in (\ref{eq:operatornondivpar}) and
\begin{equation*}
\mathcal{P}_{\text{div}} = a_0(x_d)\partial_t- D^2_d - \frac{\alpha}{x_d}D_d - \sum_{i=1}^{d-1}\sum_{j=1}^d D_i(a_{ij}(x_d)D_j\;) + \lambda c(x_d).
\end{equation*}

Similar to Theorem \ref{thm:solvell}, the uniqueness of solutions in the parabolic case follows from mollification in the $t,x'$ variables and the uniqueness part of Theorem \ref{thm:solvnondivpar}.

We now consider $\lambda = 1$ and $\beta \in I_p(\alpha) \cap [p-1,2p-1)$, where $I_p(\alpha)$ is from (\ref{eq:rangeIp}). As $F_d=0$ and $a_{dj}=\delta_{dj}$, equations (\ref{equationparabolicdivnonhom})-(\ref{dirichletparabolic}) are equivalent to 
\begin{equation} \label{eq:parabolicnonhom}
\begin{cases}
\mathcal{P}_{\text{div}}u = -\sum_{i=1}^{d-1}D_i F_i+f & \text{in $\R^{d+1}_+$,} \\
u(t,x',0)=U(x') & \text{on $\R^{d}$}.
\end{cases}
\end{equation}

Let $W = (\partial_t-\Delta_{x'}+1)^{-\frac{1}{2}} U \in B^{1-\frac{\beta+1}{2p},2-\frac{\beta+1}{p}}_{p,p}(\mathbb{R}^{d})$. By Lemma \ref{lem:isomorphisms},
\begin{equation*}
\Vert W \Vert_{B^{1-\frac{\beta+1}{2p},2-\frac{\beta+1}{p}}_{p,p}(\mathbb{R}^{d})} \leq N_{d,p,\beta} \Vert U \Vert_{B^{\frac{1}{2}-\frac{\beta+1}{2p},1-\frac{\beta+1}{p}}_{p,p}(\mathbb{R}^{d})}.
\end{equation*}
Using Theorem \ref{thm:solvnondivpar}, we find unique strong solutions $u_0,u_1,\dots,u_{d-1}, w \in W^{1,2}_{p,\alpha,\beta}(\R^{d+1}_+)$ to 
\begin{equation*}
	\begin{cases}
	\mathcal{P}_{\text{non-div}}w = 0, \\
	w(t,x',0) = W,
	\end{cases} \quad
	\begin{cases}
        \mathcal{P}_{\text{non-div}}u_0 = f, \\
        u_0(t,x',0)=0,
    \end{cases} \quad
    \begin{cases}
        \mathcal{P}_{\text{non-div}}u_i = F_i, \\
        u_i(t,x',0)=0,
    \end{cases}
\end{equation*}
where $i=1,\dots,d-1$. We now show that
\[
u = -\sum_{i=1}^{d-1}D_i u_i+u_0+(\partial_t-\Delta_{x'}+1)^{\frac{1}{2}}w \in W^{0,1}_{p,\beta}(\R^{d+1}_+)
\]
is a weak solution to (\ref{eq:parabolicnonhom}). 
For compactness of notation, denote
\begin{gather*}
\Vert \cdot \Vert = \Vert \cdot \Vert_{L_{p,\beta}(\R^{d+1}_+)}, \quad [\; \cdot \;] = [\; \cdot \;]_{B^{1-\frac{\beta+1}{2p},2-\frac{\beta+1}{p}}_{p,p}(\R^d)}, \\
    \vvvert{\cdot}_{\lambda} = \Vert \partial_t (\cdot) \Vert + \Vert D_{x'} D (\cdot)\Vert + \Vert (D^2_d  + \alpha x_d^{-1} D_d) (\cdot) \Vert + \lambda^{1/2} \Vert D(\cdot) \Vert + \lambda \Vert \cdot \Vert.
\end{gather*}
From Theorem \ref{thm:solvnondivpar}, we have
\begin{equation} \label{eq:estimatesuiu0}
    \vvvert{u_0}_{1} \leq N\Vert f \Vert, \quad \vvvert{u_i}_{1} \leq N\Vert F_i \Vert, \qquad i=1,\dots,d-1.
\end{equation}
From estimates (\ref{eq:estimatesuiu0}), we have
\begin{equation*}
\begin{array}{rlr}
    \Vert D D_i u_i \Vert + \Vert D_i u_i \Vert \leq & N\Vert F_i \Vert, & i=1,\dots,d-1, \\
    \Vert D u_0 \Vert + \Vert u_0 \Vert \leq & N\Vert f \Vert. &
\end{array}
\end{equation*}
To estimate $(\partial_t-\Delta_{x'}+1)^{\frac{1}{2}}w$, note that
\begin{gather*}
    \Vert D_{x'}(\partial_t-\Delta_{x'}+1)^{\frac{1}{2}}w \Vert + \Vert (\partial_t-\Delta_{x'}+1)^{\frac{1}{2}}w \Vert \leq N \Vert w \Vert_{L_{p,\beta}(\R_+;W^{1,2}_p(\R^d))} \\
     \leq N \vvvert{w}_1 \leq N \Bigl(\Vert W \Vert_{L_p(\R^{d})} + [W]\Bigr).
\end{gather*}
The estimate of $D_{d}(\partial_t-\Delta_{x'}+1)^{\frac{1}{2}}w$ follows from Lemma \ref{lem:embedsoltnop} (using $\beta \in I_p(\alpha)$):
\begin{gather*}
    \Vert D_{d}(\partial_t-\Delta_{x'}+1)^{\frac{1}{2}}w \Vert \leq N \Vert w \Vert_{H^{1/2,1}_p(\mathbb{R}^d; W^1_{p,\beta}(\mathbb{R}_+))} \leq N (\Vert W \Vert_{L_p(\R^{d})} + [W]).
\end{gather*}
Combining the estimates above, we get the desired estimate
\begin{equation} \label{eq:estimateparaboliclambda1}
    \Vert D u \Vert + \Vert u \Vert \leq N\sum_{i=1}^{d-1} \Vert F_i \Vert + \Vert f \Vert+ N\bigl(\Vert (\partial_t -\Delta_{x'}+1)^{-\frac{1}{2}} U \Vert_{L_p(\R^{d})} + [(\partial_t -\Delta_{x'}+1)^{-\frac{1}{2}} U]\bigr).
\end{equation}

For $\lambda = 1$ and the uncovered range $\alpha \in (0,1), \beta \in [p-1,(1+\alpha)p-1]$, the existence of solutions and estimate (\ref{eq:estimateparaboliclambda1}) follow by interpolation. Indeed, the case $\beta \in ((1+\alpha)p-1,2p-1)$ is covered above, and the case $\beta \in (\alpha p-1, p-1)$ is covered in Proposition \ref{prop:solvabilityparnondeg}. 

For general $\lambda > 0$, the existence of solutions and estimate (\ref{eq:estsolpar}) follow by scaling.
\end{proof}

\subsection{Existence of trace}

The proof of Theorem \ref{thm:tracepar} uses a duality argument, similar to the one used in the proof of Theorem \ref{thm:traceell}. However, a family of test functions needs to be used due to the presence of the non-constant coefficient $a_0(x_d)$ in the $u_t$ term. We use the duality argument in the case $\beta=p-1$ as well instead of trace interpolation, due to the lack of regularity in the time variable.

Let $a_0(x_d):\mathbb{R}_+ \to \mathbb{R}$ be a measurable function such that $\kappa \leq a_0(x_d) \leq \kappa^{-1}$ for some $\kappa \in (0,1]$. For each $V \in \mathcal{S}(\mathbb{R}^{d})$ and $h \geq 0$, let $v^h$ be the solution to
\begin{equation} \label{eq:definitionvh}
\left\{
\begin{array}{rll}
   a_0(x_d+h) v^h_t-\Delta v^h + v^h  = & 0 & \text{in $\mathbb{R}^{d+1}_+$}, \\
    v^h(t, x',0) = & V(t,x') & \text{on $\mathbb{R}^{d}$},
\end{array}\right.
\end{equation}
provided by Theorem \ref{thm:solvnondivpar}. For $h = 0$, we write $v$ instead of $v^0$.

From Theorems \ref{thm:solvpar} and \ref{thm:solvnondivpar}, Lemma \ref{lem:solutionoperatorwelldefined}, and Proposition \ref{prop:solvabilityparnondeg} with $\alpha=0$, it follows that for all $q \in (1,\infty)$ and $\gamma \in (-1,q-1]$,
\begin{equation} \label{eq:basicestimatevh}
\begin{gathered}
\Vert v^h \Vert_{W^{1,2}_{q,\gamma}(\mathbb{R}^{d+1}_+)} \leq N \Vert V \Vert_{B^{1-\frac{\gamma+1}{2q},2-\frac{\gamma+1}{q}}_{q,q}(\mathbb{R}^{d})}, \\
\Vert D v^h \Vert_{L_{q,\gamma}(\mathbb{R}^{d+1}_+)} + \Vert v^h \Vert_{L_{q,\gamma}(\mathbb{R}^{d+1}_+)} \leq N \Vert V \Vert_{B^{\frac{1}{2}-\frac{\gamma+1}{2q},1-\frac{\gamma+1}{q}}_{q,q}(\mathbb{R}^{d})},
\end{gathered}
\end{equation}
where $N=N(d,q,\gamma,\kappa) > 0$. Moreover, $v^h, D_d v^h \in C((0,\infty),L_q(\R^d))$.

In the remaining lemmas in this subsection, we show some continuity estimates of $v^h$ in $h$.
\begin{lem} \label{lem:continuityoftranslationsL_p}
    Let $q \in [1,\infty)$ and $\gamma \in \mathbb{R}$. For any $a \in L_{\infty}(\mathbb{R}_+)$ and $w \in L_{q,\gamma}(\mathbb{R}^{d+1}_+)$,
    \begin{equation*}
        \Vert (a(x_d+h)-a(x_d))w \Vert_{L_{q,\gamma}(\mathbb{R}^{d+1}_+)} \to 0 \qquad \text{as $h \to 0^+$}.
    \end{equation*}
\end{lem}
The proof follows by a standard argument using the density of $C_0(\R^{d+1}_+)$ in $L_{q,\gamma}(\mathbb{R}^{d+1}_+)$.

\begin{lem} \label{lem:continuityv^h}
    Let $V \in \mathcal{S}(\mathbb{R}^d), h > 0$, and $v^h(t,x)$ as defined in \eqref{eq:definitionvh}. For any $q \in (1,\infty)$ and $\gamma \in (-1,q-1)$, we have
    \begin{equation*}
        \Vert v^h-v \Vert_{W^{1,2}_{q,\gamma}(\mathbb{R}^{d+1}_+)} \to 0 \qquad \text{as $h \to 0^+$}.
    \end{equation*}
    If $a_0(x_d)$ is differentiable with $a_0' \in BUC(\R_+)$, then
    \begin{equation*}
        \Vert v^h-v \Vert_{W^{0,1}_{q,\gamma}(\mathbb{R}^{d+1}_+)} \leq N \Tilde{\omega}(h)\Vert V \Vert_{B^{\frac{1}{2}-\frac{\gamma+1}{2q},1-\frac{\gamma+1}{q}}_{q,q}(\mathbb{R}^{d})},
    \end{equation*}
    where $\Tilde{\omega}(h) = h+\omega_{a_0'}(h)$ and $N=N(d,q,\gamma,\kappa,\Vert a_0' \Vert_{L_{\infty}(\R_+)}) > 0$.
\end{lem}
Recall that $\omega_{a_0'}(h)$ is the modulus of continuity of $a_0'$.
\begin{proof}
    Let $w^h = v^h - v$. Then $w^h$ is a strong solution to
    \begin{equation*}
    \left\{
    \begin{array}{rll}
        a_0(x_d+h) w^h_t-\Delta w^h + w^h  = & -(a_0(x_d+h)-a_0(x_d))v_t & \text{in $\mathbb{R}^{d+1}_+$}, \\
        w^h(t, x',0) = & 0 & \text{on $\mathbb{R}^{d}$}.
    \end{array}\right.
    \end{equation*}
    By Theorem \ref{thm:solvnondivpar} with $\alpha=0$ and Lemma \ref{lem:continuityoftranslationsL_p},
    \begin{equation*}
        \Vert v^h-v \Vert_{W^{1,2}_{q,\gamma}(\mathbb{R}^{d+1}_+)} \leq N \Vert (\tau_h a_0 - a_0)v_t \Vert_{L_{q,\gamma}(\mathbb{R}^{d+1}_+)} \to 0 \qquad \text{as $h \to 0^+$},
    \end{equation*}
    where $N=N(d,q,\gamma,\kappa)>0$.

    Now suppose that $a_0' \in BUC(\R_+)$. Denote $\Lambda = \Vert a_0' \Vert_{L_{\infty}(\R_+)}$, and $\omega(h) = \omega_{a_0'}(h)$. Then
    \begin{gather*}
        -(a_0(x_d+h)-a_0(x_d))v_t = k_h(x_d)(\Delta v - v),
    \end{gather*}
    where
    \[
    k_h(x_d) = -\frac{a_0(x_d+h)-a_0(x_d)}{a_0(x_d)}.
    \]
    Therefore,
    \begin{gather*}
        a_0(x_d+h) w^h_t-\Delta w^h + w^h = k_h(x_d)(\Delta v - v) 
        = \sum_{j=1}^{d}D_i(k_h D_i v)-(k_h v + k_h' D_d v).
    \end{gather*}
    Note that
    \begin{gather*}
        k_h'(x_d) = -\frac{(a_0'(x_d+h)-a_0'(x_d))a_0(x_d)-(a_0(x_d+h)-a_0(x_d))a_0'(x_d)}{a_0^2(x_d)}.
    \end{gather*}
    It is straightforward to verify that
    \begin{gather*}
        \Vert k_h \Vert_{L_{\infty}(\R_+)} \leq \kappa^{-1}\Lambda h, \\
        \Vert k_h' \Vert_{L_{\infty}(\R_+)} \leq \omega(h)\kappa^{-1}+\Lambda^2\kappa^{-2}h \leq N_{\kappa,\Lambda}(h+\omega(h)).
    \end{gather*}
    From the estimate (\ref{eq:estdivA_p}), it now follows that
    \begin{gather*}
        \Vert v^h-v \Vert_{W^{0,1}_{q,\gamma}(\mathbb{R}^{d+1}_+)} \leq N(\Vert k_h D v\Vert + \Vert k_h v\Vert + \Vert k_h' D_{d} v\Vert) \\
        \leq N (\Vert k_h \Vert_{L_{\infty}(\R_+)}+\Vert k_h' \Vert_{L_{\infty}(\R_+)})\Vert v \Vert_{W^{0,1}_{q,\gamma}(\mathbb{R}^{d+1}_+)} \leq N \Tilde{\omega}(h)\Vert V \Vert_{B^{\frac{1}{2}-\frac{\gamma+1}{2q},1-\frac{\gamma+1}{q}}_{q,q}(\mathbb{R}^{d})},
    \end{gather*}
    where $\Vert \cdot \Vert = \Vert \cdot \Vert_{L_{q,\gamma}(\R^{d+1}_+)}$ and $N = N(d,q,\gamma,\kappa,\Lambda)>0$.
\end{proof}

\begin{lem} \label{lem:shiftintegralcontinuity}
    Let $p \in (1,\infty)$, $\beta \in (p-1,2p-1), V \in \mathcal{S}(\R^d), h > 0$ and $v^h(t,x)$ as defined in \eqref{eq:definitionvh}. Suppose $u \in L_{p,\beta}(\mathbb{R}^{d+1}_+)$, and let
    \begin{equation*}
    r(h) = \int_0^{\infty}s|\langle \tau_h u, v^h\rangle - \langle u, v\rangle| ds + \sum_{i=1}^d  \int_0^{\infty}s|\langle \tau_h u, v^h_{x_i}\rangle - \langle u, v_{x_i}\rangle| ds.
    \end{equation*}
    Then
    \begin{equation*}
    	\lim_{h \to 0^+} r(h) = 0.
	\end{equation*}
	If $a_0(x_d)$ is differentiable with $a_0' \in BUC(\R_+)$, then
	\begin{equation} \label{eq:estimaterh}
	r(h) \leq N\left(\Vert \tau_h u - u\Vert_{L_{p,\beta}(\mathbb{R}^{d+1}_+)}+\Tilde{\omega}(h)\Vert u\Vert_{L_{p,\beta}(\mathbb{R}^{d+1}_+)}\right)\Vert V \Vert_{B^{\frac{1}{2}-\frac{\gamma+1}{2p'},1-\frac{\gamma+1}{p'}}_{p',p'}(\mathbb{R}^{d})},
	\end{equation}
	where $\gamma$ is from \eqref{eq:defgamma}, $\Tilde{\omega}(h)=\omega_{a_0'}(h)+h$, and $N=N(d,p,\beta,\kappa,\Vert a_0' \Vert_{L_{\infty}(\R_+)}) > 0$.
\end{lem}

\begin{proof}
	Using 
	\begin{equation*}
		|\langle \tau_h u, v^h\rangle - \langle u, v\rangle| \leq |\langle \tau_h u - u, v^h\rangle| + \langle u, v^h-v\rangle|,
	\end{equation*}
	and the same inequality for the derivatives $v_{x_i}$, together with Lemma \ref{lem:holderweighted}, give
	\[
		r(h) \leq N_d \Vert \tau_h u - u\Vert_{L_{p,\beta}(\mathbb{R}^{d+1}_+)}\Vert v^h\Vert_{W^{0,1}_{p',\gamma}(\mathbb{R}^{d+1}_+)} 
		+ N_d \Vert u\Vert_{L_{p,\beta}(\mathbb{R}^{d+1}_+)}\Vert v^h-v\Vert_{W^{0,1}_{p',\gamma}(\mathbb{R}^{d+1}_+)}.
	\]
	By Lemma \ref{lem:continuityv^h} and the continuity of shifts $\tau_h$ in $L_{p,\beta}(\mathbb{R}^{d+1}_+)$ with $\beta \geq 0$ (see Lemma \ref{lem:contshiftslp}), we get $\lim_{h \to 0^+} r(h) = 0$. If $a_0' \in BUC(\R_+)$, we get (\ref{eq:estimaterh}) from Lemma \ref{lem:continuityv^h} and (\ref{eq:basicestimatevh}).
\end{proof}

\begin{proof}[Proof of Theorem {\rm \ref{thm:tracepar}}]
Without loss of generality, assume that $a_{dj} = 0$ for $1 \leq j \leq d-1$ and $a_{dd}=1$. The general case follows from the change of variables \eqref{eq:coordchange}. Everywhere below in the proof $N = N(d,\alpha,p,\beta,\kappa) > 0$ unless otherwise stated.

{(i)} We first consider the case $\beta \in (p-1,2p-1)$. Recall that $\gamma \in (-1,p'-1)$ is defined by (\ref{eq:defgamma}).

Suppose that $u \in \mathcal{H}^1_{p,\beta}(\mathcal{P}, \mathbb{R}^{d+1}_+)$. Then, $u$ satisfies
\begin{equation} \label{eq:equationutrace} 
\mathcal{P}u = a_0(x_d)u_t - u_{x_d x_d}-\frac{\alpha}{x_d}u_{x_d} = \sum_{i=1}^{d-1}(F_i)_{x_i}+f
\end{equation}
for some $F_i, f \in L_{p,\beta}(\mathbb{R}^{d+1}_+)$ with
\[
\sum_{i=1}^{d-1}\Vert F_i \Vert_{L_{p,\beta}(\mathbb{R}^{d+1}_+)} + \Vert f \Vert_{L_{p,\beta}(\mathbb{R}^{d+1}_+)} \leq \Vert u \Vert_{\mathcal{H}^1_{p,\beta}(\mathcal{P}, \mathbb{R}^{d+1}_+)}.
\]
Fix any $V \in \mathcal{S}(\mathbb{R}^{d})$ and $h \geq 0$. Let $v^h \in W^{0,1}_{p',\gamma}(\mathbb{R}^{d+1}_+)$ be defined by (\ref{eq:definitionvh}). For $x_d > 0$, define
\begin{gather*}
z^h(x_d) = \langle \tau_h u, v^h \rangle = \langle u(\cdot,x_d+h), v^h(\cdot,x_d) \rangle = \int_{\R^d} u(t,x',x_d+h) v^h(t,x',x_d) dtdx'.
\end{gather*}
Using the Leibnitz rule, we compute the first two derivatives of $z^h$:
\begin{gather*}
(z^h)'(x_d) = \langle \tau_h u_{x_d}, v^h \rangle + \langle \tau_h u, v^h_{x_d} \rangle, \\
(z^h)''(x_d) = \langle  \tau_h u_{x_d x_d}, v^h \rangle + 2 \langle  \tau_h u_{x_d}, v^h_{x_d} \rangle + \langle  \tau_h u, v^h_{x_d x_d} \rangle.
\end{gather*}
From equations (\ref{eq:equationutrace}) and (\ref{eq:definitionvh}), it follows that $z^h$ satisfies the following ODE:
\begin{equation} \label{eq:ODEz^h}
    (z^h)''(x_d) +\frac{\alpha}{h+x_d}(z^h)'(x_d)= g^h(x_d),
\end{equation}
where 
\begin{equation} \label{eq:g^hformula}
\begin{gathered}
g^h(x_d)= \langle \tau_h(u_{x_dx_d}+\alpha x_d^{-1}u_{x_d}),v^h \rangle + \langle \tau_h u, v^h_{x_dx_d}\rangle + 2 \langle  \tau_h u_{x_d}, v^h_{x_d} \rangle 
+ \alpha \langle \tau_h (u/x_d), v^h_{x_d}\rangle \\
= \langle \tau_h(a_0 u_t - \sum_{i=1}^{d-1} (F_i)_{x_i}-f), v^h\rangle + \langle \tau_h u, (\tau_h a_0)v^h_t-\Delta_{x'}v^h+v^h\rangle 
+ 2 \langle  \tau_h u_{x_d}, v^h_{x_d} \rangle \\ + \alpha \langle \tau_h (u/x_d), v^h_{x_d}\rangle = 
\Bigl[\text{using $\langle \tau_h u_t, v^h\rangle + \langle \tau_h u, v^h_t\rangle=0$ and integration by parts in $x'$} \Bigr] \\
= -\langle \tau_h f, v^h \rangle + \sum_{i=1}^{d-1}\langle  \tau_h F_{i}, v^h_{x_i} \rangle + \langle \tau_h u, v^h \rangle + \sum_{i=1}^{d-1}\langle  \tau_h u_{x_i}, v^h_{x_i} \rangle 
+ 2 \langle  \tau_h u_{x_d}, v^h_{x_d} \rangle + \alpha \langle \tau_h (u/x_d), v^h_{x_d}\rangle.
\end{gathered}
\end{equation}

Using Lemma \ref{lem:holderweighted} and Hardy's inequality with $\beta > p-1$, we estimate
\begin{equation} \label{eq:z^h''L1estimate}
\begin{gathered}
\Vert x_d g^h(x_d)\Vert_{L_1(0,\infty)} \leq N \Bigl\{ \sum_{i=1}^{d-1}\Vert \tau_h F_i \Vert_{L_{p,\beta}(\mathbb{R}^{d+1}_+)} + \Vert \tau_h f \Vert_{L_{p,\beta}(\mathbb{R}^{d+1}_+)} + \Vert \tau_h u \Vert_{W^{0,1}_{p,\beta}(\mathbb{R}^{d+1}_+)} \\
+ \Vert \tau_h (u/x_d) \Vert_{L_{p,\beta}(\mathbb{R}^{d+1}_+)}\Bigr\} \Vert v^h \Vert_{W^{0,1}_{p',\gamma}(\mathbb{R}^{d+1}_+)} \leq N \Vert u \Vert_{\mathcal{H}^1_{p,\beta}(\mathcal{P}, \mathbb{R}^{d+1}_+)} \Vert V \Vert_{B^{\frac{1}{2}-\frac{\gamma+1}{2p'},1-\frac{\gamma+1}{p'}}_{p',p'}(\mathbb{R}^{d})}.
\end{gathered}
\end{equation}

For $h=0$, we write $z(s)$ and $g(s)$ instead of $z^0(s)$ and $g^0(s)$. From Lemma \ref{lem:1dtrace} with $h=0$, it follows that for all $V \in \mathcal{S}(\mathbb{R}^{d})$,
\begin{equation*}
\langle U, V \rangle := \lim_{x_d\to 0^+}\langle u(\cdot, x_d), v(\cdot, x_d) \rangle = \frac{1}{1-\alpha}\int_0^{\infty} s (z''(s)+\frac{\alpha}{s}z'(s)) ds.
\end{equation*}
Moreover, by (\ref{eq:z^h''L1estimate})
\begin{equation*}
| \langle U, V \rangle | \leq N \Vert u \Vert_{\mathcal{H}^1_{p,\beta}(\mathcal{P}, \mathbb{R}^{d+1}_+)} \Vert V \Vert_{B^{\frac{1}{2}-\frac{\gamma+1}{2p'},1-\frac{\gamma+1}{p'}}_{p',p'}(\mathbb{R}^{d})}.
\end{equation*}
From Lemma \ref{lem:dualitybesov}, it now follows that
\begin{equation} \label{eq:tracespaceandestimate}
\begin{gathered}
U \in B^{\frac{1}{2}-\frac{\beta+1}{2p},1-\frac{\beta+1}{p}}_{p,p}(\mathbb{R}^{d}) = (B^{\frac{1}{2}-\frac{\gamma+1}{2p'},1-\frac{\gamma+1}{p'}}_{p',p'}(\mathbb{R}^{d}))^*, \\
\Vert U \Vert_{B^{\frac{1}{2}-\frac{\beta+1}{2p},1-\frac{\beta+1}{p}}_{p,p}(\mathbb{R}^{d})} \leq N \Vert u \Vert_{\mathcal{H}^1_{p,\beta}(\mathcal{P}, \mathbb{R}^{d+1}_+)}.
\end{gathered}
\end{equation}
We now show the weak convergence
\begin{equation} \label{eq:weakconvpar}
u(\cdot,h) \overset{\text{wk}}{\rightharpoonup} U \qquad \text{in $B^{\frac{1}{2}-\frac{\beta+1}{2p},1-\frac{\beta+1}{p}}_{p,p}(\mathbb{R}^{d})$ as $h \to 0^+$}.
\end{equation}
Similarly as above for $h=0$, from the ODE (\ref{eq:ODEz^h}) and Lemma \ref{lem:1dtrace}, it follows that for all $h>0$ and $V \in \mathcal{S}(\mathbb{R}^d)$, we have
\begin{equation*}
    \langle u(\cdot,h), V \rangle = z^h(0) = \frac{1}{1-\alpha}\int_0^{\infty}K(s/h)sg^h(s)ds
\end{equation*}
and the uniform estimate
\begin{equation} \label{eq:unifestbesovpar}
    \Vert u(\cdot,h) \Vert_{B^{\frac{1}{2}-\frac{\beta+1}{2p},1-\frac{\beta+1}{p}}_{p,p}(\mathbb{R}^{d})} \leq N \Vert u \Vert_{\mathcal{H}^1_{p,\beta}(\mathcal{P}, \mathbb{R}^{d+1}_+)},
\end{equation}
where $K(s) = 1 - \frac{(1+s)^{\alpha}-1}{s}$. We now estimate
\begin{gather*}
    | \langle u(\cdot,h) - U, V \rangle| = \frac{1}{1-\alpha}\left|\int_0^{\infty}s(K(s/h)g^h(s)-g(s))ds\right| \\
    \leq N_{\alpha} \int_0^\infty sK(s/h)|g^h(s)-g(s)| ds + N_{\alpha}\int_0^\infty |K(s/h)-1|s|g(s)| ds.
\end{gather*}
We denote $\mathcal{W}:=(u, u/x_d, Du, F_i, f)_{i=1}^{d-1}$. To estimate the second term above, by Lemma \ref{lem:holderweighted} and (\ref{eq:g^hformula}), we have
\begin{equation} \label{eq:firsttermestimate}
\begin{gathered}
\int_0^\infty |K(s/h)-1|s|g(s)| ds \leq N \Vert (K(x_d/h)-1)\mathcal{W}\Vert_{L_{p,\beta}(\mathbb{R}^{d+1}_+)} \Vert v \Vert_{W^{0,1}_{p',\gamma}(\mathbb{R}^{d+1}_+)} \\
\leq N \Vert (K(x_d/h)-1)\mathcal{W}\Vert_{L_{p,\beta}(\mathbb{R}^{d+1}_+)} \Vert V \Vert_{B^{\frac{1}{2}-\frac{\gamma+1}{2p'},1-\frac{\gamma+1}{p'}}_{p',p'}(\mathbb{R}^{d})}.
\end{gathered}
\end{equation}
By the dominated convergence theorem,
\begin{equation*}
    \lim_{h \to 0^+} \Vert (K(x_d/h)-1)\mathcal{W}\Vert_{L_{p,\beta}(\mathbb{R}^{d+1}_+)} = 0.
\end{equation*}

To estimate the remaining term, note that from (\ref{eq:g^hformula}), we have
\begin{gather*}
g^h(s)-g(s) = -(\langle \tau_h f, v^h \rangle-\langle f, v \rangle) + \sum_{i=1}^{d-1}(\langle \tau_h F_{i}, v^h_{x_i} \rangle-\langle F_{i}, v_{x_i} \rangle) + (\langle \tau_h u, v^h \rangle- \langle u, v \rangle) \\
+ \sum_{i=1}^{d-1}(\langle \tau_h u_{x_i}, v^h_{x_i} \rangle-\langle u_{x_i}, v_{x_i} \rangle)
+ 2(\langle  \tau_h u_{x_d}, v^h_{x_d} \rangle-\langle  u_{x_d}, v_{x_d} \rangle) + \alpha (\langle \tau_h (u/x_d), v^h_{x_d}\rangle-\langle u/x_d, v_{x_d}\rangle).
\end{gather*}
Lemma \ref{lem:shiftintegralcontinuity}, applied to each term in curly brackets in the above formula, gives
\begin{equation*} 
\begin{gathered}
    \left|\int_0^{\infty}sK(s/h)(g^h(s)-g(s))ds\right| \leq N_{\alpha}\int_0^{\infty}s|g^h(s)-g(s)|ds \to 0 \qquad \text{as $h \to 0^+$}.
\end{gathered}
\end{equation*}
Moreover, when $a_0' \in BUC(\R_+)$, we have
\begin{equation}\label{eq:secondtermestimate}
\begin{gathered}
\int_0^{\infty}s|g^h(s)-g(s)|ds \leq N\left(\Vert (\tau_h-\operatorname{id})\mathcal{W}\Vert_{L_{p,\beta}(\mathbb{R}^{d+1}_+)}+\Tilde{\omega}(h)\Vert \mathcal{W}\Vert_{L_{p,\beta}(\mathbb{R}^{d+1}_+)}\right)\Vert V \Vert_{B^{\frac{1}{2}-\frac{\gamma+1}{2p'},1-\frac{\gamma+1}{p'}}_{p',p'}(\mathbb{R}^{d})},
\end{gathered}
\end{equation}
where $\Tilde{\omega}(h)=\omega_{a_0'}(h)+h$ and $N=N(d,p,\beta,\kappa,\Vert a_0' \Vert_{L_{\infty}(\R_+)}) > 0$.

We have shown that
\begin{equation*}
    | \langle u(\cdot, h) - U, V \rangle| \to 0 \qquad \text{as $h \to 0^+$}
\end{equation*}
for all $V \in \mathcal{S}(\mathbb{R}^d)$. The weak convergence (\ref{eq:weakconvpar}) now follows from the uniform estimate (\ref{eq:unifestbesovpar}) and the denseness of $\mathcal{S}(\mathbb{R}^d)$ in $B^{\frac{1}{2}-\frac{\gamma+1}{2p'},1-\frac{\gamma+1}{p'}}_{p',p'}(\mathbb{R}^{d})$. When $a_0' \in BUC(\R_+)$, the strong convergence
$$\lim_{h \to 0^+} \Vert u(\cdot, h) - U\Vert_{B^{\frac{1}{2}-\frac{\beta+1}{2p},1-\frac{\beta+1}{p}}_{p,p}(\mathbb{R}^{d})}=0$$
follows from (\ref{eq:firsttermestimate}) and (\ref{eq:secondtermestimate}).

{(ii)} We now consider the case $\beta = p-1$ and $\gamma = p'-1$. The argument here is similar to Part (i), but a different test function will be used, due to the absence of Hardy's inequality. Using Theorem \ref{thm:solvpar} with $\gamma = p'-1$, for a fixed $V \in \mathcal{S}(\R^d)$, let $v \in W^{0,1}_{p',\gamma}(\R^{d+1}_+)$ be the solution to
\begin{equation*}
\left\{
\begin{array}{rll}
   a_0(x_d) v_t-\Delta v -\frac{\alpha}{x_d}D_d v+ v & = 0 & \text{in $\mathbb{R}^{d+1}_+$}, \\
    v(t,x',0) & = V(t,x') & \text{on $\mathbb{R}^{d}$}.
\end{array}\right.
\end{equation*}

Again, consider the function $z(x_d) = \langle u(\cdot,x_d), v(\cdot, x_d)\rangle$. It satisfies the ODE
\begin{equation*}
z''(x_d)+\frac{\alpha}{x_d}z'(x_d) = g(x_d) = -\langle f, v \rangle + \sum_{i=1}^{d-1}\langle  F_{i}, v_{x_i} \rangle + \langle u, v \rangle + \sum_{i=1}^{d-1}\langle  u_{x_i}, v_{x_i} \rangle + 2 \langle u_{x_d}, v_{x_d} \rangle.
\end{equation*}
We then again define $\langle U, V \rangle := z(0)$ and, as in Part (i), show (\ref{eq:tracespaceandestimate}).

To finish the proof, fix a cut-off function $\chi \in C_0^{\infty}(\R)$ such that $\chi(x_d) = 1$ for $|x_d| \leq 1$. Then $u\chi \in \mathcal{H}^1_{p,p-1+\varepsilon p}(\mathcal{P}, \mathbb{R}^{d+1}_+)$ for all $\varepsilon > 0$. Since $u=u\chi$ for $x_d \leq 1$, from Part (i) it follows that
\begin{equation*}
u(\cdot,x_d) \overset{\text{wk}}{\rightharpoonup} U \quad \text{weakly in $B^{-\varepsilon/2,-\varepsilon}_{p,p}(\R^d)$ as $x_d \to 0^+$.}
\end{equation*}
If $a_0' \in BUC(\R_+)$, the convergence above is strong. The proof is complete.
\end{proof}

\begin{remark}
In the above theorem for $\beta=p-1$, we do not prove the convergence $u(\cdot,x_d) \to u(\cdot,0)$ in the exact trace space norm $B^{0,0}_{p,p}(\R^d)$ even for $\alpha \leq \frac{p-1}{2p-1}$ as in Theorem \ref{thm:traceell}. The reason is that in the parabolic case, the assumptions in Theorem \ref{thm:tracepar} do not provide sufficient time regularity to apply the trace interpolation method. Assume for simplicity $\alpha=0, a_0=1$, and $a_{dj}=\delta_{dj}$. The given assumption is $u \in \mathcal{H}^1_{p,p-1}(\partial_t-D^2_d,\R^{d+1}_+)$. Equivalently, $u \in L_p(\R_+,x_d^{p-1}dx_d;W^{0,1}_p(\R^d)), D_d u \in L_p(\R_+,x_d^{p-1}dx_d;L_p(\R^d))$, and $u_t-D^2_d u \in L_p(\R_+,x_d^{p-1}dx_d;(1-\Delta_{x'})^{1/2}L_p(\R^d))$. Whereas to apply trace interpolation and recover the space $B^{0,0}_{p,p}(\R^d)$, the conditions $u \in L_{p}(\R_+,x_d^{p-1}dx_d;H^{1/2,1}_p(\R^d))$ and $D^2_d u \in L_p(\R_+,x_d^{p-1}dx_d; H^{-1/2,-1}_p(\R^d))$ are needed.
\end{remark}

\begin{remark}
It should be noted that the condition $a_{ij}=a_{ij}(x_d)$ is not essential and can be generalized to $a_{ij}=a_{ij}(t,x',x_d)$ for some $i,j$ in Theorems \ref{thm:existsolell}, \ref{thm:existsolpar}, \ref{thm:solvnondivpar}, \ref{thm:solvnondivell}, and Propositions \ref{prop:solvabilityellnondeg}, \ref{prop:solvabilityparnondeg}, with the additional condition that $\lambda$ is large enough. Indeed, in Theorems \ref{thm:existsolell} and \ref{thm:existsolpar}, we can take $a_{ij} \in BMO_{t,x'}$ for $1 \leq i, j \leq d-1$, as their proofs rely only on \cite[Theorem 2.4]{DongPhanIndiana2023}. We still have $a_{ij}=a_{ij}(x_d)$ for $\max\{i,j\} = d$, as this was used in the derivation of  equation (\ref{eq:equationw1}). In Theorems \ref{thm:solvnondivpar} and \ref{thm:solvnondivell}, we can take $a_{ij} \in BMO_{t,x'}$ for $1 \leq i \leq d-1$ and $1 \leq j \leq d$, as their proofs rely on respectively Theorems 2.3 and 2.4 from \cite{DongPhan2023JFA}. In Propositions \ref{prop:solvabilityellnondeg} and \ref{prop:solvabilityparnondeg}, we can take $a_{ij} \in BMO_{t,x'}$ for all $i,j=1,\dots,d$, as their proofs rely on respectively Theorems 2.8 and 2.5 from \cite{DongPhan2021TAMS}. Here $a_{ij} \in BMO_{t,x'}$ stands for the partially BMO type condition on $a_{ij}$, i.e., mere measurability in the $x_d$ variable and smallness of certain local $BMO$-seminorms of $a_{ij}$ with respect to $t$ and $x'$. See the aforementioned papers for precise definitions.

In contrast, the condition $a_{ij}=a_{ij}(x_d)$ is significant in Theorems \ref{thm:solvell} and \ref{thm:solvpar}, since differentiation in $t$ and $x'$ was taken in their proofs. It is of interest to investigate under what conditions on $a_{ij}$ Theorems \ref{thm:solvell} and \ref{thm:solvpar} can be extended to the case $a_{ij}=a_{ij}(t,x',x_d)$. 
\end{remark}

\appendix
\section{} \label{appendix}

\begin{defn}
Let $\mu$ be a non-negative Borel measure on $\R^d_+$ such that $0 < \mu(B_r^+(x)) < \infty$ for all $x \in \overline{\R^d_+}$ and $r>0$. For each $p \in (1,\infty)$, a measurable function $w$ such that $0 < w(x) < \infty$ for a.e. $x \in \R^d_+$ is said to be in the $A_p(\R^d_+,\mu)$ Muckenhoupt class of weights if $[w]_{A_p(\R^d_+,\mu)} < \infty$, where
\begin{equation*}
[w]_{A_p(\R^d_+,\mu)} = \sup_{r>0, x_0 \in \overline{\R^d_+}} \left[ \fint_{B_r^+(x_0)} w(x) d\mu\right]\left[ \fint_{B_r^+(x_0)} w(x)^{-\frac{1}{p-1}} d\mu\right]^{p-1}.
\end{equation*}
\end{defn}
Here $B_r^+(x_0)=\{x \in \R^d: |x-x_0|<r \} \cap \R^d_+$.
\begin{lem} \label{lem:finiteintegral}
    Let $a,b \in \R$. Then
    \begin{equation*}
        \int_{\R^d_+} (1+|x|^2)^{\frac{a}{2}}x_d^{b} dx < +\infty
    \end{equation*}
    if and only if
    \begin{gather*}
        b > -1 \quad \text{and} \quad a + b < -d.
    \end{gather*}
\end{lem}
\begin{proof} Denote $w(x) = (1+|x|^2)^{\frac{a}{2}}x_d^{b}$. Split $\int_{\R^d_+} w = \int_{B_1^+} w + \int_{\R^d_+ \setminus B_1^+} w$. The first term $\int_{B_1^+} w$ is finite if and only if $b > -1$. For $|x| \geq 1$, we have $1+|x|^2 \sim |x|^2$. The second term
\begin{gather*}
\int_{\R^d_+ \setminus B_1^+} w \sim \int_{\R^d_+ \setminus B_1^+} |x|^a x_d^b dx = \sum_{i=0}^{\infty} \int_{B_{2^{i+1}}^+\setminus B_{2^i}^+} |x|^a x_d^b dx = \sum_{i=0}^{\infty} 2^{i(a+b+d)}\int_{B_{2}^+\setminus B_{1}^+} |x|^a x_d^b dx
\end{gather*}
is finite if and only if $b > -1$ and $a+b+d<0$.
\end{proof}
\begin{lem} \label{lem:A_pparameters}
    Let $\gamma_0 > -1, p \in (1,\infty), a,b \in \R$. Then
    \begin{equation*}
        w(x) = w_{a,b}(x) = (1+|x|^2)^{\frac{a}{2}}x_d^{b} \in A_p(\R^d_+, x_d^{\gamma_0}dx)
    \end{equation*}
    if and only if
    \begin{equation} \label{eq:A_pparameters}
    \begin{gathered}
        -(1+\gamma_0) < b < (1+\gamma_0)(p-1), \\
        -(d+\gamma_0) < a + b < (d+\gamma_0)(p-1).
    \end{gathered}
    \end{equation}
\end{lem}
\begin{proof} Let $d\mu = x_d^{\gamma_0}dx$. Denote
\[
A(x_0,r) = A_{a,b}(x_0,r) = \left[ \fint_{B_r^+(x_0)} w_{a,b}(x) d\mu\right]\left[ \fint_{B_r^+(x_0)} w_{a,b}(x)^{-\frac{1}{p-1}}  d\mu\right]^{p-1}
\]
and $A(r) = A(0,r)$. By definition of $A_p$ weights, $w_{a,b} \in A_p(\R^d_+, x_d^{\gamma_0}dx)$ if and only if 
$$\sup\limits_{r>0, x_0 \in \overline{\R^d_+}}A_{a,b}(x_0,r) < \infty.$$

First consider the case $a=0$. From scaling, it follows that for $x_0=0$,
\begin{equation*}
A_{0,b}(r) = \left[ \int_{B_r^+} x_d^b x_d^{\gamma_0} dx\right]\left[ \int_{B_r^+} x_d^{-\frac{b}{p-1}} x_d^{\gamma_0} dx\right]^{p-1} \left[\int_{B_r^+} x_d^{\gamma_0} dx\right]^{-p}
\end{equation*}
does not depend on $r>0$ and is finite if and only if $b+\gamma_0 > -1$ and $-b/(p-1)+\gamma_0>-1$, which is equivalent to the first condition in (\ref{eq:A_pparameters}). For general $x_0$, we can reduce to $x_0=0$. Indeed, since both $w_{0,b}(x)$ and $d\mu$ do not depend on $x'$, we may assume $x_0'=0$. If $x_{0d} > 2r$, then $w_{0,b}(x) \sim x_{0d}^b$ in $B_r^+(x_0)$ and therefore, $A_{0,b}(x_0,r) \lesssim 1$. If $x_{0d} \leq 2r$, $A_{0,b}(x_0,r) \lesssim A_{0,b}(5r)$ by the doubling property of $d\mu$.

For general $a$, we now show that $w_{a,b} \in A_p(\R^d_+, x_d^{\gamma_0}dx)$ if and only if $\sup_{r>0}A_{a,b}(r) < \infty$ and $w_{0,b} \in A_p(\R^d_+, x_d^{\gamma_0}dx)$. Indeed,

if $x_{0d} > 2r$, then $w_{a,b}(x) \sim (1+|x_0|^2)^{\frac{a}{2}}x_{0d}^b$ in $B_r^+(x_0)$ and therefore $A(x_0,r) \lesssim 1$.

If $|x_0'| > 2r$, then $w_{a,b}(x) \sim (1+|x_0|^2)^{\frac{a}{2}}x_{d}^b$ in $B_r^+(x_0)$. Hence $A_{a,b}(x_0,r) \sim A_{0,b}(x_0,r)$ is bounded for $|x_0'|>2r$ if and only if $w_{0,b} \in A_p(\R^d_+, x_d^{\gamma_0}dx)$.

If $|x_0'| \leq 2r, x_{0d} \leq 2r$, from the doubling property of $d\mu$ it follows that $A_{a,b}(x_0,r) \lesssim A_{a,b}(5r)$.

We now show that $\sup_{r>0}A_{a,b}(r) < \infty$ is equivalent to (\ref{eq:A_pparameters}). For $r \leq 1$, we have $A_{a,b}(r) \sim A_{0,b}(r)$ is bounded if and only if $-(1+\gamma_0) < b < (1+\gamma_0)(p-1)$. For $r > 1$, by the doubling property of $d\mu$, we may without loss assume $r=2^n$, where $n \in \mathbb{N}$. Then by splitting the integrals as $\int_{B_{2^n}^+} = \int_{B_{1}^+} + \sum_{i=1}^{n-1} \int_{B_{2^{i}}^+\setminus B_{2^{i-1}}^+}$ and using scaling, we get
\begin{gather*}
A_{a,b}(2^n) \sim \sum_{i=0}^{n-1} 2^{i(a+b+\gamma_0+d)} \left(\sum_{j=0}^{n-1}2^{j(-\frac{a+b}{p-1}+\gamma_0+d)}\right)^{p-1} 2^{-np(\gamma_0+d)}.
\end{gather*}
This expression is uniformly bounded as $n \to \infty$ if and only if $a+b+\gamma_0+d>0$ and $-(a+b)/(p-1)+\gamma_0+d>0$, which is equivalent to the second condition in (\ref{eq:A_pparameters}).
\end{proof}
The statement of Lemma \ref{lem:A_pparameters} for $a=0$ was previously used in \cite[Remark 2.5]{DongPhan2023JFA}. See also \cite{MiaoZhaoArxiv2023} for criteria for anisotropic weights of the form $|x'|^{\theta_1}|x|^{\theta_2}|x_d|^{\theta_3}$ to be of the class $A_p(\R^d)$.

To prove Lemma \ref{lem:solutionoperatorwelldefined}, we first prove
\begin{lem} \label{lem:bddrhs}
Let $ \alpha \in (-\infty, 1), p_i \in (1,\infty)$, and $\beta_i \in (\alpha^+ p_i-1, 2p_i-1)$, where $i=1,2$. Suppose that $f \in L_{\infty}(\R^{d+1}_+)$ has compact support in $\{x_d>0\}$. Let $u_i \in W^{1,2}_{p_i,\alpha,\beta_i}(\R^{d+1}_+)$, $i=1,2$, be the solutions to \eqref{eq:nondivpar} with $u(t,x',0)=0$ given by Theorem {\rm \ref{thm:solvnondivpar}}. Then $u_1=u_2$.
\end{lem}

\begin{proof}
Let $\gamma_0 = 1- \alpha$. Fix any $q \in (\max(p_1,p_2),\infty)$. Then, fix any $a,b \in \R$ such that
\begin{gather*}
-(1+\gamma_0) < b < q\min_{i=1,2}\frac{\beta_i-\alpha p_i +1}{p_i}-(1+\gamma_0), \\
q\max_{i=1,2}\frac{\beta_i-\alpha p_i +d}{p_i}-(d+\gamma_0) < a + b < (q-1)(d+\gamma_0).
\end{gather*}
Then, by using Lemmas \ref{lem:finiteintegral}, \ref{lem:A_pparameters}, and H\"older's inequality, one may verify that the weight
\begin{equation*}
w(x) = (1+|x|^2)^{\frac{a}{2}}x_d^b \in A_q(\R^d_+, x_d^{\gamma_0}dx)
\end{equation*}
and that
\[
L_q(\R^{d+1}_+,x_d^{\alpha q+\gamma_0}w(x) dxdt) \subset L_{p_i,\beta_i}(\R^{d+1}_+), \qquad i=1,2.
\]
From \cite[Theorem 2.3]{DongPhan2023JFA}, it follows that there exists a unique solution
\[
u \in W^{1,2}_{q,\alpha}(\R^{d+1}_+,x_d^{\alpha q+\gamma_0}w(x) dxdt) \subset W^{1,2}_{p_i,\alpha,\beta_i}(\R^{d+1}_+), \qquad i=1,2,
\]
to (\ref{eq:nondivpar}) with $u(t,x',0)=0$. Here $W^{1,2}_{q,\alpha}(\R^{d+1}_+,x_d^{\alpha q+\gamma_0}w(x) dxdt)$ is defined similarly to $W^{1,2}_{q,\alpha,\beta}(\R^{d+1}_+)$, except that the measure $x_d^{\alpha q+\gamma_0}w(x) dxdt$ is used instead of $x_d^{\beta}dxdt$. From Theorem \ref{thm:solvnondivpar}, it follows that $u_1=u$ and $u_2=u$ and therefore, $u_1=u_2$.
\end{proof}

\begin{proof}[Proof of Lemma {\rm \ref{lem:solutionoperatorwelldefined}}]
    Let $E_{\text{par}}$ be the extension operator from Proposition \ref{prop:extensionoperator}. Note that $v_i:= u_i - E_{\text{par}}U \in W^{1,2}_{p_i,\alpha,\beta_i}(\R^{d+1}_+), \; i=1,2$, satisfy (\ref{eq:nondivpar}) with the zero boundary condition for some
    \[
    f \in L_{p_1,\beta_1}(\R^{d+1}_+) \cap L_{p_2,\beta_2}(\R^{d+1}_+).
    \]
    Taking a sequence of $f_n \in L_{\infty}(\R^{d+1}_+)$ with compact supports in $\{x_d > 0\}$ such that
    \[
    f_n \to f \quad \text{in} \quad L_{p_1,\beta_1}(\R^{d+1}_+) \cap L_{p_2,\beta_2}(\R^{d+1}_+),
    \]
    and using Lemma \ref{lem:bddrhs}, give $v_1=v_2$. Hence $u_1=u_2$.
\end{proof}

\bibliographystyle{amsplain}

\begin{thebibliography}{10}

\bibitem{amann2019linear}
H.~Amann, \emph{Linear and quasilinear parabolic problems: Volume {II}:
  Function spaces}, Monographs in Mathematics, Springer International
  Publishing, 2019.

\bibitem{audrito2024higher}
Alessandro Audrito, Gabriele Fioravanti, and Stefano Vita, \emph{Higher order
  {S}chauder estimates for parabolic equations with degenerate or singular
  weights}, arXiv:2403.08575.

\bibitem{audrito2024schauder}
Alessandro Audrito, Gabriele Fioravanti, and Stefano Vita, \emph{Schauder
  estimates for parabolic equations with degenerate or singular weights}, Calc.
  Var. Partial Differential Equations \textbf{63} (2024), no.~8, 204.

\bibitem{BartonMayboroda}
A.~Barton and S.~Mayboroda, \emph{Layer potentials and boundary-value problems
  for second order elliptic operators with data in {B}esov spaces}, American
  Mathematical Society, Providence, Rhode Island, 2016.

\bibitem{BogachevSmolyanov}
Vladimir~I. Bogachev and Oleg~G. Smolyanov, \emph{Real and functional
  analysis}, Moscow Lectures, vol.~4, Springer, Cham, [2020] \copyright 2020,
  Expanded and revised version of the Russian original. \MR{4292288}

\bibitem{CaffarelliSilvestre2007CPDE}
Luis Caffarelli and Luis Silvestre, \emph{An extension problem related to the
  fractional {L}aplacian}, Comm. Partial Differential Equations \textbf{32}
  (2007), no.~8, 1245--1260.

\bibitem{CwikelEinavJFA2019}
Michael Cwikel and Amit Einav, \emph{Interpolation of weighted {S}obolev
  spaces}, J. Funct. Anal. \textbf{277} (2019), no.~7, 2381--2441.

\bibitem{DongJeonVita2023arxiv}
Hongjie Dong, Seongmin Jeon, and Stefano Vita, \emph{Schauder type estimates
  for degenerate or singular elliptic equations with {D}ini mean oscillation
  coefficients with application}, arXiv:2311.06846, to appear in Calc. Var.
  Partial Differential Equations.

\bibitem{DongPhan2021TAMS}
Hongjie Dong and Tuoc Phan, \emph{Parabolic and elliptic equations with
  singular or degenerate coefficients: the {D}irichlet problem}, Trans. Amer.
  Math. Soc. \textbf{374} (2021), no.~9, 6611--6647. \MR{4302171}

\bibitem{DongPhanIndiana2023}
Hongjie Dong and Tuoc Phan, \emph{On parabolic and elliptic equations with
  singular or degenerate coefficients}, Indiana Univ. Math. J. \textbf{72}
  (2023), 1461--1502.

\bibitem{DongPhan2023JFA}
\bysame, \emph{Weighted mixed-norm {$L_p$} estimates for equations in
  non-divergence form with singular coefficients: The {D}irichlet problem}, J.
  Funct. Anal. \textbf{285} (2023), no.~2, 109964.

\bibitem{Grisvard}
P.~Grisvard, \emph{Espaces interm\'ediaires entre espaces de {Sobolev} avec
  poids}, Annali della Scuola Normale Superiore di Pisa - Scienze Fisiche e
  Matematiche \textbf{3e s{\'e}rie, 17} (1963), no.~3, 255--296 (fr).
  \MR{160104}

\bibitem{Grubb2018JFA}
Gerd Grubb, \emph{Regularity in {$L_p$} {S}obolev spaces of solutions to
  fractional heat equations}, J. Funct. Anal. \textbf{274} (2018), no.~9,
  2634--2660.

\bibitem{JonesJMM1968}
B.~Frank Jones, \emph{Lipschitz spaces and the heat equation}, J. Math. Mech.
  \textbf{18} (1968), no.~5, 379--409.

\bibitem{JungKim2024arxiv}
Pilgyu Jung and Doyoon Kim, \emph{${L}_p$-estimates for parabolic equations in
  divergence form with a half-time derivative}, arXiv:2407.11305.

\bibitem{KimKHKimWoo2024}
Doyoon Kim, Kyeong-Hun Kim, and Kwan Woo, \emph{Trace theorem and non-zero
  boundary value problem for parabolic equations in weighted {S}obolev spaces},
  Stochastics and Partial Differential Equations: Analysis and Computations
  \textbf{12} (2024), no.~1, 134--172.

\bibitem{KozlovNazarovMathNach2009}
Vladimir Kozlov and Alexander Nazarov, \emph{The {D}irichlet problem for
  non-divergence parabolic equations with discontinuous in time coefficients},
  Math. Nachr. \textbf{282} (2009), no.~9, 1220--1241. \MR{2561181}

\bibitem{KrylovCPDE1999}
N.~V. Krylov, \emph{Weighted {S}obolev spaces and {L}aplace's equation and the
  heat equations in a half space}, Comm. Partial Differential Equations
  \textbf{24} (1999), no.~9-10, 1611--1653. \MR{1708104}

\bibitem{LLRVarxiv2024}
Nick Lindemulder, Emiel Lorist, Floris Roodenburg, and Mark Veraar,
  \emph{Functional calculus on weighted {S}obolev spaces for the {L}aplacian on
  the half-space}, arXiv:2406.03297.

\bibitem{LindemulderVeraar2020JDE}
Nick Lindemulder and Mark Veraar, \emph{The heat equation with rough boundary
  conditions and holomorphic functional calculus}, J. Differential Equations
  \textbf{269} (2020), no.~7, 5832--5899.

\bibitem{metafune2023JDE}
G.~Metafune, L.~Negro, and C.~Spina, \emph{A unified approach to degenerate
  problems in the half-space}, J. Differential Equations \textbf{351} (2023),
  63--99.

\bibitem{metafune2024regularity}
Giorgio Metafune, Luigi Negro, and Chiara Spina, \emph{Regularity theory for
  parabolic operators in the half-space with boundary degeneracy},
  arXiv:2309.14319.

\bibitem{MetafuneNegroSpina2023TokyoJMath}
Giorgio Metafune, Luigi Negro, and Chiara Spina, \emph{Anisotropic {S}obolev
  spaces with weights}, Tokyo J. Math. \textbf{46} (2023), no.~2, 313--337.
  \MR{4690601}

\bibitem{MiaoZhaoArxiv2023}
Changxing Miao and Zhiwen Zhao, \emph{On a class of anisotropic {M}uckenhoupt
  weights and their applications to $p$-{L}aplace equations}, arXiv:2310.01359.

\bibitem{SireTerraciniVitaCPDE2021}
Yannick Sire, Susanna Terracini, and Stefano Vita, \emph{Liouville type
  theorems and regularity of solutions to degenerate or singular problems part
  {I}: even solutions}, Comm. Partial Differential Equations \textbf{46}
  (2021), no.~2, 310--361.

\bibitem{SireTerraciniVitaMathEng2021}
\bysame, \emph{Liouville type theorems and regularity of solutions to
  degenerate or singular problems part {II}: odd solutions}, Math. Eng.
  \textbf{3} (2021), no.~1, Paper No. 5, 50. \MR{4144100}

\bibitem{StingaTorrea2017SIAM}
Pablo~Ra\'{u}l Stinga and Jos\'{e}~L. Torrea, \emph{Regularity theory and
  extension problem for fractional nonlocal parabolic equations and the master
  equation}, SIAM J. Math. Anal. \textbf{49} (2017), no.~5, 3893--3924.

\bibitem{Taibleson1964Indiana}
Mitchell~H. Taibleson, \emph{On the theory of {L}ipschitz spaces of
  distributions on {E}uclidean n-space: I. principal properties}, J. Math.
  Mech. \textbf{13} (1964), no.~3, 407--479.

\bibitem{Triebel3}
H.~Triebel, \emph{Theory of function spaces {III}}, Monographs in Mathematics,
  Birkh{\"a}user Basel, 2006.

\bibitem{Triebel}
Hans Triebel, \emph{Interpolation theory, function spaces, differential
  operators}, North-Holland Pub. Co Amsterdam; New York, 1978 (English).

\end{thebibliography}
\providecommand{\bysame}{\leavevmode\hbox to3em{\hrulefill}\thinspace}
\providecommand{\MR}{\relax\ifhmode\unskip\space\fi MR }
\providecommand{\MRhref}[2]{%
  \href{http://www.ams.org/mathscinet-getitem?mr=#1}{#2}
}
\providecommand{\href}[2]{#2}

\end{document}